\documentclass{memo-l}

\usepackage{}
\usepackage{comment}
\usepackage{color}
\DeclareMathOperator*{\esssup}{ess\,sup}

\newtheorem{assumption}{Assumption}

\newtheorem{theorem}{Theorem}[chapter]
\newtheorem{lemma}[theorem]{Lemma}

\theoremstyle{definition}
\newtheorem{definition}[theorem]{Definition}

\theoremstyle{remark}
\newtheorem{remark}[theorem]{Remark}

\numberwithin{section}{chapter}
\numberwithin{equation}{chapter}

\makeindex

\begin{document}

\frontmatter

\title{Utility maximization in constrained and unbounded financial markets:\\[0.2em]\smaller{}Applications to indifference valuation, regime switching, consumption and Epstein-Zin recursive utility\thanks{The work was presented
		at the Workshop on BSDEs, SPDEs and their Applications, Edinburgh,
		U.K., and the seminars at Fudan University, Oxford University, and
		Universit{\'e} Rennes 1. The authors would like to thank the
		audience for helpful comments and suggestions.} }


\author{Ying Hu}
\address{Universit{\'e} Rennes, CNRS, IRMAR-UMR6625, F-35000, Rennes, France.}
\curraddr{}
\email{ying.hu@univ-rennes1.fr}
\thanks{Partially supported by Lebesgue
	Center of Mathematics ``Investissements d'avenir"
	program-ANR-11-LABX-0020-01, by ANR CAESARS (Grant No. 15-CE05-0024)
	and by ANR MFG (Grant No. 16-CE40-0015-01). }

\author{Gechun Liang}
\address{Department of Statistics, University of Warwick, Coventry, CV4 7AL,
	U.K.}
\curraddr{}
\email{g.liang@warwick.ac.uk}
\thanks{Partially supported by Royal Society International Exchange (Grant No. 170137).}

\author{Shanjian Tang}
\address{Department of Finance and Control Sciences, School of Mathematical
	Sciences, Fudan University, Shanghai 200433, China. }
\curraddr{}
\email{sjtang@fudan.edu.cn}
\thanks{Partially
	supported by National Science Foundation of China (Grant No.
	11631004) and Science and Technology Commission of Shanghai
	Municipality (Grant No. 14XD1400400). }

\date{}




\maketitle

\tableofcontents

%

\chapter*{Abstract}

This memoir presents a systematic study of the utility maximization problem of an investor in a constrained and unbounded financial market. Building upon the work of Hu et al. (2005) [{\it Ann. Appl. Probab.}, {\bf 15}, 1691--1712] in a bounded framework, we extend our analysis to the more challenging unbounded case. Our methodology combines both methods of quadratic backward stochastic differential equations with unbounded solutions and convex duality. Central to our approach is the verification of the finite entropy condition, which plays a pivotal role in solving the underlying utility maximization problem and establishing the martingale property and the convex duality representation of the conditional value process. Through four distinct applications, we first study the utility indifference valuation of financial derivatives with unbounded payoffs, uncovering novel asymptotic behaviors as the risk aversion parameter approaches zero or infinity. Furthermore, we study the regime switching market model with unbounded random endowments and consumption-investment problems with unbounded random endowments, both constrained to portfolios chosen from a convex and closed set. Finally, we investigate the investment-consumption problem of an investor with the Epstein-Zin recursive utility in an unbounded financial market.\\

\noindent\textit{Keywords}:
    Quadratic BSDE,  convex duality, exponential utility maximization,  utility indifference valuation, finite entropy condition, regime switching, investment-consumption, Epstein-Zin recursive utility.\\

\noindent\textit{Mathematics Subject Classification (2020)}: 91G10, 91G80, 60H10, 60H30.



%
%
%

\chapter{Introduction}

This memoir contributes to the systematic study of the utility maximization problem for an investor in a \textit{constrained} and \textit{unbounded} financial market. Specifically, the investor's  portfolios are constrained within a convex and closed set. This addresses portfolio constraints more effectively, reflecting real-world scenarios. Additionally, the investor may receive a random endowment, which is not necessarily bounded and thus 
includes the typical European call options. 
Moreover, the market coefficients, such as the market price of risk, may also be unbounded.

The financial market is incomplete due to the inability to fully hedge risks arising from portfolio constraints, random market coefficients, and exposures to non-traded assets. This incompleteness has received great attention in the context of utility maximization problems.
One notable application is in the pricing and hedging of derivatives in incomplete markets through {utility indifference valuation}. Here, the random endowment is treated as the payoff of the derivative.

The concept of {utility indifference valuation} was first introduced by Hodges and Neuberger \cite{HodgesNeuberger} and has since been further developed by researchers such as Becherer \cite{Bec0}, Davis \cite{Davis}, Henderson \cite{Henderson1}, Musiela and Zariphopoulou \cite{Musiela}, and more recently by Biagini et al.~\cite{Biagini}, Frei et al.~\cite{FMS, FS}, Henderson and Liang \cite{Henderson3}, Owen and Zitkovic \cite{Gordan}, among others. For a comprehensive overview of utility indifference valuation and related topics, the reader is referred  to the monograph edited by Carmona \cite{Carmona} and the additional references provided therein.

The unboundedness and constraints inherent in financial models capture market realities, yet they also present mathematical challenges. For instance, traditional convex duality methods have limitations due to the portfolio constraints, as the control domains are not necessarily  the entire space. Moreover, since the market is not necessarily Markovian, with market coefficients and random endowments being general stochastic processes, the standard Hamilton-Jacobi-Bellman (HJB) equation approach does not work.
These challenges call for new mathematical techniques that can address the complexities arising from both unboundedness and constraints in utility maximization modeling. 

{Nobel laureate Harry Markowitz commented in the foreword to Sethi's monograph~\cite{Sethi}: ``{\it It remains to be seen whether the introduction of realistic investor constraints is an impenetrable barrier to analysis, or a golden opportunity for someone with a novel approach.}"  This memoir provides a rigorous mathematical analysis  on how to overcome the barriers arising from several popular  realistic modeling factors (like constraints, unboundedness and regime-switching of financial markets, and consumption and Epstein-Zin recursive utility of investors) by introducing techniques from quadratic backward stochastic differential equations (quadratic BSDEs) with unbounded solutions, alongside the convex duality methods and the verification of the finite entropy condition.}

When the random endowment and market coefficients are bounded,  which defines a \textit{bounded financial market}, addressing portfolio constraints in the corresponding utility maximization problem has been completely solved in Hu et al.~\cite{him}. The solution heavily relies on 
the application of the martingale optimality principle and quadratic BSDEs with \textit{bounded} terminal data. Remarkably, both the value function and the associated optimal trading strategy can be fully characterized in terms of the bounded solution to a quadratic BSDE with bounded terminal data. The solution provides valuable insights into managing portfolio constraints within bounded financial markets. 

It has been twenty years since the work \cite{him} was published, which has since inspired numerous follow-up research utilizing the quadratic BSDE method. For instance, Mania and Schweizer \cite{MS} and Morlais \cite{Morlais} generalized the approach to consider general continuous semimartingale models, while Becherer \cite{Bec} extended it to markets with jumps. Cheridito and Hu \cite{CH} expanded the framework to include both consumption and investment, paving the way to applications such as solution of the optimal  consumption-investment problem  with Epstein-Zin utility (see Xing \cite{Xing2017}). Most existing studies on utility maximization problems with the quadratic BSDE method primarily focus on bounded financial markets, with \cite{Xing2017} being a notable exception.
For example, all of Becherer \cite{Bec}, Cheridito and Hu \cite{CH}, El Karoui and Rouge \cite{El}, Henderson and Liang \cite{Henderson3}, Mania and Schweizer \cite{MS}, and Morlais \cite{Morlais}, among others, deal with the bounded case. These results heavily rely on the BMO martingale property of the solution to quadratic BSDEs with bounded terminal data.

Regarding the theory of quadratic BSDEs, Kobylanski \cite{Kobylanski} first established the existence and uniqueness of \textit{bounded} solutions in a Brownian setting. Subsequently, Briand and Hu \cite{BH1, BH2} solved the case of \textit{unbounded} solutions, further developed in Delbaen et al.~\cite{Delbaen0, dhr} with an additional convexity assumption on generators to strengthen the uniqueness of the solution. In the semimartingale case, Morlais \cite{Morlais} and Tevzadze \cite{Tevzadze} investigated bounded solutions, with the former  extending the main results of \cite{him} and \cite{Kobylanski}, and the latter introducing a fixed-point argument. Additionally, Mocha and Westray \cite{Mocha} extended the semimartingale analysis to unbounded solutions, while El Karoui et al.~\cite{Matoussi2016} and Matoussi and Salhi \cite{Matoussi2019} extended it to a class of exponential quadratic BSDEs with jumps. Briand and Elie \cite{Briand} discussed delayed generators, while useful convexity bounds for quadratic BSDE solutions were derived by Frei et al.~\cite{FMS}. Moreover, Barrieu and El Karoui \cite{Barrieu} introduced the notion of quadratic semimartingales to study the stability of solutions, and Delbaen et al.~\cite{Delbaen0} explored the case with superquadratic generators. 
For more recent developments on quadratic BSDEs and beyond, see Fan et al.~\cite{Fan} for a comprehensive introduction. Additionally, the interested reader may find further technical details in Chapter 7 of the recent monograph by Zhang \cite{ZhangJF}.

Despite the growing interest in unbounded random endowments and market coefficients, relatively few results are available. The well known \textit{one-dimensional case} stands as an exception. In a Markovian framework with a derivative written on a single non-traded asset, Henderson \cite{Henderson1} and Musiela and Zariphopoulou \cite{Musiela} utilized the Cole-Hopf transformation to linearize the HJB equations for the value functions, enabling them to deal with unbounded payoffs in a special case. On the other hand, Xing \cite{Xing2017} addressed the general case of unbounded market coefficients, and tackled the optimal consumption-investment problem with Epstein-Zin utility. This work marks one of the first attempts at utilizing quadratic BSDEs with \textit{unbounded} solutions, as developed in \cite{BH1, BH2}, for utility maximization.

On the other hand, when dealing with \textit{subspace} portfolio constraints, the minimal entropy representation as the convex dual has been extensively utilized to investigate the exponential utility maximization problem in incomplete markets. The works of 
Karatzas et al.~\cite{KS, Kaz} and Kramkov and Schachermayer \cite{KS1999, KS2003} are notable in this regard. Their application to utility indifference valuation is available in Becherer \cite{Bec0}, Delbaen et al.~\cite{dgr}, and Frei and Schweizer \cite{FS}, among others. In these works, the unbounded random endowment can be eliminated through a straightforward change of probability measure. However, this trick does not prove as effective in the Brownian setting. Additionally, the presence of non-subspace portfolio constraints, such as convex and closed constraints, prohibits the application of convex duality methods. This is because the duality between the minimal entropy representation and the exponential utility maximization problem no longer holds.\\

In the rest of this chapter, we briefly sketch the contributions of each following chapter.

In \textbf{Chapter \ref{section:main}}, we build upon and extend the theory of quadratic BSDEs with unbounded terminal data as developed in \cite{BH1, BH2, dhr}, combining it with the convex duality method to solve a general exponential utility maximization problem with unbounded random endowments under general portfolio constraints.
While the existence of a solution to the corresponding quadratic BSDE follows from the localization argument of \cite{BH1}, one  main difficulty lies in the verification of the martingale property of the conditional value process when implementing the optimal trading strategy. If the payoff is unbounded, the BMO martingale property of the solution is lost. To overcome this difficulty, we explore the duality between the quadratic BSDE with unbounded terminal data and the associated optimal density process. Furthermore, we verify the \textit{finite entropy condition} of the optimal density process, which guarantees the martingale property of the conditional value process. As a result, we only require the payoffs to be exponentially integrable, which covers call option payoffs as random endowments.

In \cite{him}, a technical Class (D) condition on the exponential utility of the investor's wealth is imposed, which might appear unnecessary. We remove such a technical assumption when the payoff is bounded from below and exponentially integrable. Instead of the Class (D) condition, we impose an \textit{equivalent minimal martingale measure condition}, which seems to fit better with the pricing and hedging of derivatives. The idea is to approximate the original payoff from below by a sequence of bounded payoffs and construct a sequence of trading strategies that satisfy the Class (D) condition as approximate trading strategies.

The idea of verifying the finite entropy condition is central in our work, and we further develop it in four applications in the next chapters. In \textbf{Chapter \ref{Chapter5}}, we apply to utility indifference valuation of unbounded payoffs. The quadratic BSDE and convex duality methods provide a new \textit{convex dual representation} of the utility indifference price, which seem to generalize the existing minimal entropy representation. This result is motivated by proof of the uniqueness of the solution to quadratic BSDE with unbounded terminal data, as first established in \cite{dhr,dhr2}. When establishing the convex dual representation of the utility indifference price, we focus on the duality between the dual BSDE and the corresponding optimal density process. We verify the finite entropy condition of the optimal density process, which in turn guarantees the existence of the dual optimizer.

The quadratic BSDE method also facilitates a systematic study of the asymptotic analysis of the utility indifference price concerning the risk aversion parameter, particularly when the portfolio is constrained to a convex and closed \textit{cone}. However, the more general case of convex and closed set constraints remains open. To the best of our knowledge, most existing asymptotic results deal with bounded payoffs under \textit{subspace} portfolio constraints (see \cite{dgr,El,MS}). The subspace portfolio constraints allow for the use of the minimal entropy representation, which does not hold under cone portfolio constraints. Instead, we address the primal problem by characterizing the utility indifference price as the unbounded solution to a quadratic BSDE. Utilizing the stability theory of BSDEs, we demonstrate that the utility indifference price converges, on one hand, to an expected payoff under some equivalent probability measure (not necessarily the minimal entropy martingale measure) as the risk aversion parameter tends to zero. On the other hand, it converges to the super-replication price (not necessarily under martingale measures) as the risk aversion parameter tends to infinity. Our asymptotic results appear to cover existing ones, and our method seems new.

\textbf{Chapter \ref{Chapter6}} is concerned with the application to regime switching markets, where the market coefficients follow different stochastic dynamics under different regimes. Regime switching models have demonstrated a closer reflection of reality in practice. For an introduction of regime switching diffusions and their applications, we refer to the monograph by Yin and Zhu \cite{Yin}. In a special Markovian setting with convex portfolio constraints, Heunis \cite{Heunis} tackled a problem similar to our setup. Recently, Hu et al.~\cite{HSX2022} solved the non-Markovian case with \textit{bounded} random endowments and market coefficients using multidimensional quadratic BSDEs with bounded terminal data. Additionally, Hu et al.~\cite{HLT} presented an infinite horizon version related to forward utility maximization. However, no work on \textit{unbounded} random endowments is available for utility maximization in regime-switching market models to date.

When dealing with unbounded random endowments, the original characterization of quadratic BSDE with unbounded terminal data becomes a system of equations, with each equation representing one market regime. Due to the regime switching feature, the generator of the BSDE system satisfies the monotonic condition proposed by Hu and Peng \cite{HuPeng}, enabling the application of the comparison theorem for multidimensional BSDEs. This comparison theorem is crucial for solution existence, as it allows for the determination of common upper and lower bounds for each equation. These common bounds are obtained by solving simple \textit{one-dimensional} BSDEs. Subsequently, the localization technique, similar to those employed in the one-dimensional case, can be applied to handle the unbounded terminal data. Furthermore, we observe that the corresponding multidimensional quadratic BSDE driven by Brownian motion has an equivalent formulation of a one-dimensional quadratic BSDE driven by Brownian motion and Poisson processes used to model market regimes. In the literature, the latter equation is known as an \textit{exponential quadratic BSDE} with jumps (\cite{Matoussi2016,Lin2023,Matoussi2019}). Consequently, akin to the Brownian motion case, we develop and utilize the dual relationship between the exponential quadratic BSDE with jumps and the corresponding optimal density process, as well as the duality between the dual BSDE with jumps and its corresponding optimal density. These relationships yield, respectively, the martingale property of the conditional value process and the convex dual representation in the regime switching market setting.

\textbf{Chapter \ref{Chapter7}} is concerned with a consumption-investment problem with unbounded random endowments. When the random endowment is bounded, the problem has been completely solved in \cite{CH} (see also \cite{HSX2022} in the regime switching case). This chapter extends~\cite{CH} from bounded endowments to the unbounded case. The characterizing quadratic BSDE will involve a generator depending on the solution component $Y$. However, we observe that the generator is \textit{decreasing} and Lipschitz continuous in $Y$, resulting in the inclusion of $Y$ in the generator having no impact on the exponential integrability of the solution component $Y$. On the other hand, in contrast to the investment-only case where the dual BSDE can be expressed as a conditional expectation, the consumption-investment setting renders the generator of the dual BSDE nonlinear in the solution component $Y$. Consequently, the corresponding conditional expectation becomes recursive in nature.    

\textbf{Chapter \ref{Chapter8}} is concerned with the consumption-investment problem of an investor with Epstein-Zin recursive utility. A notable departure from earlier models is the \textit{unbounded} nature of the market coefficients, particularly in scenarios where the market price of risk is unbounded. The Epstein-Zin recursive utility, initially introduced in \cite{EpsteinZin}, represents a generalization of the traditional expected power utility. It allows for the separation of time aggregation and risk aggregation, thereby permitting the capturing of intertemporal substitution and risk aversion within a more sophisticated and general framework. The associated consumption-investment problem leads to a quadratic BSDE with unbounded coefficients, as first demonstrated in Xing \cite{Xing2017}. See also Kraft et al.~ \cite{Kraft2017} for a Markovian case and more recently Feng and Tian \cite{Feng2023} for a special portfolio constraint. Additionally, see Matoussi and Xing \cite{Matoussi2018} for a convex duality representation, Aurand and Huang \cite{Huang2023} and Herdegen et al.~\cite{Herdegen2023(1), Herdegen2023(2)} for extensions to an infinite horizon framework. Notably, the aforementioned works do not cover (general) portfolio constraints, which constitutes one of the primary contributions of our work. In addition, most existing works focus on the case when both relative risk aversion and elasticity of intertemporal substitution are greater than one, excluding the traditional expected power utility maximization.

Our objective is to consider cases where relative risk aversion is allowed to be less than one, so the expected power utility maximization is covered as a special case when the relative risk aversion and the elasticity of intertemporal substitution are mutually reciprocals. Note that the case where the elasticity of intertemporal substitution is less than one is still open. The corresponding quadratic BSDE, though with null terminal value, has \textit{unbounded coefficients} in its generator, and moreover, the generator is only \textit{decreasing} in the solution component $Y$ with no Lipschitz continuity property. We show that under an exponential integrability condition on the market price of risk, the corresponding quadratic BSDE admits a unique solution bounded from above (when the relative risk aversion is greater than one) and bounded from below (when the relative risk aversion is less than one). The duality between the quadratic BSDE and the optimal density process is key to establishing the martingale property of the conditional value process, while the duality between the dual BSDE and its optimal density process is essential to establishing the convex dual representation and the uniqueness of the solution of the characterizing quadratic BSDE with unbounded coefficients.\\

The rest of the memoir is organized as follows: In Chapter \ref{Chapter2}, we introduce the market model and the utility maximization problem. Chapter \ref{Chapter3} provides background results on the theory of quadratic BSDEs with unbounded solutions. Chapters 4 to 8 present the main results of this paper. Finally, Chapter 9 concludes with some future research directions.

\chapter{The optimal investment model}\label{Chapter2}

\section{Preliminaries}
Let us fix a real number $T>0$ as the maturity.
Let $B=(B_t)_{t\ge 0}$
be a standard $m$-dimensional Brownian motion defined on a
complete probability space $(\Omega,{\mathcal F},\mathbb P)$, and
$\{{\mathcal F}_t\}_{t\ge 0}$ be the augmented natural filtration of $B$
which satisfies the \emph{usual conditions}. The sigma field of predictable subsets of
$[0,T]\times\Omega$ is denoted by ${\mathcal P}$, and a stochastic
process is called predictable if it is measurable with respect to
$\mathcal{P}$.

Consider a financial market consisting of one risk-free bond with
interest rate zero and $d \leq m$ stocks. In the case $d<m$, we face
an incomplete market. The price process of the stock $i$ evolves
according to the equation
\begin{equation} \label{price}
	\frac{dS^i_t}{S^i_t}= b^i_t dt + \sigma^i_t dB_t, \quad
	i=1,\ldots,d,
\end{equation}
where $b^i$ (resp. $\sigma^i$) is an $\mathbb R$-valued (resp.
$\mathbb R^{m}$-valued) predictable bounded stochastic process. 
The volatility matrix $\sigma_t =
(\sigma^1_t,\ldots,\sigma^d_t)^{tr}$ has full rank, i.e.
$\sigma_t\sigma_t^{tr}$ is invertible, $\mathbb P$-a.s., for
$t\in[0,T]$. Define the market price of risk/market price of risk as an $\mathbb R^m$-valued
predictable process
\begin{equation}\label{risk_premium}
\theta_t = \sigma_t^{tr} (\sigma_t \sigma_t ^{tr})^{-1}b_t, \quad t \in [0,T],
\end{equation}
and assume that $\theta$ is also \emph{bounded}\footnote{In Chapter \ref{Chapter8}, we will extend to the case of unbounded market price of risk.}.  Throughout, we will be
using $A^{tr}$ to denote the transpose of matrix $A$. Consequently,
the market price of risk $\theta$ defined in (\ref{risk_premium}) solves $$\min_{\theta_t}|\theta_t|^2,$$ subject to the equation for the market price of risk:
\[\sigma_t\theta_t=b_t, \quad t\in[0,T].\]

In this market environment, an investor trades dynamically among the
risk-free bond and the risky assets. For $1\le i\le d$, let
$\pi^i_t$ denote the amount of money invested in stock $i$ at time
$t$, so the number of shares is $\frac{\pi^i_t}{S^i_t}$.  An
$\mathbb{R}^d$-valued predictable process $\pi=(\pi_t)_{0 \leq t\leq
	T}$ is called a {self-financing trading strategy} if
$\int_0^{\cdot} \pi^{tr}_t \frac{dS_t}{S_t}$ is well defined, for
example, $\sigma^{tr}\pi\in L^2[0,T]$, i.e. $\int_0^T | \sigma_t^{tr} \pi_t|^2dt < \infty$,
$\mathbb{P}$-a.s., and the corresponding wealth process $X^{\pi}$
with initial capital $x$ satisfies the equation
\begin{equation}\label{wealth}
	X^{\pi}_t= x +\sum_{i=1}^d \int_0^t
	\frac{\pi^i_{u}}{S^i_{u}} dS^{i}_{u} =x + \int_0^t \pi_u^{tr}
	\sigma_u (dB_u + \theta_u du),\quad x\in\mathbb{R}.
\end{equation}

The investor has an exponential utility with respect to their
terminal wealth $X_T^{\pi}$. We recall that, for $\alpha>0$, an
exponential utility function is defined as
\[ U(x) = -\exp(- \alpha x),\quad x\in\mathbb R. \]
In addition to the terminal wealth $X_T^{\pi}$, the investor
also pays or receives an $\mathcal{F}_T$-measurable random
endowment/payoff $F$ at maturity $T$. $F$ means a
payment when it is nonnegative, and  an income when it is negative.

The investor chooses an \emph{admissible}
self-financing trading strategy $\pi^{\star}$ so as to maximize
the expected utility of the net wealth at maturity $T$:
\begin{equation}\label{expoopt} V(0,x):=\sup_{\substack{\{\pi_u,\ u\in[0,T]\}\\{\text{admissible}}}} E \left[ -\exp \left(-\alpha\left( x + \int_0^T
	\pi_t^{tr} \frac{dS_t}{S_t} -F \right) \right)  \right],
\end{equation}
where $V(0,\cdot)$ is called the \emph{value function} at initial
time $0$. To solve (\ref{expoopt}), we need to further choose an
admissible set from which we select the optimal trading strategy
$\pi^{\star}$. Different admissible sets and different assumptions
on the payoff $F$ may lead to different solutions.

With the help of the theory of quadratic BSDEs with
\emph{bounded} terminal data and martingale optimality principle, Hu et al.~\cite{him} solved the above
optimization problem (\ref{expoopt}) under the assumption that $F$
is bounded, and $\pi$ takes values in the admissible set
$\mathcal{A}_D$ defined as follows.

\begin{definition}\label{admiss}
	[{\bf Admissible strategies with constraints ${\mathcal
		A}_D$}]
	
	Let ${\mathcal C}$ be a closed set in $\mathbb R^{d}$ with $0\in\mathcal{C}$.
The set of admissible trading strategies ${\mathcal
		A}_D$ consists of all $\mathbb{R}^{d}$-valued predictable processes
	$\pi\in L^2[0,T]$, which satisfy $\pi_t \in
	{\mathcal C}$, $\mathbb P$-a.s., for $t\in[0,T]$. Moreover, the
	following Class (D) condition holds: the family of random variables 
	$$\{
	\exp(-\alpha X_{\tau}^{\pi}):\ \tau\ \mbox{is a stopping time taking
		values in}\ [0,T]\}
	$$ are uniformly integrable.
\end{definition}

On one hand, the assumption of boundedness on the payoff $F$ leads to the exclusion of various interesting cases, such as call options. On the other hand, the technical nature of the Class (D) condition within the admissible set $\mathcal{A}_D$ appears unnecessary, particularly when $F$ is bounded from below and exponentially integrable.

We aim to relax the above two assumptions, which
are crucial to the proofs in \cite{him}, by leveraging the elements from
the theory of quadratic BSDE with \emph{unbounded} terminal data, coupled with the convex duality theory.

\section{Unbounded payoffs}

We observe that the minimal condition on the payoff $F$ should
guarantee that the expectation in (\ref{expoopt}) is finite with
$\pi_t\equiv0$, namely,
\[E[e^{\alpha F}]<+\infty.\]
Intuitively, the investor should have a finite expected utility when
 all the money is put in the risk-free bond, so $F$ is not so
bad that doing nothing leads to a prohibitive punishment, i.e. $E[U(F)]=-\infty$ is excluded. 
 
Moreover, it is natural to require
$$E^{\mathbb{Q}}[|F|]<+\infty$$ under different equivalent
probability measures $\mathbb{Q}$, i.e. the expected payoff (under
different equivalent probability measures) should be finite, so $F$
is not too good to be true. The above discussions motivate us to
impose the following assumption on $F$, which will hold throughout
the paper.

\begin{assumption}\label{assumption1} The payoff $F$ satisfies the
	exponential integrability condition
	\begin{equation}\label{condition}
		E[e^{p\alpha F^{+}}]<+\infty; \quad E[e^{\varepsilon F^{-}}]<+\infty
	\end{equation}
	for some integer $p>1$ and positive number $\varepsilon>0$, where
	$F^{+}=\max\{F,0\}$ and $F^{-}=\max\{-F,0\}$.
\end{assumption}

{Typical examples include call option payoffs of the form $\pm (V_T - K)^+$ written on a stochastic factor process $V$ governed by 
$$
dV_{t}=\eta \left( V_{t}\right) dt+\kappa dB_{t}.
$$
for some deterministic function $\eta$ and constant vector $\kappa$.
See Chapter \ref{section:example} for more details. }

\begin{remark} By requiring $p>1$, we require more exponential integrability
	on $F^+$. This is to guarantee the finite entropy condition for both
	$L_T^{q^{\star}}$ (see Lemma \ref{lemma}), which will in turn be used to verify the
	Class (D) condition.
	
	On the other hand, H\"older's inequality might indicate that $F^-$
	being $L^{p}$-integrable is sufficient to guarantee that $F^{-}$ is
	integrable under different equivalent probability measures since 
 \begin{align*}
 E^{\mathbb{Q}}[F^-]&=E\left[\frac{d\mathbb{Q}}{d\mathbb{P}}F^-\right]\\
 &\leq E\left[\left(\frac{d\mathbb{Q}}{d\mathbb{P}}\right)^{\frac{p}{p-1}}\right]^{\frac{p-1}{p}} E\left[(F^-)^{p}\right]^{\frac{1}{p}}.
 \end{align*}
 However,
	the assertion relies on $L^{\frac{p}{p-1}}$-integrability of the
	corresponding Radon-Nikodym density process, which does not always
	hold (e.g. the density process $L^{q}$ in Theorem
	\ref{theorem_dual}). For this reason, we require exponential
	integrability of $F^{-}$. Then, we only need the density process
	with finite entropy.
	
	A similar type of asymmetric exponential integrability condition on
	the terminal data $F$ also appears in Delbaen et al \cite{dhr} and Frei and Schweizer \cite{FS}.
	Nonetheless, it might be possible to relax the assumption $p>1$ on
	$F^+$ by adapting the argument used in Delbaen et al \cite{dhr2}. Such an
	extension is left for the interested reader.
\end{remark}

\section{The underlying admissible strategies}
Since the payoff $F$ satisfies the exponential integrability
condition (\ref{condition}) only,  we need to further strengthen the
Class (D) condition in the admissible set $\mathcal{A}_D$ in order
to solve the optimization problem (\ref{expoopt}).

Let $(\pi_u)_{u\in[0,T]}$ be a given self-financing trading strategy. For any random variable $\xi\in\mathcal{F}_t$ such that ${\xi}$ satisfying Assumption \ref{assumption1}, 
we define the reward functional: 
\begin{equation}\label{expoopt_22}
	V\left(t,\xi;(\pi_u)_{u\in[t,T]}\right):=E \left[\left. -\exp
	\left(-\alpha\left( \xi + \int_t^T \pi_u \frac{dS_u}{S_u} -F
	\right) \right)\right|\mathcal{F}_t\right],
\end{equation}
for $t\in[0,T]$, and the associated \emph{conditional value process } as
\begin{equation}\label{expoopt_2}
	V(t,\xi):= \esssup_{\substack{\{\pi_u,\ u\in[t,
			T]\}\\\text{admissible}}}
	V\left(t,\xi;(\pi_u)_{u\in[t,T]}\right).
\end{equation}

In the linear-quadratic setting, the conditional value process and its definition have been discussed in detail by Tang \cite{Tang2015}. 
Note that the conditional value process  includes the value function in
(\ref{expoopt}) as a special case when $t=0$. Moreover, by taking
$\pi_u\equiv 0$ for $u\in[t,T]$, we obtain a lower bound of the
conditional value process, $V(t,\xi)\geq e^{-\alpha \xi}E[-e^{\alpha
	F}|\mathcal{F}_t]>-\infty$, $\mathbb{P}$-a.s., for $t\in[0,T]$, so
the conditional value process  is always finite.

To solve (\ref{expoopt_2}), we look for $V(\cdot,\cdot)$ such that
$V(t,X_t^{\pi})$, $t\in[0,T]$, is a supermartingale for any
admissible $\pi$, and there exists an admissible $\pi^{\star}$ such
that $V(t,X_t^{\pi^{\star}})$, $t\in[0,T]$, is a martingale, 
which is termed as \textit{martingale optimality principle}. It is
thus natural to impose some integrability conditions on
$V(\cdot,X^{\pi}_{\cdot})$ in the admissible set.

\begin{definition}\label{admiss2}
	[{\bf Admissible strategies with constraints ${\mathcal A}_D^{conv}$}]
	
	The set of admissible trading strategies ${\mathcal A}_D^{conv}$ is the same
	as $\mathcal{A}_D$ in Definition \ref{admiss}, except that the constraint set ${\mathcal C}$
assumed to be a closed and \emph{convex} set in $\mathbb R^{d}$ and $0\in\mathcal{C}$, and the Class
	(D) condition is replaced by the following Class
	(D) condition on the conditional value process,
	$$\{ V(\tau,X_{\tau}^{\pi}):\ \tau\ \mbox{is a stopping time
		taking values in}\ [0,T]\}
	$$ is a uniformly integrable family.
\end{definition}

We use the superscript $^{conv}$ to emphasize the dependency on the convexity of $\mathcal{C}$.
The convexity of $\mathcal{C}$ is crucial in our subsequent analysis which relies on the convex duality theory.
The admissible set $\mathcal{A}_D^{conv}$ depends on the
integrability of
$V(\cdot,X^{\pi}_{\cdot})$. Therefore, in some sense, the admissible set $\mathcal{A}_{D}^{conv}$ also constitutes a part of the solution to be determined.

However, this does not mean there is a loop of dependency herein. In fact,
by the definition of $V(t,\xi)$, it is immediate to check that the following homogeneity property 
\[V(t,X_t^{\pi})=\exp(-\alpha
X_{t}^{\pi})V(t,0).\] In the proof of Theorem \ref{thm:opti2}, we
shall show that $V(t,0)=-e^{\alpha Y_t}$ with $Y$ solving an
upcoming BSDE (\ref{BSDEexpo}). Thus, $\mathcal{A}_D^{conv}$ is
equivalent to say $e^{-\alpha X^{\pi}_{\cdot}}e^{\alpha Y_{\cdot}}$
is in Class (D), which \emph{in prior} has nothing to do with the
optimization problem (\ref{expoopt}).

On the other hand, if $F$ is bounded, then the admissible set
$\mathcal{A}_D^{conv}=\mathcal{A}_D$, which is independent of the
conditional value process. Indeed, if $F$ is bounded, then $Y$ is also
bounded, so are $V(t,0)$ and $\frac{1}{V(t,0)}$ (see Theorem 7 of
\cite{him} or Lemma \ref{lemma_boundedcase} for its proof). In this case, the Class (D) condition on
the conditional value process  $V(\cdot,X_{\cdot}^{\pi})$ is equivalent
to the Class (D) condition on the exponential utility of the wealth
$e^{-\alpha X_{\cdot}^{\pi}}$, and therefore,
$\mathcal{A}^{conv}_{D}$ coincides with $\mathcal{A}_{D}$.

Our utility maximization problem is therefore formulated as follows:

\emph{Solving the optimization problem (\ref{expoopt}) with the
	payoff $F$ satisfying Assumption \ref{assumption1} and  admissible
	set $\mathcal{A}_D^{conv}$ as in Definition \ref{admiss2}.}

%
%
%


\chapter{Preliminaries on quadratic BSDEs}\label{Chapter3}

The theory of quadratic BSDEs will play a crucial role in the solution of  utility maximization problems. A BSDE with terminal condition $F$ and generator $f$ is an equation
of the following type
\begin{equation}\label{BSDEexpo}
	Y_t=F+\int_t^T f(s,Y_s, Z_s)ds-\int_t^T Z_s^{tr}dB_s,\quad t\in[0,T],
\end{equation}
and is often denoted by BSDE$(F,f)$. Recall that a generator is a
random function $f:[0,T]\times \Omega\times\mathbb R^m\rightarrow
\mathbb R$, which is measurable with respect to ${\mathcal
	P}\otimes{\mathcal B}(\mathbb R^m)$, and a terminal condition is a
real-valued ${\mathcal F}_T$-measurable random variable $F$.

By a solution to BSDE$(F,f)$, we mean a pair of predictable
processes $(Y,Z)=(Y_t,Z_t)_{t\in [0,T]}$, with values in $\mathbb
R\times\mathbb R^m$ such that $\mathbb P$-a.s., $t\mapsto Y_t$ is
continuous, $t\mapsto Z_t$ belongs to $L^2(0,T)$, i.e.
$\int_0^T|Z_t|^2dt<+\infty$, $t\mapsto f(t,Z_t)$ belongs to
$L^1(0,T)$, i.e. $\int_0^T|f(t,Z_t)|dt<\infty$, and $(Y,Z)$ satisfies (\ref{BSDEexpo}).

\section{A convex generator for a quadratic BSDE}

It turns out we need to solve a quadratic BSDE$(F,f)$ with its generator $f$ given by
\begin{align}\label{driver}
		f(t,z) &= \frac{\alpha}{2}\min_{\pi\in {\mathcal C}}
		\left|\sigma_t^{tr}\pi-(z +
		\frac{1}{\alpha} \theta_t)\right|^2 - z^{tr}\theta_t - \frac{1}{2 \alpha}
		|\theta_t|^2\notag\\
          &=\frac{\alpha}{2}\mbox{dist}_{\sigma_t^{tr}\mathcal{C}}^2\Big(z+\frac{1}{\alpha}\theta_t\Big)
          - z^{tr}\theta_t - \frac{1}{2 \alpha}
		|\theta_t|^2,\quad (t,z)\in[0,T]\times\mathbb{R}^d,
	\end{align}
where $\mbox{dist}_{\sigma_t^{tr}\mathcal{C}}(\cdot)$ is the
distance function of $\sigma_t^{tr}\mathcal{C}\in\mathbb{R}^{m}$, which is closed and convex due to the full rank condition on $\sigma_t$.
Let $\pi^*_t$ be such that
\begin{equation}\label{pi_2}	\sigma_t^{tr}\pi^{\star}_t=\mbox{Proj}_{\sigma_t^{tr}\mathcal{C}}\left(z+\frac{\theta_t}{\alpha}\right),\
		\mbox{$\mathbb{P}$-a.s.},\ \text{for}\ t \in [0,T],
	\end{equation}
	where $\mbox{Proj}_{\sigma_t^{tr}\mathcal{C}}(\cdot)$ is the
	projection operator on
	$\sigma_t^{tr}\mathcal{C}$.
Then, the generator $f$ is equivalent to
\begin{equation*}\label{driver1}
	f(t,z) = \frac{\alpha}{2}
	\left|\mbox{Proj}_{\sigma_t^{tr}\mathcal{C}}(z +
	\frac{1}{\alpha} \theta_t)-(z +
	\frac{1}{\alpha} \theta_t)\right|^2 - z^{tr}\theta_t - \frac{1}{2 \alpha}
	|\theta_t|^2.
\end{equation*}

The following lemma about the growth property and local Lipschitz continuity of the generator $f$ is immediate.

\begin{lemma}\label{lemma_generator} Let $\mathcal{C}$ be a closed and convex set in $\mathbb{R}^d$ satisfying $0\in\mathcal{C}$. Then, the generator $f$ admits the following properties:
for
$(t,z,\bar{z})\in[0,T]\times\mathbb{R}^{m}\times\mathbb{R}^m$,

(i) (Growth property) \begin{equation*}\label{boundforgenerator}
-z^{tr}\theta_t-\frac{1}{2\alpha}|\theta_t|^2\le
f(t,z)\le \frac{\alpha}{2}|z|^2.
\end{equation*}

(ii) (Local Lipschitz continuity)
\begin{equation*}\label{Lipforgeneraor}
|f(t,z)-f(t,\bar{z})|\leq \Big((1+\frac{2}{\alpha})|\theta_t|+|z|+|\bar{z}|\Big)|z-\bar{z}|.
\end{equation*}
\end{lemma}

While the properties discussed above for the generator function $f$ are sufficient to address bounded payoffs, the case of unbounded payoffs crucially depends on the convexity of $f$. We observe that when $\mathcal{C}$ is convex, the generator $f(t,z)$, defined in (\ref{driver}), is convex in
$z$. We can therefore introduce the convex dual of $f(t,z)$,
\begin{equation}\label{dual}
	{f}^{\star}(t,q):=\sup_{z\in\mathbb{R}^m}\left({z}^{tr}q-f(t,z)\right),
\end{equation}
for $(t,q)\in[0,T]\times\mathbb{R}^m$. Note that $f^{\star}$ is
valued in $\mathbb{R}\cup\{+\infty\}$.  Using the upper bound for the generator $f$ in Lemma \ref{lemma_generator}, we obtain a lower bound of
	$f^{\star}(t,q)$, namely, for any $z\in\mathbb{R}^m$,
	\begin{align}\label{lowerbound}
		f^{\star}(t,q)&\geq z^{tr}q-f(t,z)\notag\\
		&\geq z^{tr}q-\frac{\alpha}{2}|z|^2\geq \frac{|q|^2}{2\alpha},
	\end{align} by taking $z=q/\alpha$ in the last inequality.

Since $f(t,z)$ is continuous and convex in $z$, the Fenchel-Moreau theorem yields that
\begin{equation}\label{FM_theorem}
	f(t,z)=\sup_{q\in\mathbb{R}^m}({z}^{tr}q-f^{\star}(t,q)),
\end{equation}
for $(t,z)\in[0,T]\times\mathbb{R}^m$. Moreover,
$q^{\star}\in\partial f_{z}(t,z)$, which is the subdifferential of
${z}\mapsto f(t,z)$ at ${z}\in\mathbb{R}^m$, achieves the supremum
in (\ref{FM_theorem}),
\begin{equation}\label{FM_theorem_2}
	f(t,z)={z}^{tr}q^{\star}-f^{{\star}}(t,q^{\star}).
\end{equation}
Recall that $\partial f_{z}(t,z)$ is defined as
\[\partial f_z(t,z)=\{q\in\mathbb{R}^m: f(t,\bar{z})-f(t,z)\geq q^{tr} (\bar{z}-z),\ \text{for}\ \bar{z}\in\mathbb{R}^m\}.\]

We calculate $\partial f_z(t,z)$ in the following lemma.

\begin{lemma}\label{lemma_convex} Let $\mathcal{C}$ be a closed and convex set in $\mathbb{R}^d$ satisfying $0\in\mathcal{C}$. Then, for $(z +
		\frac{1}{\alpha} \theta_t)\notin\sigma_t^{tr}\mathcal{C}$,
the subdifferential of $z\mapsto f(t,z)$ satisfies
\begin{equation}\label{subdifferential}
\alpha(z-\sigma_t^{tr}\pi_t^{\star})\in \partial f_z(t,z),
\end{equation}
with $\sigma_t^{tr}\pi^*_t$ given in (\ref{pi_2}). Otherwise, for $(z +
		\frac{1}{\alpha} \theta_t)\in\sigma_t^{tr}\mathcal{C}$, $\partial f_z(t,z)=-\theta_t$.
\end{lemma}

\begin{proof}  Note that
\begin{equation}\label{derivative}
\partial_zf(t,z)=\alpha \mbox{dist}_{\sigma_t^{tr}\mathcal{C}}\Big(z+\frac{1}{\alpha}\theta_t\Big)\partial_{z}\mbox{dist}_{\sigma_t^{tr}\mathcal{C}}\Big(z+\frac{1}{\alpha}\theta_t\Big)-\theta_t.
\end{equation}
Hence, it is sufficient to prove the case $(z +
		\frac{1}{\alpha} \theta_t)\notin\sigma_t^{tr}\mathcal{C}$.

Introduce the tanget cone of $\sigma_t^{tr}\mathcal{C}$ at the projection point $\sigma_t^{tr}\pi^{\star}_t$:
\begin{equation}
    N^{\star}_{\sigma_t^{tr}\mathcal{C}}(\sigma_t^{tr}\pi_t^{\star})=
    \Big\{y\in\mathbb{R}^{m}:\Big(z+\frac{\theta_t}{\alpha}-\sigma_t^{tr}\pi_t^{\star}\Big)^{tr}\Big(y-\sigma_t^{tr}\pi_t^{\star}\Big)\leq 0\Big\}.
\end{equation}
Note that $0\in \sigma_t^{tr}\mathcal{C}\subset  N^{\star}_{\sigma_t^{tr}\mathcal{C}}(\sigma_t^{tr}\pi_t^{\star})$.
In turn, for $y\notin N^{\star}_{\sigma_t^{tr}\mathcal{C}}(\sigma_t^{tr}\pi_t^{\star})$,
\begin{equation*}
\mbox{dist}_{\sigma_t^{tr}\mathcal{C}}(y)
\geq \mbox{dist}_{N^{\star}_{\sigma_t^{tr}\mathcal{C}}}(y)
=\frac{\Big(z+\frac{\theta_t}{\alpha}-\sigma_t^{tr}\pi_t^{\star}\Big)^{tr}\Big(y-\sigma_t^{tr}\pi_t^{\star}\Big)}{ \mbox{dist}_{\sigma_t^{tr}\mathcal{C}}\Big(z+\frac{1}{\alpha}\theta_t\Big)}.
\end{equation*}
On the other hand, for $y\in N^{\star}_{\sigma_t^{tr}\mathcal{C}}(\sigma_t^{tr}\pi_t^{\star})$,
\begin{equation*}
\mbox{dist}_{\sigma_t^{tr}\mathcal{C}}(y)\geq 0\geq\frac{\Big(z+\frac{\theta_t}{\alpha}-\sigma_t^{tr}\pi_t^{\star}\Big)^{tr}\Big(y-\sigma_t^{tr}\pi_t^{\star}\Big)}{ \mbox{dist}_{\sigma_t^{tr}\mathcal{C}}\Big(z+\frac{1}{\alpha}\theta_t\Big)}.
\end{equation*}
In both situations, for $y\in\mathbb{R}^m$, we have
\begin{align*}
\mbox{dist}_{\sigma_t^{tr}\mathcal{C}}(y)&\geq \frac{\Big(z+\frac{\theta_t}{\alpha}-\sigma_t^{tr}\pi_t^{\star}\Big)^{tr}\Big(y-
(z+\frac{1}{\alpha}\theta_t)+(z+\frac{1}{\alpha}\theta_t)-
\sigma_t^{tr}\pi_t^{\star}\Big)}{ \mbox{dist}_{\sigma_t^{tr}\mathcal{C}}\Big(z+\frac{1}{\alpha}\theta_t\Big)}\\
&= \frac{\Big(z+\frac{\theta_t}{\alpha}-\sigma_t^{tr}\pi_t^{\star}\Big)^{tr}}{ \mbox{dist}_{\sigma_t^{tr}\mathcal{C}}\Big(z+\frac{1}{\alpha}\theta_t\Big)}\Big(y-
(z+\frac{1}{\alpha}\theta_t)\Big)+\mbox{dist}_{\sigma_t^{tr}\mathcal{C}}\Big(z+\frac{1}{\alpha}\theta_t\Big),
\end{align*}
which shows that
\begin{equation*}
 \frac{\Big(z+\frac{\theta_t}{\alpha}-\sigma_t^{tr}\pi_t^{\star}\Big)}{ \mbox{dist}_{\sigma_t^{tr}\mathcal{C}}\Big(z+\frac{1}{\alpha}\theta_t\Big)}\in \partial_{z}\mbox{dist}_{\sigma_t^{tr}\mathcal{C}}\Big(z+\frac{1}{\alpha}\theta_t\Big).
\end{equation*}
The conclusion then follows by subsiting the above expression into (\ref{derivative}).
\end{proof}

\section{Quadratic BSDE with bounded terminal data}
We first recall the existence and uniqueness theorem for quadratic BSDE with bounded terminal data. The existence follows from Theorem 2.3 of Kobylanski \cite{Kobylanski} and uniqueness follows from Theorem 7 of Hu et al.~\cite{him}.

Let $\mathcal{S}^{\infty}$ be the space of real-valued, adapted and c\`adl\`ag bounded processes with its norm denoted by $||\cdot||_{\mathcal{S}^{\infty}}$. For $p\ge 1$, ${\mathcal S}^p$ denotes the
space of real-valued, adapted and c\`adl\`ag processes $(Y_t)_{t\in
	[0,T]}$ such that
$$||Y||_{{\mathcal S}^p}:=E\left[\sup_{t\in [0,T]}|Y_t|^p\right]^{1/p}<+\infty,$$
and $M^p$ denotes the space of (equivalent classes of)
$\mathbb{R}^m$-valued predictable processes $(Z_t)_{t\in [0,T] }$
such that
$$||Z||_{M^p}:= E\Bigr[\Bigr(\int_0^T|Z_s|^2ds\Bigr)^{p/2}\Bigr]^{1/p}<+\infty.$$

\begin{lemma}[\cite{him}]\label{theorem: boundedBSDE} Suppose that the payoff $F$ is bounded. Then, BSDE$(F,f)$ with $f$ given by (\ref{driver}) admits a uniqueness solution $(Y,Z)$, where $(Y,Z)\in\mathcal{S}^{\infty}\times M^2$.
\end{lemma}

\begin{remark}\label{remark_BMO} Due to the local Lipschitz continuity of $f(t,z)$ in $z$ (see Lemma \ref{lemma_generator}(ii)),  it is shown in Lemma 12 of \cite{him} that the stochastic process $\int_0^{s}Z_t^{tr}dW_t$, $s\in[0,T]$, is a BMO martingale (see Kazamaki \cite{Kaz}). In other words,
$$||\int_0^{\cdot}Z^{tr}_tdB_t||_{BMO}=\sup_{\tau\ \text{is an}\ \mathbb{F}\text{-stopping time}}\mathop{\text{esssup}}_\omega\  {E}\left[\int_{\tau}^T|Z_t|^2dt\Bigm |\mathcal{F}_{\tau}\right]^{1/2}<+\infty.
$$
As a consequence,  the stochastic exponential $L=\mathcal{E}(\int_0^{\cdot}Z_t^{tr}dB_t)$ is a uniformly integrable martingale. The BMO martingale property has been extensively used in the solution of  the utility maximization problem with bounded payoff $F$ in \cite{him}, which will be reviewed in Chapter \ref{section:review}.
\end{remark}

The stability property for quadratic BSDE with bounded terminal data plays a pivotal role when passing to the limit within a sequence of approximating equations. It follows from Proposition 2.4 in Kobylanski \cite{Kobylanski} (see also Lemma 2 in Briand and Hu \cite{BH1} and Lemma 3 in Briand and Hu \cite{BH2}).

\begin{lemma}[\cite{Kobylanski}]\label{lemma_stability}
Consider the following sequence of quadratic BSDE with bounded data:
\begin{equation}\label{QBSDE_bounded}
	Y_t=F^n+\int_t^T f^n(s,Y_s, Z_s)ds-\int_t^T Z_s^{tr}dB_s,\quad t\in[0,T],
\end{equation}
satisfying the following conditions: 

(i) $(F^n)_{n\geq 1}$ is a sequence of $\mathcal{F}_T$-measurable uniformly bounded random variables; 

(ii) $(f^{n})_{n\geq 1}$ is a sequence of continuous generators satisfying the monotonic condition in $y$: 
$$y(f^n(t,y,z)-f(t,0,z))\leq \beta |y|^2$$
and the growth condition in $z$:
$$|f^n(t,y,z)|\leq |\theta_t|+\phi(|y|)+\frac{\alpha}{2}|z|^2,$$
for two constants $\alpha,\beta\geq 0$ and a deterministic continuous and increasing function $\phi$ with $\phi(0)=0$, and an $\mathbb{F}$-progressively measurable process $\theta$ such that $\int_0^T|\theta_t|dt$ is bounded. 

Suppose that $F_n\rightarrow F$, $\mathbb{P}$-a.s., and $f^n(t,z_n)\rightarrow f(t,z)$ when $z_{n}\rightarrow z$, $\mathbb{P}$-a.s., for $t\in[0,T]$. If BSDE$(F^n,f^n)$ admits a solution $(Y^n,Z^n)\in\mathcal{S}^{\infty}\times {M}^2$ such that $Y^n$ is increasing (resp. decreasing) and $\sup_{n\geq 1}||Y_n||_{S^{\infty}}\leq C$, then there exists a limit pair $(Y,Z)\in \mathcal{S}^{\infty}\times {M}^2$ as a solution to BSDE$(F,f)$ such that $Y^n\rightarrow Y$ uniformly in $[0,T]$ in probability and $Z^n\rightarrow Z$ in ${M}^2$.
\end{lemma}

\begin{remark}\label{remark_approx}
When dealing with the bounded terminal data $F$ and a generator $f$ as described in (\ref{driver}), we typically resort to approximating it through \textit{inf convolution} as introduced in \cite{Lepeltier1997}: For $n\geq |\theta_t|$, consider
\begin{equation}\label{inf_convolution}
f^n(t,z)=\inf_{q}\{f(t,q)+n|q-z|\}.
\end{equation}
Then, $f^n$ is well defined and globally Lipschitz continuous with constant $n$, so the comparison theorem for Lipschitz BSDE holds for BSDE$(F,f^n)$. 
This point is particularly useful when the standard comparison theorem for quadratic BSDE cannot be applied directly, such as when the generator involves a non-Lipschitz continuous term in 
$y$ or when the equation becomes a system.
Moreover, $f^n$ is increasing and converges point-wisely to $f$. In turn, Dini's theorem implies that convergence is also uniform on compact sets. Finally, we have the following growth for $f^n$:
$$-|\theta_t||z|-\frac{1}{2\alpha}|\theta_t|^2\le
f^n(t,z)\leq f(t,z)\le \frac{\alpha}{2}|z|^2.$$

On other hand, in the case of an unbounded terminal data $F$, it is common to employ an approximation technique through \textit{truncation}. For instance, we approximate $F$ using truncation as follows: $F^{n,k}=F^+\wedge n - F^{-}\wedge k$, where $n$ and $k$ are positive integers. However, the stability property mentioned above does not apply because $F^{n,k}$ is not uniformly bounded. Therefore, a localization argument, as introduced by Briand and Hu \cite{BH1, BH2}, needs to be employed. This localization argument will be discussed in further detail in the next section.
\end{remark}

\section{Quadratic BSDEs with unbounded terminal data}

Next, we present the existence and uniqueness theorem for a
quadratic BSDE with the terminal data satisfying the exponential
integrability condition (\ref{condition}). It will be subsequently
used to solve the optimization problem (\ref{expoopt}).


\begin{theorem}\label{theorem:BSDE1}
	Suppose that Assumption \ref{assumption1} holds. Then BSDE$(F,f)$ with the generator $f$ given in (\ref{driver}) admits a unique solution $(Y,Z)$, where $e^{\alpha
		Y^+}\in\mathcal{S}^{p}$, $e^{\varepsilon Y^-}\in\mathcal{S}^{1}$,
	and $Z\in M^2$, i.e.
	\begin{equation*}
		E\left[e^{p\alpha Y^{+}_{\star}}+e^{\varepsilon
			Y^{-}_{\star}}+\int_0^T|Z_s|^2ds\right]<+\infty,
	\end{equation*}
	where $Y_{\star}=\sup_{t\in[0,T]}Y_t$ is the running maximum of a
	stochastic process $Y$.
	
	Moreover, if $E[e^{p^{\prime} |F|}]<+\infty$ for any $p^{\prime}\ge
	1$, then $e^{\alpha Y}\in\mathcal{S}^{p^{\prime}}$, and $Z\in
	M^{p^{\prime}}$, i.e.
	$$E\left[e^{p^{\prime}\alpha Y_{\star}}+\left(\int_0^T |Z_s|^2ds\right)^{p^{\prime}/2}\right]<+\infty.$$
\end{theorem}

\begin{proof}

The main ideas for establishing the existence of the solution to BSDE (\ref{BSDEexpo}) are adapted from Briand and Hu \cite{BH1, BH2}, while the uniqueness of the solution to BSDE (\ref{BSDEexpo}) is inspired by Delbaen et al.~\cite{dhr}.

To this end, we truncate the terminal data $F$ as in Remark \ref{remark_approx},
$$F^{n,k}=F^+\wedge n-F^-\wedge k,\quad \mbox{for integers}\
n,k\geq 1.$$ Then, it follows that $|F^{n,k}|\le \max\{n,k\},$ and
\[-F^-\leq F^{n,k+1}\leq F^{n,k}\leq F^{n+1,k}\leq F^+.\]
Moreover, $\lim_{n,k\rightarrow\infty}F^{n,k}=F$.

With $F^{n,k}$ at hand, we first consider the following truncated BSDE$(F^{n,k},f)$
\begin{equation}\label{BSDEn*}
	Y_t^{n,k}=F^{n,k}+\int_t^T f(s,Z_s^{n,k})ds-\int_t^T
	(Z_s^{n,k})^{tr}dB_s,\quad t\in [0,T].
\end{equation}
It admits a unique solution
$(Y^{n,k},Z^{n,k})\in\mathcal{S}^{\infty}\times M^{2}$ by Lemma \ref{theorem: boundedBSDE}. Moreover,
applying the comparison
theorem for quadratic BSDE with bounded terminal data (see Theorem
2.6 in \cite{Kobylanski}) to (\ref{BSDEn*}) 
yields that $Y^{n,k}$ is increasing in $n$ and decreasing in $k$ and, moreover,
\begin{equation}\label{bound}
	\underline{Y}_t\le Y_t^{n,k}\le
	\overline{Y}_t,
\end{equation}
where $\underline{Y}\in\mathcal{S}^{\infty}$ solves
BSDE$(F^{n,k},-z^{tr}\theta-\frac{1}{2\alpha}|\theta|^2)$, namely,
\begin{equation}\label{BSDEl*}
	\underline{Y}_t=F^{n,k}+\int_t^T
	(-\underline{Z}_s^{tr}\theta_s-\frac{1}{2\alpha}|\theta_s|^2
	)ds-\int_t^T \underline{Z}_s^{tr}dB_s,\ t\in[0,T],
\end{equation}
and $\overline{Y}\in\mathcal{S}^{\infty}$ solves BSDE$(F^{n,k},\frac{\alpha}{2}|z|^2)$, namely,
\begin{equation}\label{BSDE2*}
	\overline{Y}_t=F^{n,k}+\int_t^T
	\frac{\alpha}{2}|\overline{Z}_s|^2ds-\int_t^T
	\overline{Z}_s^{tr}dB_s,\ t\in[0,T].
\end{equation}

It is routine to check that both $\underline{Y}$ and $\overline{Y}$
have the explicit expressions
\begin{align*}
	\underline{Y}_t&=E^{\mathbb{Q}^{\theta}}\left[F^{n,k}-\int_t^T
	\frac{1}{2\alpha}|\theta_s|^2ds\Bigm |{\mathcal F}_t\right]\geq -E^{\mathbb{Q}^{\theta}}\left[F^{-}+\int_t^T
	\frac{1}{2\alpha}|\theta_s|^2ds\Bigm |{\mathcal F}_t\right];\\
	\overline{Y}_t&=\frac{1}{\alpha}\ln E\left [e^{\alpha F^{n,k}}\bigm |\mathcal{F}_t\right]\leq \frac{1}{\alpha}\ln E\left[e^{\alpha F^+}\bigm |\mathcal{F}_t\right],
\end{align*}
where $\mathbb{Q}^{\theta}$ is the minimal local martingale measure (MLMM) given in (\ref{EMM}) in the next Chapter.


Next, we pass to the limit in (\ref{BSDEn*}). Since the solution component
$Y^{n,k}$ is only locally bounded, we need to apply the localization
method introduced in \cite{BH1} (see also \cite{BH2}). For integer
$j\ge 1$, we introduce the following stopping time
\[
\tau_j = T \wedge \inf \left\{ t \in [0,T]: \max \left\{ \frac{1}{\alpha} \ln \mathbb{E} \left[ e^{\alpha F^+} \big| \mathcal{F}_t \right], \mathbb{E}^{\mathbb{Q}^{\theta}} \left[ F^- + \int_t^T \frac{1}{2\alpha} |\theta_s|^2 ds \big| \mathcal{F}_t \right] \right\} > j \right\}.
\]

Then $(Y_j^{n,k}(t),Z_j^{n,k}(t)):=(Y_{t\wedge
	\tau_j}^{n,k},Z_t^{n,k} 1_{\{t\le\tau_j\}})$, $t\in[0,T]$, satisfies
\begin{equation}\label{BSDEkn*}
	Y_j^{n,k}(t)=F_j^{n,k}+\int_t^T 1_{\{s\le\tau_j\}}
	f(s,Z_j^{n,k}(s))ds-\int_t^T \left(Z_j^{n,k}(s)\right)^{tr}dB_s,
\end{equation}
where $F_j^{n,k}=Y_j^{n,k}(T)=Y^{n,k}_{\tau_j}$.

For fixed $j$, $Y^{n,k}_{j}(\cdot)$ is increasing in $n$ and
decreasing in $k$, while it remains bounded by $j$. Hence,
setting $Y_{j}(t)=\inf_{k}\sup_{n}Y_{j}^{n,k}(t)$, it follows from
the stability property of quadratic BSDE with bounded terminal data
(see Lemma \ref{lemma_stability}) that there
exists $Z_{j}(\cdot)\in M^2$ such that
$\lim_{k\rightarrow\infty}\lim_{n\rightarrow\infty}Z_{j}^{n,k}(\cdot)=Z_{j}(\cdot)$
in $M^2$, and $\left(Y_{j}(\cdot),Z_{j}(\cdot)\right)$ satisfies
\begin{equation}\label{BSDEj*}
	Y_j(t)=F_{j}+\int_t^{\tau_{j}}f(s,Z_j(s))ds-\int_t^{\tau_{j}}
	\left(Z_j(s)\right)^{tr}dB_s,
\end{equation}
where $F_{j}=Y_{j}(T)=\inf_{k}\sup_{n}Y_{\tau_{j}}^{n,k}$.

We now let \( j \) tend to infinity in (\ref{BSDEj*}), following the procedure outlined in Section 4 of \cite{BH1} (see also Section 2 of \cite{BH2}).
 By construction, $\tau_{j}\leq \tau_{j+1}$. Hence, $Y^{n,k}_{j+1}(t\wedge \sigma_{j})=Y^{n,k}_{j}(t)$, and so we have the localization property:
$$Y_{j+1}(t\wedge \sigma_{j})=Y_{j}(t),\ \ \ Z_{j+1}(t\wedge \sigma_{j})=Z_{j}(t).$$
So if we set $\tau_0=0$ and define the processes $(Y,Z)$ by 
\begin{align*}
Y_t:=Y_{1}(0)+\sum_{j=1}^{\infty}Y_{j}(t)1_{(\tau_{j-1},\tau_{j}]}(t), \quad 
Z_t:=\sum_{j=1}^{\infty}Z_{j}(t)1_{(\tau_{j-1},\tau_{j}]}(t); \quad t\in [0,T], 
\end{align*}
we deduce that $Y$ is a
continuous process, $Z\in L^{2}(0,T)$, and
$$\lim_{j\rightarrow\infty} \sup_{t\in [0,T]} |Y_{j}(t)-Y_t|=0,\quad \lim_{j\rightarrow\infty} \int_0^T |Z_j(t)-Z_t|^2dt=0,\quad \mathbb{P}\text{-a.s.}$$
Moreover, $(Y,Z)$ satisfies BSDE$(F,f)$ (cf. (\ref{BSDEexpo})).

To finish the proof of existence, it remains to prove that the
solution $(Y,Z)$ stays in an appropriate space, which is established in the following lemma.
\end{proof}

\begin{lemma}\label{lemma0} The solution $(Y,Z)$ to BSDE$(F,f)$ constructed in Theorem \ref{theorem:BSDE1}
	satisfies
	\begin{equation}\label{solutionclassY}
		E\left[e^{p\alpha Y^{+}_{\star}}+e^{\varepsilon
			Y^{-}_{\star}}+\int_0^{T}|Z_s|^2ds\right]<+\infty,
	\end{equation}
	where $p>1$ and $\varepsilon>0$ are given in (\ref{condition}).
	Moreover, if $E[e^{p^{\prime} |F|}]<+\infty$ for any $p^{\prime}\ge
	1$, then
	\begin{equation}\label{solutionclassZ}
		E\left[e^{p^{\prime}\alpha Y_{\star}}+\left(\int_0^T
		|Z_s|^2ds\right)^{p^{\prime}/2}\right]<+\infty.
	\end{equation}
\end{lemma}

\begin{proof} Firstly, the construction of the limiting process $(Y,Z)$
	implies that $Z\in L^2(0,T)$, and moreover, by the inequality
	(\ref{bound}), it also holds that $\underline{Y}_t\leq Y_t\leq
	\overline{Y}_t$. Hence, using the upper bound of $Y_t$, we obtain
	$$e^{p\alpha Y_t^+}\leq \left(E[e^{\alpha F^+}|\mathcal{F}_t]\right)^{p}.$$
	In turn, Doob's inequality further yields
	\begin{align}\label{inequality11}
		E[e^{p\alpha Y^+_{\star}}]\leq
		E\left[\sup_{t\in[0,T]}\left(E[e^{\alpha
			F^+}|\mathcal{F}_t]\right)^{p}\right]\ 
		\leq\left(\frac{p}{p-1}\right)^{p}E[e^{p\alpha F^+}]<+\infty,
	\end{align}
	that is,  $e^{\alpha Y^+}\in\mathcal{S}^{p}$.
	
	To prove $e^{\varepsilon Y^-}\in\mathcal{S}^1$, let us fix
	$q^{\prime}>q>1$. Using $Y_t^{-}\leq -\underline{Y}_t$ and Jensen's
	inequality we deduce that
	\begin{align*}
		e^{\frac{\varepsilon}{q^{\prime}}Y_t^{-}}&\leq
		\exp\left(\frac{\varepsilon}{q^{\prime}}E^{\mathbb{Q}^{\theta}}\left[F^{-}+\int_t^T
		\frac{1}{2\alpha}|\theta_s|^2ds|{\mathcal F}_t\right]\right)\\
		&\leq
		E^{\mathbb{Q}^{\theta}}\left[\exp\left(\frac{\varepsilon}{q^{\prime}}(F^{-}+\int_t^T
		\frac{1}{2\alpha}|\theta_s|^2ds)\right)|\mathcal{F}_t\right].
	\end{align*}
	Then, H\"older's inequality and the boundedness of $\theta$ further
	yield
	\begin{align*}
		e^{\frac{\varepsilon}{q^{\prime}}Y_t^{-}}&\leq C
		\left(E\left[\exp\left(\frac{\varepsilon
			q}{q^{\prime}}(F^{-}+\int_t^T
		\frac{1}{2\alpha}|\theta_s|^2ds)\right)|\mathcal{F}_t\right]\right)^{\frac{1}{q}}\\
		&\leq C \left(E[e^{\frac{\varepsilon
				q}{q^{\prime}}F^-}|\mathcal{F}_t]\right)^{\frac{1}{q}}.
	\end{align*}
	In turn, it follows from Doob's inequality that
	\begin{align}\label{inequality111}
		E[e^{\varepsilon Y^{-}_{\star}}]&\leq CE\left[\sup_{t\in[0,T]}
		\left(E[e^{\frac{\varepsilon
				q}{q^{\prime}}F^-}|\mathcal{F}_t]\right)^{\frac{q^{\prime}}{q}}
		\right]\notag\\
		&\leq
		C\left(\frac{q^{\prime}}{q^{\prime}-q}\right)^{q^{\prime}/q}E[e^{\varepsilon
			F^-}]<+\infty.
	\end{align}
	
	Next, we show that $Z\in M^2$. To this end, following \cite{bdh}
	(see also \cite{BH1} and \cite{BH2}), we introduce
	$u_{\varepsilon}(x):=\frac{\exp(\varepsilon x)-\varepsilon
		x-1}{\varepsilon^2}$ for $x\geq 0$, and apply It\^o-Tanaka  formula
	to $u_{\varepsilon}(Y_{t}^{-})$. It follows that
	\begin{align}\label{Itoformula}
		&\int_0^{t\wedge\tau_{j}^{\varepsilon}}1_{\{Y_s\leq
			0\}}\left(\frac12
		u^{''}_{\varepsilon}(Y_s^{-})|Z_s|^2+u^{'}_{\varepsilon}(Y_s^{-})f(s,Z_s)\right)ds\notag\\
		=&\
		u_{\varepsilon}(Y_{t\wedge\tau_j^{\varepsilon}}^{-})-u_{\varepsilon}(Y_0^{-})
		-\int_0^{t\wedge\tau_{j}^{\varepsilon}}\frac12u_{\varepsilon}^{'}(Y_s^{-})dL_s\notag\\
		&+\int_0^{t\wedge\tau_{j}^{\varepsilon}}1_{\{Y_s\leq
			0\}}u_{\varepsilon}^{'}(Y_s^{-})Z_s^{tr}dB_s
	\end{align}
	where $L$ is the local time of $Y$ at the level $0$, and
	$\tau_{j}^{\varepsilon}$ is defined as
	$$\tau_j^{\varepsilon}=T\wedge \inf\left\{t\in [0,T]: \int_0^{t}|u_{\varepsilon}^{\prime}(Y_s^{-})Z_s|^2ds>j\right\}.$$

	Since $f(s,Z_s)\geq -Z_s^{tr}\theta_s-\frac{1}{2\alpha}|\theta_s|^2$
	and $u_{\varepsilon}^{\prime}(Y_s^{-})\geq 0$, we have
	\begin{equation*}
		u_{\varepsilon}^{\prime}(Y_s^{-})f(s,Z_s)
		\geq-u^{'}_{\varepsilon}(Y_s^{-})\frac{\varepsilon}{2}|Z_s|^2-
		u_{\varepsilon}^{'}(Y_s^{-})\left(\frac{1}{2\varepsilon}+\frac{1}{2\alpha}\right)|\theta_s|^2.
	\end{equation*}
	With the choice of the function $u_{\varepsilon}(\cdot)$, we also
	have
	$$u_{\varepsilon}(Y_{t\wedge\tau_j}^{-})\leq \frac{1}{\varepsilon^2}\, {e^{\varepsilon
			Y_{t\wedge\tau_j}^{-}}}\leq \frac{1}{\varepsilon^2}\, {e^{\varepsilon
			Y^{-}_{\star}}},\quad -u_{\varepsilon}(Y_0^{-})\leq
	0, \quad \text{and} \quad 
	-\int_0^{t\wedge\tau_j}u_{\varepsilon}^{'}(Y_s^{-})dL_s\leq 0.$$
	Hence, (\ref{Itoformula}) further yields
	\begin{align*}
		&\int_0^{t\wedge\tau_j^{\varepsilon}}1_{\{Y_s\leq 0\}}\frac
		{u_{\varepsilon}^{''}(Y_s^-)-\varepsilon u_{\varepsilon}^{'}(Y_s^-)}{2}|Z_s|^2ds\\[3mm]
		\leq&\ \frac{e^{\varepsilon
				Y_{\star}^-}}{\varepsilon^2}+\int_0^{T}1_{\{Y_s\leq
			0\}}u_{\varepsilon}^{'}(Y_s^-)(\frac{1}{2\varepsilon}+\frac{1}{2\alpha})|\theta_s|^2ds\\[3mm]
		&\ +\int_0^{t\wedge\tau_j^{\varepsilon}}1_{\{Y_s\leq
			0\}}u_{\varepsilon}^{'}(Y_s^-)Z_s^{tr}dB_s.
	\end{align*}
	Since 
	$$u^{''}_{\varepsilon}(Y_s^-)-\varepsilon
	u_{\varepsilon}^{'}(Y_s^-)=1\quad \text{and}\quad u_{\varepsilon}^{'}(Y_s^-)\leq
	\frac{e^{\varepsilon Y_s^-}}{\varepsilon}\leq\frac{e^{\varepsilon
			Y_{\star}^-}}{\varepsilon},$$
		 we have
	\begin{align}\label{inequality22}
		&\ \frac12\int_0^{t\wedge\tau_j^{\varepsilon}}1_{\{Y_s\leq 0\}}|Z_s|^2ds\notag\\
		\leq&\ (\frac{1}{\varepsilon^2}+C)e^{\varepsilon
			Y_{\star}^{-}}+\int_0^{t\wedge\tau_j^{\varepsilon}}1_{\{Y_s\leq
			0\}}u_{\varepsilon}^{'}(Y_s^-)Z_s^{tr}dB_s,
	\end{align}
	and therefore, $E[\int_0^T1_{\{Y_s\leq 0\}}|Z_s|^2ds]<+\infty$.
	
	Similarly, applying It\^o-Tanaka  formula to
	$u_{\alpha}(Y_{t}^{+})$, since $f(s,Z_s)\leq
	\frac{\alpha}{2}|Z_s|^2$, we obtain
	\begin{align}\label{inequality2222}
		\ \frac12\int_0^{t\wedge\tau_j^{\alpha}}1_{\{Y_s> 0\}}|Z_s|^2ds\ 
		\leq\ \frac{1}{\alpha^2}e^{\alpha
			Y_{\star}^{+}}-\int_0^{t\wedge\tau_j^{\alpha}}1_{\{Y_s>
			0\}}u_{\alpha}^{'}(Y_s^+)Z_s^{tr}dB_s,
	\end{align}
	and in turn, $E[\int_0^T1_{\{Y_s> 0\}}|Z_s|^2ds]<+\infty$.
	
	Next, we show that $e^{\alpha Y}\in\mathcal{S}^{p^{\prime}}$ and
	$Z\in M^{p^{\prime}}$ for any $p^{\prime}\geq 1$, when $F$ has
	exponential moment of any order. The first part about $Y$ follows
	along the same arguments as in (\ref{inequality11}) and
	(\ref{inequality111}), so we only prove the second part about $Z$.
	To this end, we send $j\rightarrow\infty$ in (\ref{inequality22})
	and (\ref{inequality2222}), and add the two inequalities,
	\begin{align*}
		\frac{1}{2}\int_0^t|Z_s|^2ds\leq&\
		(\frac{1}{\varepsilon^2}+C)e^{\varepsilon
			Y_{\star}^{-}}+\frac{1}{\alpha^2}e^{\alpha Y_{\star}^{+}}\\
		&\ +\sup_{t\in[0,T]}\left|\int_0^t\left(1_{\{Y_s\leq
			0\}}u_{\varepsilon}^{'}(Y_s^-)-1_{\{Y_s>
			0\}}u_{\alpha}^{'}(Y_s^+)\right)Z_s^{tr}dB_s\right|.
	\end{align*}
	Consequently, there exists a constant $K>0$ such that
	\begin{align*}
		&\ E\left[\left(\int_0^T|Z_s|^2ds\right)^{p^{\prime}/2}\right]\leq
		K\, E\left[e^{\frac{p^{\prime}\varepsilon}{2}Y_{\star}^-}\right]+K\, E\left[e^{\frac{p^{\prime}\alpha}{2}Y_{\star}^+}\right]\\
		&\qquad\qquad +K E\left[\sup_{t\in[0,T]}\left|\int_0^{t} \left(1_{\{Y_s\leq
			0\}}u_{\varepsilon}^{'}(Y_s^-)-1_{\{Y_s>
			0\}}u_{\alpha}^{'}(Y_s^+)\right)
		Z_s^{tr}dB_s\right|^{p^{\prime}/2}\right]
	\end{align*}
	for any $p^{\prime}\geq 1$.
	
	Next, we apply the B-D-G inequality to the last term in the above
	inequality, and obtain (with constant $C$ varying from line to line)
	\begin{align*}
		&KE\left[\sup_{t\in[0,T]}\left|\int_0^{t} \left(1_{\{Y_s\leq
			0\}}u_{\varepsilon}^{'}(Y_s^-)-1_{\{Y_s>
			0\}}u_{\alpha}^{'}(Y_s^+)\right)
		Z_s^{tr}dB_s\right|^{p^{\prime}/2}\right]\\
		\leq &\ CE\left[\left(\int_0^T(|u^{'}_{\varepsilon}(Y_s^-)|^2+|u^{'}_{\alpha}(Y_s^+)|^2)|Z_s|^2ds\right)^{p^{\prime}/4}\right]\\
		\leq &\
		CE\left[\left(e^{\frac{p^{\prime}\varepsilon}{2}Y_{\star}^-}+
		e^{\frac{p^{\prime}\alpha}{2}Y_{\star}^+}
		\right)\left(\int_0^T|Z_s|^2ds\right)^{p^{\prime}/4}\right]\\
		\leq &\
		\frac12E\left[\left(\int_0^T|Z_s|^2ds\right)^{p^{\prime}/2}\right]+CE\left [e^{p^{\prime}\varepsilon
			Y_{\star}^-}+e^{p^{\prime}\alpha Y_{\star}^+}\right ].
	\end{align*}
	Hence, $E[(\int_0^T|Z_s|^2ds)^{p^{\prime}/2}]< +\infty$, and we
	conclude the proof.
\end{proof}

Finally, the uniqueness of the solution to BSDE (\ref{BSDEexpo})
will be established in the proof of Theorem \ref{theorem_dual} in Chapter 5
as a corollary of the convex dual representation of the solution component
$Y$.

%
%
%


\chapter{Utility maximization with unbounded payoffs}\label{section:main}

\section{Review of the case of bounded payoffs}\label{section:review}

We first review the case of bounded payoffs, which was studied in Hu et al.~\cite{him}. For its application to utility indifference valuation,  see~\cite{MS}; for an extension to continuous martingales,  see~\cite{Morlais}; and for an extension to jumps, see~\cite{Bec}. The problem is directly solved by utilizing the \emph{martingale optimality principle}. In essence, the conditional value process is a supermartingale for any admissible trading strategy and is a martingale for an optimal trading strategy. It is important to note that duality does not play a role in this context, which explains why the constraint set
$\mathcal{C}$ does not need to be convex when $F$ is bounded. The following result is available in Theorem 7 of \cite{him}. We include it for completeness.

\begin{lemma}[\cite{him}]\label{lemma_boundedcase} Suppose that the payoff $F$ is bounded. Let $(Y,Z)$ be the unique solution
	to BSDE$(F, f)$ with $f$ given by (\ref{driver}). Then, the value function of
	the optimization problem (\ref{expoopt}) with admissible set
	$\mathcal{A}_D$ is given by
	\begin{equation}\label{lemmavaluefunction}
		V(0, x)= -\exp(-\alpha(x-Y_0)),
	\end{equation}
	and there exists an optimal trading strategy $\pi^{\star} \in {\mathcal
		A}_D$ such that 
	\begin{equation}  \label{lemmapi}
\sigma_t^{tr}\pi^{\star}_t=\mbox{Proj}_{\sigma_t^{tr}\mathcal{C}}\left(Z_t+\frac{\theta_t}{\alpha}\right),\
		\mbox{$\mathbb{P}$-a.s.},\ \text{for}\ t \in [0,T].	
	\end{equation}
\end{lemma}

\begin{proof} It suffices
	to prove that the conditional value process has the form
	$V(t,X_t^{\pi})=-e^{-\alpha (X_t^{\pi}-Y_t)}$, $t\in[0,T]$, which is
	a supermartingale for any $\pi\in\mathcal{A}_{D}$, and is a
	martingale for $\pi^*$ given in (\ref{lemmapi}), with
	$\pi^{\star}\in\mathcal{A}_D$. Consequently,
\begin{equation*}
-\exp(-\alpha(x-Y_0))=\sup_{\pi\in\mathcal{A}_{D}}E[-\exp(-\alpha(X_T^{\pi}-F))]=E[-\exp(-\alpha(X_T^{\pi^{\star}}-F))].
\end{equation*}
	
	Firstly, for $\pi\in\mathcal{A}_{D}$, an application of
	It\^o's formula to $-e^{-\alpha(X_{t}^{\pi}-Y_{t})}$ gives
	$$-e^{-\alpha(X_{t}^{\pi}-Y_{t})}=-e^{-\alpha(x-Y_0)}A_t^{\pi}L_t^{\pi},$$
	where
	$$A_t^{\pi}=\exp\left(\alpha\int_0^t
	\left(\frac{\alpha}{2}|\sigma_u^{tr}\pi_u-(Z_u +
	\frac{\theta_u}{\alpha})|^2 - Z_u^{tr}\theta_u - \frac{|\theta_u|^2}{2
		\alpha}-f(u,Z_u)\right)du
	\right),$$ and
	\begin{equation*}
		L_{t}^{\pi}=\mathcal{E}_t\left(\alpha\int_0^{\cdot}(Z_u^{tr}-\pi_u^{tr}\sigma_u)dB_u\right).
	\end{equation*}
 Since $A_t^{\pi}$ is increasing in $t$, it follows that
	$A_t^{\pi}L_t^{\pi}$, $t\in[0,T]$, is a
	local submartingale. Hence, there exist an increasing sequence of stopping times $\tau_j\uparrow T$, such that
    \begin{equation}\label{inequality_local}
    E[A_{s\wedge\tau_j}^{\pi}L_{s\wedge\tau_j}^{\pi}|\mathcal{F}_t]\geq A_{t\wedge\tau_j}^{\pi}L_{t\wedge\tau_j}^{\pi},
    \end{equation}
    for $s\geq t$.
From the definition of $\mathcal{A}_{D}$ and the boundedness of $Y$, we know that $-e^{-\alpha(X_{t}^{\pi}-Y_{t})}$, $t\in[0,T]$, is uniformly integrable. The supermartingale property of $-e^{-\alpha(X_{t}^{\pi}-Y_{t})}$, $t\in[0,T]$, then follows by taking $j\uparrow \infty$ in (\ref{inequality_local}).
	
	To prove the martingale property of
	$-e^{-\alpha(X_{t}^{\pi^{\star}}-Y_{t})}$, $t\in[0,T]$, we observe
	that if $\pi^{\star}$ is a minimizer in (\ref{pi}), then
	$A_t^{\pi^{\star}}=1$. On the other hand,
$\int_0^{t}(Z_u^{tr}-(\pi_u^{\star})^{tr}\sigma_u)dB_u$, $t\in[0,T]$, is a BMO martingale by Remark \ref{remark_BMO}, which in turn implies the uniform integrability of $L^{\pi^{\star}}$ (see \cite{Kaz}). Since, moreover, $Y$ is bounded, we obtain the class (D) property of $e^{-\alpha X^{\pi^{\star}}}$.
\end{proof}

\section{Unbounded payoffs of an exponential integrability}

We are now ready to provide one of the main results in this section, which is
the characterization of the value function and the corresponding
optimal trading strategy for the optimization problem
(\ref{expoopt}) under Assumption \ref{assumption1}.

\begin{theorem} \label{thm:opti2}
	Suppose that Assumption 1 holds. Let $(Y,Z)$ be the unique solution
	to BSDE$(F, f)$ with $f$ given by (\ref{driver}). Then, the value function of
	the optimization problem (\ref{expoopt}) with admissible set
	$\mathcal{A}_D^{conv}$ is given by
	\begin{equation}\label{valuefunction}
		V(0, x)= -\exp(-\alpha(x-Y_0)),
	\end{equation}
	and there exists an optimal trading strategy $\pi^* \in {\mathcal
		A}_D^{conv}$ such that
	\begin{equation}  \label{pi}
\sigma_t^{tr}\pi^{\star}_t=\mbox{Proj}_{\sigma_t^{tr}\mathcal{C}}\left(Z_t+\frac{\theta_t}{\alpha}\right),\
		\mbox{$\mathbb{P}$-a.s.},\ \text{for}\ t \in [0,T].
	\end{equation}
\end{theorem}

\begin{proof} The first part of the proof mirrors the proof for Lemma \ref{lemma_boundedcase}. However, given that $Y$ is unbounded, ensuring uniform integrability requires entirely different arguments, as presented below.

Note that the Class (D)
	condition in $\mathcal{A}_D^{conv}$ is equivalent to say that
	the conditional value process $-e^{-\alpha(X_{t}^{\pi}-Y_{t})}$, $t\in[0,T]$, 
	is in Class (D) and, therefore, $A^{\pi}L^{\pi}$ is in
	Class (D). Since $A_t^{\pi}$ is increasing in $t$, it follows that
	$-e^{-\alpha(X_{t}^{\pi}-Y_{t})}$, $t\in[0,T]$, is a true
	supermartingale.
	
	To prove the martingale property of
	$-e^{-\alpha(X_{t}^{\pi^{\star}}-Y_{t})}$, $t\in[0,T]$, we observe
	that if $\pi^{\star}$ is a minimizer in (\ref{pi}), then
	$A_t^{\pi^{\star}}=1$. Hence, it is sufficient to prove that, with $\pi^{\star}$ given in
	(\ref{pi_2}), the optimal density process $L^{\pi^{\star}}$ is in
	Class (D), which will further imply that $L^{\pi^{\star}}$ is a
	uniformly integrable martingale, and
	$-e^{-\alpha(X_{t}^{\pi^{\star}}-Y_{t})}$, $t\in[0,T]$, is also in Class
	(D).

We complete the proof by demonstrating that $L_T^{\pi^{\star}}$ has finite entropy, as detailed in Lemma \ref{lemma}. Subsequently, by invoking the De la Vall\'ee-Poussin theorem (See Chapter 2 in \cite{Meyer1978} by Dellacherie and Meyer), we establish its uniform integrability. By (\ref{subdifferential}) in Lemma \ref{lemma_convex}, $\alpha(Z_t-\sigma_t^{tr}\pi_t^{\star})\in\partial_z f(t,Z_t)$.  Hence,  it is sufficient to prove that $L_T^{q^{\star}}$ has finite entropy for any $q^{\star}_t\in\partial_z f(t,Z_t)$, where
\begin{equation*}		L_{t}^{q^{\star}}=\mathcal{E}_t\left(\int_0^{\cdot}q_u^{\star}dB_u\right),\quad t\in[0,T],
	\end{equation*}
which is verified in the following lemma, concluding the proof.
\end{proof}

\begin{lemma}\label{lemma}
	The optimal density process $L_T^{q^{\star}}$ has finite entropy.
	Hence, by De la Vall\'ee-Poussin theorem, $L^{q^{\star}}$ is in
	Class (D) and, therefore, it is a uniformly integrable martingale.
\end{lemma}

\begin{proof}

To show that
	$L^{q^{\star}}_T=\mathcal{E}_T(\int_0^{\cdot}(q_u^{\star})^{tr}dB_u)$
	has a  finite entropy, we introduce
	the following stopping times: for any integer $j\ge 1$,
	$$\sigma_j=T\wedge \inf\left\{t\in [0,T]: \max\left\{\int_0^t|{q_s^{\star}}|^2ds,\int_0^t|Z_s|^2ds\right\} >j\right\},$$
	so that $L^{q^{\star}}_{\cdot\wedge\sigma_j}$ is a uniformly
	integrable martingale under $\mathbb{P}$.
	
	We then define a probability measure $\mathbb{Q}^{q^{\star}}$ on
	$\mathcal{F}_{\sigma_j}$ by
	$d\mathbb{Q}^{q^{\star}}:=L_{\sigma_j}^{q^{\star}} d\mathbb{P} $,
	and an ${m}$-dimensional Brownian motion
	$B^{q^{\star}}_t:=B_t-\int_0^tq_u^{\star}du$, for
	$t\in[0,\sigma_j]$. It follows from the definitions of $\sigma_j$
	and $B^{{q}^{\star}}$ that both processes 
	$$\int_0^{\cdot\wedge\sigma_j}({q^{\star}})^{tr}dB_s^{q^{\star}} \quad 
	\text{and} \quad  \int_0^{\cdot\wedge\sigma_j}(Z_u)^{tr}dB_s^{q^{\star}}$$ are  $\mathbb{Q}^{q^{\star}}$-martingales.
	
Note that $x(\ln x-1)$ and $e^x$ are convex duals of each other. Applying the Fenchel inequality
	\begin{equation}\label{Young_inequ}
		xy=(\frac{x}{p})(p y)\leq \frac{x\ln x}{p}-\frac{x\ln p}{p}+e^{p
			y},\quad \mbox{for} \ (x,y)\in (0,\infty)\times \mathbb{R},
	\end{equation}
	to
	$L_{\sigma_{j}}^{q^{\star}}(\alpha Y_{\sigma_j})$ gives
	\begin{equation}\label{inequality4}
		E\left [L_{\sigma_{j}}^{q^{\star}}\alpha Y_{\sigma_j}\right]\leq
		\frac{1}{p}E\left [L_{\sigma_j}^{q^{\star}}\ln L_{\sigma_j}^{q^{\star}}\right]-
		\frac{1}{p}E\left [L_{\sigma_j}^{q^{\star}}\right]\ln p+E\left[e^{p\alpha
			Y^+_{\star}}\right]
	\end{equation}
	with $p>1$ given in (\ref{condition}).
	
	Furthermore, writing BSDE
	(\ref{BSDEexpo}) under $\mathbb{Q}^{q^{\star}}$, we obtain
	\begin{align}\label{inequality5}
		E[L_{\sigma_{j}}^{q^{\star}}\alpha
		Y_{\sigma_j}] & =E^{\mathbb{Q}^{q^{\star}}}[\alpha
		Y_{\sigma_{j}}]\notag\\
		&=E^{\mathbb{Q}^{q^{\star}}}\left[\alpha
		Y_0+\alpha\int_0^{\sigma_{j}}
(-f(u,Z_u)+Z_u^{tr}q_u^{\star})du+\alpha\int_0^{\sigma_j}(Z_u)^{tr}dB_u^{q^{\star}}\right]\notag\\
        &=\alpha Y_0+E^{\mathbb{Q}^{q^{\star}}}\left[\alpha\int_0^{\sigma_j}f^{\star}(u,q_u^{\star})du\right]\notag\\
		&\geq \alpha Y_0+E^{\mathbb{Q}^{q^{\star}}}\left[
		\int_0^{\sigma_{j}}\frac12|q_u^{\star}|^2du\right],
	\end{align}
	where we used the lower bound of $f^{\star}$ (cf.
	(\ref{lowerbound})) in the last inequality.
	
	Finally, combining (\ref{inequality4}) and (\ref{inequality5}), and
	observing
\begin{align*}	
E[L_{\sigma_j}^{q^{\star}}\ln
L_{\sigma_j}^{q^{\star}}]&=E^{\mathbb{Q}^{q^{\star}}}\left[\int_0^{\sigma_{j}}(q_u^{\star})^{tr}dB_u^{q^{\star}}+\frac12|q_u^{\star}|^2du\right]\\
&=E^{\mathbb{Q}^{q^{\star}}}\left[\int_0^{\sigma_{j}}\frac12|q_u^{\star}|^2du\right],
\end{align*}
	we obtain
	$$(1-\frac{1}{p})E\left[L_{\sigma_j}^{q^{\star}}\ln
	L_{\sigma_j}^{q^{\star}}\right]\leq -\alpha Y_0-\frac{\ln
		p}{p}+E\left [e^{p\alpha Y^+_{\star}}\right]<+\infty.$$ The assertion then
	follows by sending $j\rightarrow\infty$ in the above inequality, and
	using Fatou's lemma.
\end{proof}

\begin{remark}
We can offer a more direct proof without relying on the subdifferential $\partial_z f(t, Z_t)$, by instead utilizing the projection property onto a closed and convex set $\Pi \subset \mathbb{R}^m$. Specifically, for any $x \in \mathbb{R}^m$, the projection operator on a closed and convex set $\Pi$ satisfies:
\begin{equation}\label{projection_operator}
(\mbox{Proj}_\Pi(x) - x)^{tr} (\mbox{Proj}_\Pi(x) - y) \leq 0, \quad \text{for all } y \in \Pi.
\end{equation}

With this, the lower bound in (\ref{inequality5}) can be replaced as follows:
\begin{align*}
E[L_{\tau_{j}}^{\pi^{\star}} \alpha Y_{\tau_j}] &= E^{\mathbb{Q}^{\pi^{\star}}}[\alpha Y_{\tau_j}] \\
&= E^{\mathbb{Q}^{\pi^{\star}}} \left[ \alpha Y_0 - \int_0^{\tau_j} \left( \alpha f(u,Z_u) - \alpha^2 Z_u^{tr} (Z_u - \sigma_u^{tr} \pi_u^{\star}) \right) du \right] \\
&\geq \alpha Y_0 + E^{\mathbb{Q}^{\pi^{\star}}} \left[ \int_0^{\tau_j} \frac{\alpha^2}{2} |Z_u - \sigma_u^{tr} \pi_u^{\star}|^2 du \right] \\
&\quad - \alpha^2 E^{\mathbb{Q}^{\pi^{\star}}} \left[ \int_0^{\tau_j} \left( |Z_u - \sigma_u^{tr} \pi_u^{\star}|^2 - \frac{\theta_u^{tr}}{\alpha} \sigma_u^{tr} \pi_u^{\star} - Z_u^{tr} (Z_u - \sigma_u^{tr} \pi_u^{\star}) \right) du \right] \\
&= \alpha Y_0 + E^{\mathbb{Q}^{\pi^{\star}}} \left[ \int_0^{\tau_j} \frac{\alpha^2}{2} |Z_u - \sigma_u^{tr} \pi_u^{\star}|^2 du \right] \\
&\quad - \alpha^2 E^{\mathbb{Q}^{\pi^{\star}}} \left[ \int_0^{\tau_j} \left( (\pi_u^{\star})^{tr} \sigma_u - (Z_u^{tr} + \frac{\theta_u^{tr}}{\alpha}) \right) \sigma_u^{tr} \pi_u^{\star} du \right].
\end{align*}
Here, $\mathbb{Q}^{\pi^{\star}}$ is the probability measure defined by $\frac{d\mathbb{Q}^{\pi^{\star}}}{d\mathbb{P}} := L_{\tau_j}^{\pi^{\star}}$ for a localization sequence $(\tau^j)_{j\geq 1}$. Note that the above inequality follows from the expression of $f$ in (\ref{driver}):
\[
f(u, Z_u) \leq \frac{\alpha}{2} |Z_u - \sigma_u^{tr} \pi_u^{\star}|^2 - \theta_u^{tr} \sigma_u^{tr} \pi_u^{\star}.
\]
Then, by invoking (\ref{projection_operator}) and observing that $0 \in \sigma_u^{tr} \mathcal{C}$, we obtain:
\begin{align*}
&\left( (\pi_u^{\star})^{tr} \sigma_u - \left( Z_u^{tr} + \frac{\theta_u^{tr}}{\alpha} \right) \right) \sigma_u^{tr} \pi_u^{\star} \\
&= \left( \mbox{Proj}_{\sigma_u^{tr} \mathcal{C}} \left( Z_u + \frac{\theta_u}{\alpha} \right) - \left( Z_u + \frac{\theta_u}{\alpha} \right) \right)^{tr} \left( \mbox{Proj}_{\sigma_u^{tr} \mathcal{C}} \left( Z_u + \frac{\theta_u}{\alpha} \right) - 0 \right) \leq 0.
\end{align*}
Consequently, (\ref{inequality5}) can now be replaced by the following inequality:
\[
E[L_{\tau_{j}}^{\pi^{\star}} \alpha Y_{\tau_j}] \geq \alpha Y_0 + E^{\mathbb{Q}^{\pi^{\star}}} \left[ \int_0^{\tau_j} \frac{\alpha^2}{2} |Z_u - \sigma_u^{tr} \pi_u^{\star}|^2 du \right].
\]
The remainder of the proof follows similarly to the proof of Lemma \ref{lemma}.
\end{remark}

\section{Exponentially integrable payoffs bounded from below}

Our next result is about relaxing the Class (D) condition in the
admissible set $\mathcal{A}_D$ when $F$ satisfies a stronger
condition than Assumption \ref{assumption1}.

\begin{assumption}\label{assumption2} The payoff $F$ satisfies
	$E[e^{p\alpha F^+}]<+\infty$ for some integer $p>1$, and there
	exists a constant $k>0$ such that $F^{-}\leq k$.
\end{assumption}

\begin{remark}
	The boundedness from below on $F$ means that there is a uniform
	lower bound on the amount one can lose from the random endowment
	$F$, if it is negative. Very often, $F$ models a payoff, so it is
	even nonnegative.
	
	A similar type of assumption on $F$ is also imposed in Delbaen et al \cite{dgr},
	where the authors establish the minimal entropy representation as
	the dual of the utility maximization problem (\ref{expoopt}) in the
	case of subspace portfolio constraints.
\end{remark}

Since $\theta$ is bounded, we can define an equivalent minimal local
martingale measure (MLMM) $\mathbb{Q}^{\theta}$ on $\mathcal{F}_T$
by
\begin{equation}\label{EMM}
	\frac{d\mathbb{Q}^{\theta}}{d\mathbb P}:=L_T^{\theta}={\mathcal
		E}_T(-\int_0^{\cdot}\theta_t^{tr}dB_t),
\end{equation}
where $\mathcal{E}(\cdot)$ denotes the stochastic exponential.
Similarly, we define $\frac{d\mathbb{Q}^{\theta}}{d\mathbb
	P}|_{\mathcal{F}_t}:=L_t^{\theta}$, for $t\in[0,T]$. Then, under
$\mathbb{Q}^{\theta}$, the wealth process $X^{\pi}$ follows
\begin{equation}\label{wealth_Q}
	X^{\pi}_t=x + \int_0^t \pi_u^{tr}
	\sigma_ud{B}_u^{\theta},
\end{equation}
where ${B}_t^{\theta}:=B_t+\int_0^{t}\theta_udu$, $t\in[0,T]$, is an
$m$-dimensional Brownian motion under MLMM $\mathbb{Q}^{\theta}$.

In the following, we replace the Class (D) condition in the
admissible set $\mathcal{A}_D$ by an equivalent minimal martingale
measure (MMM) condition, and solve the optimization problem
(\ref{expoopt}) under Assumption \ref{assumption2}.

\begin{definition}\label{admiss3}
	[{\bf Admissible Strategies with constraints $\mathcal{A}_D^{M}$}]
	
	Let ${\mathcal C}$ be a closed and convex set in $\mathbb R^{d}$ satisfying $0\in\mathcal{C}$. The set of
	admissible trading strategies $\mathcal{A}_D^{M}$ consists of all
	$\mathbb{R}^d$-valued predictable processes $\pi\in L^2[0,T]$, which
	are self-financing and satisfy $\pi_t \in {\mathcal C}$, $\mathbb
	P$-a.s., for $t\in[0,T]$. Moreover, $(X^{\pi}_t)_{t\ge 0}$ is a
	martingale under $\mathbb{Q}^{\theta}$, that is,
	$\mathbb{Q}^{\theta}$ is an MMM.
\end{definition}

Herein, we use the superscript $^M$ to highlight the dependency
on the minimal martingale measure.
We are now in a position to present the final main result in this
section. In particular, when $F$ is bounded, our result will also
generalize Theorem 7 of Hu et al.~\cite{him} by enlarging its admissible set
from $\mathcal{A}_{D}$ to $\mathcal{A}_D^{M}$.

\begin{theorem} \label{thm:opti}
	Suppose that Assumption \ref{assumption2} holds. Let $(Y,Z)$ be the
	unique solution of BSDE$(F,f)$ with $f$ given by (\ref{driver}). Then, the value function
	$V(0,x)$ of the optimization problem (\ref{expoopt}) with admissible
	set $\mathcal{A}_D^M$ and the associated optimal trading strategy
	$\pi^*\in\mathcal{A}_D^M$ are given, respectively, as in
	(\ref{valuefunction}) and ($\ref{pi}$).
\end{theorem}


We first identify the space of solutions to the
quadratic BSDE$(F,f)$ when the terminal condition $F$ satisfies
Assumption \ref{assumption2}.

\begin{lemma}\label{lemma2} Suppose that Assumption \ref{assumption2} holds. Then
	BSDE$(F,f)$ admits a unique solution $(Y,Z)$,
	where $e^{\alpha Y^{+}}\in\mathcal{S}^{p}$,
	$Y^{-}\in\mathcal{S}^{\infty}$, and $Z\in M^{2}$.
\end{lemma}

\begin{proof} We show $Y^{-}\in\mathcal{S}^{\infty}$. Indeed, from the proof
	of Theorem \ref{theorem:BSDE1}, we have
	$$Y_t\geq \underline{Y_t}=-E^{\mathbb{Q}^{\theta}}\left[F^{-}+\int_t^T \frac{1}{2\alpha}|\theta_s|^2ds|{\mathcal F}_t\right]\geq -k-E^{\mathbb{Q}^{\theta}}\left[\int_0^{T}\frac{1}{2\alpha}|\theta_s|^2ds\right].$$
	The rest assertions have been proved in Lemma \ref{lemma0}.
\end{proof}

\textbf{Proof of Theorem \ref{thm:opti}.} For integer $n>0$, we
truncate the payoff $F$ as $F^n:=F\wedge n$, so that $-k\leq
F^{n}\leq n$. We first show that for any $\pi\in\mathcal{A}_D^M$, it
holds that
\begin{equation}\label{inequality7}
	E\left[-e^{-\alpha(X_T^{\pi}-F^n)}\right]\leq -e^{-\alpha(x-Y_0^n)},
\end{equation}
where $Y^n$ is (the first component of) the bounded solution to the
corresponding truncated quadratic BSDE$(F^n,f)$.

Note that if $E\left[e^{-\alpha X^\pi_T}\right]=+\infty$, then due
to the boundedness of $F^n$,
$$E\left[-e^{-\alpha(X_T^{\pi}-F^n)}\right]=-\infty.$$ Hence, without
loss of generality, we may assume that $E\left[e^{-\alpha
	X^\pi_T}\right]<+\infty$.

According to Theorem 7 of \cite{him}, it is clear that the
inequality (\ref{inequality7}) holds for
$\pi\in\mathcal{A}_D\subset\mathcal{A}_D^M$. Hence, to show
(\ref{inequality7}), it suffices to prove that there exists
$\pi^j\in \mathcal{A}_D$, such that
\begin{equation}\label{convergence}
	E\left[e^{-\alpha(X^{\pi^j}_T-F^n)}\right]\rightarrow
	E\left[e^{-\alpha(X^{\pi}_T-F^n)}\right], \quad \mbox{as}\
	j\rightarrow\infty, 
\end{equation}
for any $\pi\in\mathcal{A}_D^{M}$.

To this end, we define $\pi^j_t=\pi_t 1_{\{t\le \tau_j\}}$, for
$t\in[0,T]$, and integer $j\geq 1$, where $\tau_j$ is the stopping
time defined as
$$\tau_j=T\wedge\inf\{t\in[0,T]: X_t^{\pi}\le -j\}.$$
By the definition of $\tau_j$, $X^{\pi}_{\cdot\wedge\tau_j}$ is
bounded from below, so $e^{-\alpha X^{\pi}_{\cdot\wedge\tau_j}}$ is
bounded and therefore in Class (D), which means
$\pi_j\in\mathcal{A}_D$.

It remains to show the convergence in (\ref{convergence}). Note that
with $\pi^j$ defined as above, $$e^{-\alpha X^{\pi^j}_T}=e^{-\alpha
	X^{\pi}_{\tau_j}}\rightarrow e^{-\alpha X^{\pi}_T},\quad
\mbox{$\mathbb{P}$-a.s.},$$ so we only need to establish the
uniformly integrability of $e^{-\alpha X^{\pi}_{\tau_j}}$ under
$\mathbb{P}$.

We recall that, from a usual truncation argument, for any $\xi\in
L^{2}(\mathbb P)$, $E^{\mathbb{Q}^{\theta}}[\xi|{\mathcal F}_{\tau_j}]$
converges to $\xi$ in $L^{2}(\mathbb P)$, where
$\mathbb{Q}^{\theta}$ is the MLMM given in (\ref{EMM}), and
$L^{2}(\mathbb{P})$ denotes the space of square integrable random
variables under $\mathbb{P}$.

By our assumption on $\pi$, $E[e^{-\alpha X^\pi_T}]<+\infty$, so
$E^{\mathbb{Q}^{\theta}}[e^{-\frac{\alpha}{2} X_T^\pi}|{\mathcal
	F}_{\tau_j}]$ converges to $e^{-\frac{\alpha}{2} X_T^\pi}$ in
$L^{2}(\mathbb P)$. As a consequence, $\left(E^{\mathbb{Q}^{\theta}}
[e^{-\frac{\alpha}{2} X_T^\pi}|{\mathcal F}_{\tau_j}]\right)^{2}$ is
uniformly integrable under $\mathbb{P}$.

On the other hand, from Jensen's inequality and the fact that
$X^\pi$ is a martingale under $\mathbb{Q}^{\theta}$ (by the
definition of admissibility), we have
$$e^{-\frac{\alpha}{2} X_{\tau_j}^\pi}=e^{E^{\mathbb{Q}^{\theta}}[-\frac{\alpha}{2}X_T^{\pi}|\mathcal{F}_{{\tau}_j}]}\le E^{\mathbb{Q}^{\theta}}\left[e^{-\frac{\alpha}{2} X_T^\pi}|{\mathcal F}_{\tau_j}\right].$$
Thus, $e^{-\alpha X^\pi_{\tau_j}}$ is uniformly integrable under
$\mathbb {P}$, and the inequality (\ref{inequality7}) holds.

Since $F_n\leq F$, it follows from (\ref{inequality7}) that
$$
E\left[-e^{-\alpha(X^\pi_T-F)}\right]\le -e^{-\alpha(x-Y_0^n)}.
$$
Furthermore, since $Y_0^{n}\rightarrow Y_0$ (see the proof of
Theorem \ref{BSDEexpo}), sending $n\rightarrow\infty$ in the above
inequality yields
\begin{equation}\label{inequality1}
	E\left[-e^{-\alpha(X^\pi_T-F)}\right]\le -e^{-\alpha(x-Y_0)}.
\end{equation}

To prove the equality, note that the choice of the optimal trading
strategy $\pi^{\star}$ in (\ref{pi_2}) ensures that
$e^{-\alpha(X_t^{\pi^{\star}}-Y_t)}$, $t\in[0,T]$, is a positive
local martingale, hence a supermartingale, which implies
\begin{equation}\label{inequality2}
	E\left[-e^{-\alpha(X^{\pi^{\star}}_T-F)}\right]\geq \mathbb
	-e^{-\alpha(x-Y_0 )}.
\end{equation}

Finally, since $|\sigma_t^{tr}\pi_t^{\star}|\le
|Z_t|+\frac{|\theta_t|}{\alpha}$ and $Z\in M^{2}$, we have
$\sigma^{tr}\pi^{\star}\in M^{2}$. The B-D-G inequality then implies
that
\begin{align*}
	E^{\mathbb{Q}^{\theta}}\left[\sup_{t\in[0,T]}\left|\int_0^t(\sigma^{tr}_t\pi_t^{\star})^{tr}dB_t^{\theta}\right|\right]
	&\leq
	CE^{\mathbb{Q}^{\theta}}\left[\left(\int_0^T|\sigma_t^{tr}\pi_t^{\star}|^2dt\right)^{\frac12}\right]\\
	&\leq CE\left[(L_T^{\theta})^2\right]^{\frac12}
	E\left[\int_0^T|\sigma_t^{tr}\pi_t^{\star}|^2dt\right]^{\frac12}<+\infty.
\end{align*}
Consequently, $X^{\pi^{\star}}$ is a martingale under
$\mathbb{Q}^{\theta}$ and $\pi^{\star}\in\mathcal{A}_D^M$. Combining
(\ref{inequality1}) and (\ref{inequality2}), we conclude
\begin{equation*}
-\exp(-\alpha(x-Y_0))=\sup_{\pi\in\mathcal{A}_D^M}E[-\exp(-\alpha(X_T^{\pi}-F))]=E[-\exp(-\alpha(X_T^{\pi^{\star}}-F))].
\end{equation*}

\chapter{Utility indifference valuation}\label{Chapter5}

In this chapter, we apply the results obtained in Chapter
\ref{section:main} to utility indifference pricing of derivatives
with unbounded payoffs. The notion of utility indifference valuation
was proposed in Hodges and Neuberger \cite{HodgesNeuberger} and further developed in Davis
\cite{Davis}.  We refer to the monograph \cite{Carmona} and more references
therein for an overview of utility indifference pricing and related
topics.

To define the utility indifference price for a derivative with
payoff $F$, we also need to consider the optimization problem for
the investor without selling (or buying) the derivative. This
involves the investor investing only in the risk-free bond and risky
assets themselves, and the corresponding optimal trading strategy is
denoted as $\pi^{\star}(0)$. To emphasize the dependence of
$\pi^{\star}$ on $F$, we also write it as $\pi^{\star}(F)$.

\begin{definition}\label{definition1} [{\bf Utility indifference valuation and
		hedging}]
	
	Suppose that Assumption \ref{assumption1} holds. Then, the utility
	indifference price $C_0(F)$ of the derivative with payoff $F$ is
	defined by the solution to
	\begin{equation}\label{indifference_equ}
		\sup_{\pi\in\mathcal{A}_D^{conv}}E\left[-e^{-\alpha\left(X_T^{x+C_0(F)}(\pi)-F
			\right)}\right]= \sup_{\pi\in\mathcal{A}_D^M}E\left[-e^{-\alpha
			X_T^{x}(\pi)}\right],
	\end{equation}
	where $\mathcal{A}^{conv}_D$ and $\mathcal{A}_D^M$ are given in
	Definitions \ref{admiss2} and \ref{admiss3}, respectively.
	
	The hedging strategy for the derivative is defined by the difference
	in the two optimal trading strategies
	$\pi^{\star}(F)-\pi^{\star}(0)$.
\end{definition}

\begin{remark} Since the payoff $F$ is unbounded, we choose
	different admissible sets for the two optimization problems in
	(\ref{indifference_equ}) due to their different natures.
	
	If $F\geq 0$, then $C_0(F)$ is interpreted as the selling price of
	$F$. Since in this case $F$ is automatically bounded from below,
	according to Theorem \ref{thm:opti}, we can enlarge the admissible
	set from $\mathcal{A}_D^{conv}$ to $\mathcal{A}_D^M$. Therefore, the two
	optimization problems in (\ref{indifference_equ}) are solved under
	the same admissible set $\mathcal{A}_D^M$.
	
	On the other hand, if $F\leq 0$, then $-C_0(-F)$ can be interpreted
	as the buying price of $-F$.
\end{remark}

From Theorems \ref{thm:opti2} and \ref{thm:opti}, we have
$$-e^{-\alpha(x+C_0(F)-Y_0(F))}=-e^{-\alpha(x-Y_{0}(0))},$$
where $Y(F)$ is (the first component of) the unique solution to
BSDE$(F,f)$ (cf. (\ref{BSDEexpo})). Herein, we use $Y(F)$ to
emphasize the dependence of the solution component $Y(F)$ on its terminal data
$F$. We thus obtain the utility indifference price of the derivative
with the payoff $F$ as
\begin{equation}\label{primal_representation }
	C_0(F)=Y_0(F)-Y_0(0).
\end{equation}
The associated hedging strategy satisfies
\begin{align}\label{hedging}
	&\sigma_t^{tr}\pi_t^{\star}(F)-\sigma_t^{tr}\pi_t^{\star}(0)\\
	=&\
	\mbox{Proj}_{\sigma_t^{tr}\mathcal{C}^c}\left(Z_t(F)+\frac{\theta_t}{\alpha}\right)
	-\mbox{Proj}_{\sigma_t^{tr}\mathcal{C}^c}\left(Z_t(0)+\frac{\theta_t}{\alpha}\right)
,\ \hbox{$\mathbb{P}$-a.s.,}\ \text{for}\ t \in
		[0,T].\notag
\end{align}

\section[A convex dual representation]{A convex dual representation of utility indifference price}\label{section:dual}

Motivated by Delbaen et al.~\cite{dhr}, we provide a \textit{convex dual representation} of
the solution component $Y(F)$ in this section, which on one hand completes the
proof of Theorem \ref{theorem:BSDE1}, and on the other hand, gives a
convex dual representation of the utility indifference price
$C_0(F)$.

For the specific example considered in Chapter
\ref{section:example}, the convex dual representation will reduce to
the well known minimal entropy representation of the utility
indifference price (see \cite{dgr}, \cite{El} and \cite{MS} with
more references therein).

We next introduce the admissible set of the convex dual problem. For
an $\mathbb{R}^m$-valued predictable process $q\in L^2{[0,T]}$, we
define its stochastic exponential as $$
L^q:=\mathcal{E}(\int_0^{\cdot}q_u^{tr}dB_u).$$ If $L_T^q$ has finite
entropy, i.e. $E[L_T^q\ln L_T^q]<+\infty$, then De la
Vall\'ee-Poussin theorem implies that $L^q$ is in Class (D) and
therefore a uniformly integrable martingale. We can then define a
probability measure $\mathbb{Q}^q$ on $\mathcal{F}_T$ by
$\frac{d\mathbb{Q}^q}{d\mathbb{P}}:=L_T^q$, and introduce the
admissible set
\begin{align*}
	\mathcal{A}_D^{\star}=&\left\{q\in L^2[0,T]:
	L_T^q=\mathcal{E}_T(\int_0^{\cdot}q_u^{tr}dB_u)\ \mbox{has finite
		entropy,}\right.\\
	&\left.\ \text{and}\
	E^{\mathbb{Q}^q}\left[|F|+\int_0^T|f^*(s,q_s)|ds\right]<+\infty,
	\mbox{where}\ \frac{d\mathbb{Q}^q}{d\mathbb{P}}=L_T^q\right\}
	\text{.}
\end{align*}%

\begin{theorem}\label{theorem_dual}
	Suppose that Assumptions \ref{assumption1} holds. Then, the solution component
	$Y$ to BSDE$(F,f)$ for $f$ given by (\ref{driver}) admits the following
	convex dual representation
	\begin{equation}\label{dual_formula}
		Y_t=\esssup_{q\in\mathcal{A}_D^{\star}}E^{\mathbb{Q}^q}
		\left[\left.F-\int_t^{T}{f}^{\star}(s,q_s)ds\right\vert\mathcal{F}_t\right],
	\end{equation}
	where
	${f}^{\star}:[0,T]\times\Omega\times\mathbb{R}^{m}\rightarrow\mathbb{R}\cup\{+\infty\}$
	is the convex dual of $f$ in (\ref{dual}). Moreover, there exists an
	optimal density process $q^{\star}\in\mathcal{A}_D^{\star}$
	such that
	\begin{equation}\label{dual_formula_2}
		Y_t=E^{\mathbb{Q}^{q^{\star}}}
		\left[\left.F-\int_t^{T}{f}^{\star}(s,q_s^{\star})ds\right\vert\mathcal{F}_t\right].
	\end{equation}
\end{theorem}

\begin{proof}

For any $q\in\mathcal{A}_D^{\star}$, we define
\begin{equation*}
	Y_t^{q}:=E^{\mathbb{Q}^q}
	\left[\left.F-\int_t^{T}{f}^{\star}(s,q_s)ds\right\vert\mathcal{F}_t\right],
\end{equation*}
which is finite due to the integrability condition in the admissible
set $\mathcal{A}_D^{\star}$.

Note that $Y_t^{q}-\int_0^{t}f^{\star}(s,q_s)ds$, $t\in[0,T]$, is a
uniformly integrable martingale under $\mathbb{Q}^q$, so similar to
the arguments as in Chapter 5.8 of Karatzas and Shreve \cite{KS}, it follows from the
martingale representation theorem (under the filtration
$\{\mathcal{F}_t\}_{t\geq 0}$) that
\begin{equation}\label{BSDE_qq}
	Y_t^{q}-\int_0^{t}f^{\star}(s,q_s)ds=\left(F-\int_0^{T}{f}^{\star}(s,q_s)ds\right)-\int_t^{T}(Z_s^{q})^{tr}dB_s^{q},
\end{equation}
for some $\mathbb{R}^m$-valued predictable density process $Z^{q}\in
L^2[0,T]$, where $B^{q}_t:=B_t-\int_0^{t}q_sds$, $t\in[0,T]$, is an
$m$-dimensional Brownian motion under $\mathbb{Q}^q$.

On the other hand, we rewrite BSDE (\ref{BSDEexpo}) under
$\mathbb{Q}^q$ as
\begin{equation}\label{BSDE_q}
	Y_{t}=F+\int_{t}^{T}\left(f(s,Z_{s})-Z_s^{tr}q_s\right)ds-\int_{t}^{T}Z_{s}^{tr}dB_{s}^q,\quad
	t\in[0,T].
\end{equation}

For integer $j\ge 1$, we introduce the following stopping time
$$\sigma_j=T\wedge\inf\left\{t\in [0,T]: \max\left\{\int_0^t|Z_s|^2ds, \int_0^t|Z^q_s|^2ds\right\} >j\right\},$$
so that both $\int_0^{\cdot\wedge\sigma_j}Z_s^{tr}dB_s^{q}$ and
$\int_0^{\cdot\wedge\sigma_j}(Z_s^{q})^{tr}dB_s^q$ are martingales
under $\mathbb{Q}^q$.

Combining (\ref{BSDE_qq}) and (\ref{BSDE_q}) and taking the
conditional expectation with respect to $\mathcal{F}_t$ give
\begin{equation*}
	Y_{t}-Y_t^{q}=E^{\mathbb{Q}^q}\left[\left.Y_{\sigma_{j}}-Y_{\sigma_{j}}^q+
	\int_t^{\sigma_j}\left(f(s,Z_{s})-Z_s^{tr}q_s+{f}^{\star}(s,q_s)\right)ds\right\vert\mathcal{F}_t\right].
\end{equation*}
By the Fenchel-Moreau theorem, we then deduce that, for any
$q\in\mathcal{A}_D^{\star}$,
\begin{equation*}
	f(s,Z_{s})-Z_s^{tr}q_s+{f}^{\star}(s,q_s)\geq 0,
\end{equation*}
and thus
\begin{equation}\label{inequality3}
	Y_t-Y_t^{q}\geq
	E^{\mathbb{Q}^q}\left[\left.Y_{\sigma_{j}}-Y_{\sigma_{j}}^q\right\vert\mathcal{F}_t\right].
\end{equation} Note that $Y_{\sigma_{j}}\rightarrow F$,
$Y_{\sigma_{j}}^q\rightarrow F$, $\mathbb{Q}^q$-a.s., and $Y^{q}$ is
uniformly integrable under $\mathbb{Q}^q$, so we only need to prove
that $Y$ is uniformly integrable under $\mathbb{Q}^q$ in order to
pass to the limit in (\ref{inequality3}). For this purpose, we
recall that $e^{\alpha Y^+}\in\mathcal{S}^{p}$ and $e^{\varepsilon
	Y^-}\in\mathcal{S}^{1}$ (cf. Theorem \ref{theorem:BSDE1}). Then, by
applying the inequality (\ref{Young_inequ}) to $L_T^{q}Y_{\star}^+$
(with $p$ replaced by $p\alpha$) and to $L_T^{q}Y_{\star}^-$ (with
$p$ replaced by $\varepsilon$), we get
\begin{align}\label{inequality6}
	E^{\mathbb{Q}^q}[Y_{\star}]\leq&\ E[L^q_TY_{\star}^+]+E[L^q_TY_{\star}^-]\notag\\
	\leq&\ \frac{E[L^q_{T}\ln L_T^q]}{p\alpha}-\frac{E[L_T^q]\ln
		p\alpha}{p\alpha}+E[e^{p\alpha Y_{\star}^+}]\notag\\
	&+\frac{E[L^q_{T}\ln L_T^q]}{\varepsilon}-\frac{E[L_T^q]\ln
		\varepsilon}{\varepsilon}+E[e^{\varepsilon Y_{\star}^-}]<+\infty.
\end{align}
Hence, by letting $j\rightarrow\infty$ in (\ref{inequality3}), we
conclude that $Y_t\geq Y_t^{q}$ for any
$q\in\mathcal{A}_D^{\star}$.

To prove the equality, we set $q^{\star}_s\in\partial f_{z}(s,Z_s)$,
for $s\in[0,T]$. Then, (\ref{FM_theorem_2}) implies
\begin{equation}\label{dual_formula_3}
	f(s,Z_{s})-Z_s^{tr}q_s^{\star}+{f}^{\star}(s,q_s^*)=0,
\end{equation}
from which we obtain $Y_t=Y_t^{q^{\star}},$ for $t\in[0,T]$. It
remains to prove $q^{\star}\in\mathcal{A}_D^{\star}$, which is
established in the following lemma.
\end{proof}

\begin{lemma}\label{lemma000} The optimal density process
	$q^{\star}\in\mathcal{A}_D^{\star}$.
\end{lemma}

\begin{proof}
	We first show $q^{\star}\in L^2[0,T]$. Indeed, following
	(\ref{lowerbound}) and (\ref{dual_formula_3}), we have
	\begin{align*}
		f(s,Z_{s})&=Z_s^{tr}q_s^{\star}-{f}^{\star}(s,q_s^{\star})\\
		&\leq Z_s^{tr}q_s^{\star}-\frac{|q_s^{\star}|^2}{2\alpha}
		\leq\frac{|q_s^{\star}|^2}{4\alpha}+\alpha|Z_s|^2-\frac{|q_s^{\star}|^2}{2\alpha},
	\end{align*}
	and together with the lower bound of the generator $f$ in Lemma \ref{lemma_generator}, i.e.
	$f(s,Z_s)\geq-Z_s^{tr}\theta_s-\frac{|\theta_s|^2}{2\alpha}$, we
	further obtain
	$$\frac{1}{4\alpha}|q_s^{\star}|^2\leq \alpha|Z_s|^2+Z_s^{tr}\theta_s+\frac{|\theta_s|^2}{2\alpha}.$$
	The assertion then follows by noting that $Z\in M^2\subset L^2[0,T]$
	and $\theta$ is bounded.
	
	Next, Lemma \ref{lemma} implies that
	$L^{q^{\star}}_T=\mathcal{E}_T(\int_0^{\cdot}(q_u^{\star})^{tr}dB_u)$
	has finite entropy.
	
	We conclude the proof by verifying that both $f^{\star}$ and $F$ are
	intergrade under $\mathbb{Q}^{q^{\star}}$. Firstly, the lower bound
	of $f^{\star}$ in (\ref{lowerbound}) yields
	$$E^{\mathbb{Q}^{q^{\star}}}\left[\int_0^Tf^{\star}(u,q_u^{\star})du\right]\geq E^{\mathbb{Q}^{q^{\star}}}\left[\int_0^{T}\frac{1}{2\alpha}|q_u^{\star}|^2du\right]\geq 0.$$
	On the other hand, it follows from BSDE (\ref{BSDE_qq}) under
	$\mathbb{Q}^{q^{\star}}$ and the fact $Y=Y^{q^{\star}}$ that
	\begin{align*} E^{\mathbb{Q}^{q^{\star}}}\left[\int_0^{\sigma_j}f^{\star}(u,q_u^{\star})du\right]=E^{\mathbb{Q}^{q^{\star}}}\left[Y_{\sigma_{j}}^{q^{\star}}-Y_0^{q^{\star}}\right]\leq
		E^{\mathbb{Q}^{q^{\star}}}\left[Y_{\star}\right]-Y_0.
	\end{align*}
	According to the inequality (\ref{inequality6}) (with $q$ replaced
	by $q^{\star}$), $E^{\mathbb{Q}^{q^{\star}}}[Y_{\star}]<+\infty$, so
	we have verified the integrability of $f^{\star}$. Similarly,
	applying the inequality (\ref{Young_inequ}) to
	$L_{\sigma_j}^{q^{\star}}F$ yields that $F$ is integrable under
	$\mathbb{Q}^{q^{\star}}$:
	\begin{align*}
		E^{\mathbb{Q}^{q^{\star}}}[|F|]=&\ E[L_T^{q^{\star}}F^+]+E[L_T^{q^{\star}}F^-]\notag\\
		\leq&\ \frac{E[L^{q^{\star}}_{T}\ln
			L_T^{q^{\star}}]}{p\alpha}-\frac{E[L_T^{q^{\star}}]\ln
			p\alpha}{p\alpha}+E[e^{p\alpha F^+}]\notag\\
		&+\frac{E[L^{q^{\star}}_{T}\ln
			L_T^{q^{\star}}]}{\varepsilon}-\frac{E[L_T^{q^{\star}}]\ln
			\varepsilon}{\varepsilon}+E[e^{\varepsilon F^-}]<+\infty,
	\end{align*}
	and the proof is complete.
\end{proof}

As a direct consequence of Theorem \ref{theorem_dual}, we obtain a
convex dual representation of the utility indifference price
\begin{equation}\label{dual_representation}
	C_0(F)=\sup_{q\in\mathcal{A}_D^{\star}}E^{\mathbb{Q}^q}
	\left[F-\int_0^{T}{f}^*(s,q_s)ds\right]-\sup_{q\in\mathcal{A}_D^{\star}}E^{\mathbb{Q}^q}
	\left[-\int_0^{T}{f}^*(s,q_s)ds\right].
\end{equation}
Note that, in general, the above convex dual representation may not
be the minimal entropy representation (see, for example, \cite{dgr},
\cite{El} and \cite{MS}) due to the appearance of non-subspace
portfolio constraints. If it is a subspace portfolio constraint, the
above convex dual representation is precisely the minimal entropy
representation as shown in Chapter \ref{section:example}.

\section{Asymptotics for the risk aversion parameter
}\label{section:asymptotic}

We study the asymptotics of the utility indifference price for the
risk aversion parameter $\alpha$ in this section. We make the
following assumption on the payoff $F$.

\begin{assumption}\label{assumption3}
	The payoff $F$ satisfies $E[e^{p^{\prime}F^+}]<\infty$ for any
	$p^{\prime}\geq 1$, and there exists a constant $k>0$ such that
	$F^{-}\leq k$.
\end{assumption}

\begin{remark} To address the asymptotics as the risk aversion parameter $\alpha\to +\infty$,
	it is obvious that we need to assume  $F^+$ is exponentially
	integrable with any order.
	
	Under Assumption \ref{assumption3}, by Theorem
	\ref{theorem:BSDE1} and Lemma \ref{lemma2}, BSDE$(F,f)$ with $f$ given by (\ref{driver}) admits a unique solution $(Y,Z)$, such that
	$e^{\alpha Y^+}\in\mathcal{S}^{p^{\prime}}$,
	$Y^{-}\in\mathcal{S}^{\infty}$, and $Z\in M^{p^{\prime}}$.
\end{remark}

The existing literature (e.g. \cite{dgr}, \cite{El} and \cite{MS})
addresses only the case of a subspace constraint on portfolios,
which allows them to work on the dual problem, and to use the
corresponding minimal entropy representation in the asymptotic
analysis of the utility indifference price. In our more general case of portfolio constraints, the minimal
entropy representation does not hold anymore, as demonstrated in the
convex dual representation (\ref{dual_representation}). We shall
work on the primal problem by considering the BSDE representation
(\ref{primal_representation }) of the utility indifference price. To
facilitate our subsequent discussion,  we further impose the
following \emph{cone condition} on the constraint set.

\begin{definition}\label{admiss4}
	[{\bf Admissible strategies with constraints $\mathcal{A}_D^{cone}$}]
	
	The set of admissible trading strategies $\mathcal{A}_D^{cone}$
	is the same as $\mathcal{A}_D^{conv}$ in Definition \ref{admiss2}, except
	that the constraint set $\mathcal{C}$ is replaced by
	$\mathcal{C}^c$, where ${\mathcal C}^{c}$ is a closed and convex cone in
	$\mathbb R^{d}$ with $0\in {\mathcal C}^{c}$.
\end{definition}

We use the superscript $^{cone}$ to highlight the dependency on the cone condition of $\mathcal{C}^c$.
To emphasize the dependence on the risk aversion parameter $\alpha$,
we write BSDE$(F,f^{\alpha})$ with its solution
$(Y^{\alpha}(F),Z^{\alpha}(F))$, and the corresponding utility
indifference price and associated hedging strategy as
$C_0^{\alpha}(F)$ and $\pi^{\alpha,\star}(F)-\pi^{\alpha,\star}(0)$.

Recall the generator $f^{\alpha}$ in (\ref{driver}) (with
$\mathcal{C}$ replaced by the cone $\mathcal{C}^{c}$) has the form
\begin{equation*}\label{driver1}
	f^{\alpha}(t,z) = \frac{\alpha}{2}
	\left|\mbox{Proj}_{\sigma_t^{tr}\mathcal{C}^c}(z +
	\frac{1}{\alpha} \theta_t)-(z +
	\frac{1}{\alpha} \theta_t)\right|^2 - z^{tr}\theta_t - \frac{1}{2 \alpha}
	|\theta_t|^2.
\end{equation*}

Since $\mathcal{C}^c$ is a cone, it follows that
$\alpha\mbox{Proj}_{\sigma_t^{tr}\mathcal{C}^c}(z +
\frac{1}{\alpha} \theta_t)=\mbox{Proj}_{\sigma_t^{tr}\mathcal{C}^c}(\alpha z +
\theta_t)$, and therefore,
\begin{equation}\label{invariant1}
	\alpha f^{\alpha}(t,z)=f^{1}(t,\alpha
	z).
\end{equation}
In turn, if $(Y^{\alpha}(F),Z^{\alpha}(F))$ is the unique solution to
BSDE$(F, f^{\alpha})$, then
\begin{align*}
	\alpha Y_t^{\alpha}(F)&=\alpha F+\int_t^{T}\alpha
	f^{\alpha}(u,Z_u^{\alpha}(F))du-\int_t^{T}(\alpha
	Z_u^{\alpha}(F))^{tr}dB_u\\
	&=\alpha F+\int_t^{T}f^{1}(u,\alpha
	Z_u^{\alpha}(F))du-\int_t^{T}(\alpha Z_u^{\alpha}(F))^{tr}dB_u,
\end{align*}
so $(\alpha Y^{\alpha}(F),\alpha Z^{\alpha}(F))$ solves BSDE$(\alpha
F,f^1)$. But according to Theorem \ref{theorem:BSDE1} and Lemma
\ref{lemma2} in Chapter \ref{section:main}, BSDE$(\alpha F,f^1)$
admits a unique solution under Assumption \ref{assumption3}, so we
must have $(\alpha Y^{\alpha}(F),\alpha Z^{\alpha}(F))=(Y^{1}(\alpha
F), Z^{1}(\alpha F))$, in particular,
\begin{equation}\label{invariant2}
	(\alpha Y^{\alpha}(0),\alpha Z^{\alpha}(0))=(Y^{1}(0), Z^{1}(0)).
\end{equation}
The above scaling property is crucial to the asymptotic analysis of
the utility indifference price $C^{\alpha}_0(F)$.

Next, we define
\begin{equation}\label{construction}
	(C_t^{\alpha}(F),H_t^{\alpha}(F)):=(Y_t^{\alpha}(F)-Y_t^{\alpha}(0),Z_t^{\alpha}(F)-Z_t^{\alpha}(0)).
\end{equation}
Then it is immediate to check that
\begin{align}\label{equation2}
	C_t^{\alpha}(F)=&\ F+\int_t^{T}\left(f^{\alpha}(u,H_u^{\alpha}(F)+Z_u^{\alpha}(0))-f^{\alpha}(u,Z_u^{\alpha}(0))\right)du\notag\\
	&-\int_t^T(H^{\alpha}_u(F))^{tr}dB_u\notag\\
	=&\ F+\int_t^{T}\frac{1}{\alpha}\left(f^{1}(u,\alpha H_u^{\alpha}(F)+Z_u^{1}(0))-f^{1}(u,Z_u^{1}(0))\right)du\notag\\
	&-\int_t^T(H^{\alpha}_u(F))^{tr}dB_u,\quad t\in[0,T],
\end{align}
where we used (\ref{invariant1}) and (\ref{invariant2}) in the last
equality.

Furthermore, we introduce the generator
$$g^{\alpha}(t, h, Z_t^{1}(0))=
\frac{\mbox{dist}^2_{\sigma_t^{tr}\mathcal{C}^c}(\alpha
	h+Z_t^1(0)+\theta_t)-\mbox{dist}^2_{\sigma_t^{tr}\mathcal{C}^c}(Z_t^1(0)+\theta_t)}{2\alpha},$$
where $\mbox{dist}_{\sigma_t^{tr}\mathcal{C}^{c}}(\cdot)$ is the
distance function of $\sigma_t^{tr}\mathcal{C}^{c}$. With the
generator $g^{\alpha}$, we rewrite (\ref{equation2}) as
\begin{align}\label{BSDE_indifference_price}
	C_t^{\alpha}(F)=&\
	F+\int_t^Tg^{\alpha}(u,H_u^{\alpha}(F),Z_u^{1}(0))du\notag\\
	&\
	-\int_t^T(H^{\alpha}_u(F))^{tr}\theta_udu-\int_t^T(H^{\alpha}_u(F))^{tr}dB_u,\
	t\in[0,T].\end{align}

Note that if $F$ satisfies Assumption \ref{assumption3}, the above
BSDE (\ref{BSDE_indifference_price}) actually admits a unique
solution $(C^{\alpha}(F),H^{\alpha}(F))$, where $e^{\alpha
	C^{\alpha}(F)^+}\in\mathcal{S}^{p{\prime}}$,
$C^{\alpha}(F)^{-}\in\mathcal{S}^{\infty}$, and $H^{\alpha}(F)\in
M^{p^{\prime}}$ for any $p^{\prime}\geq 1$. Indeed, if
$(C^{\alpha}(F),H^{\alpha}(F))$ is a solution to
(\ref{BSDE_indifference_price}), then with
$(Y^{\alpha}(0),Z^{\alpha}(0))\in\mathcal{S}^{\infty}\times
M^{p^{\prime}}$ as the unique solution to BSDE$(0,f^{\alpha})$, it
is clear that
$(C^{\alpha}(F)+Y^{\alpha}(0),H^{\alpha}(F)+Z^{\alpha}(0))$ solves
BSDE$(F,f^{\alpha})$. But according to Theorem \ref{theorem:BSDE1}
and Lemma \ref{lemma2}, BSDE$(F,f^{\alpha})$ admits a unique
solution. Thus, $(C^{\alpha}(F),H^{\alpha}(F))$ must be the unique
solution to (\ref{BSDE_indifference_price}).

In the following, we will study the asymptotics of the utility
indifference price $C^{\alpha}_0(F)$ via BSDE
(\ref{BSDE_indifference_price}).

\subsection{Asymptotics as $\alpha\rightarrow 0$}

\begin{theorem}\label{thm:asy0} Suppose that Assumption \ref{assumption3} holds, and the admissible set is
	$\mathcal{A}_{D}^{cone}$ as in Definition \ref{admiss4}. Then,
	$\lim_{\alpha\rightarrow 0}C_0^{\alpha}(F)=C^0_0({F})$, where
	$C^0(F)$ is (the first component of) the unique solution to the
	linear BSDE
	\begin{align}\label{linear_BSDE}
		C_t^{0}(F)=&\
		F+\int_t^T(H_u^{0}(F))^{tr}\left(Z_u^{1}(0)-\mbox{Proj}_{\sigma_u^{tr}\mathcal{C}^{c}}(Z_u^{1}(0)+\theta_u)\right)du\notag\\
		&\ -\int_t^{T}(H^{0}_u(F))^{tr}dB_u,\quad t\in[0,T].
	\end{align}
\end{theorem}


\begin{remark}\label{remark} If $\mathcal{C}^c$ is a subspace of
	$\mathbb{R}^d$, we also have the convergence of the optimal trading
	strategy. In fact, by the linearity of the projection operator on
	the subspace $\sigma_t^{tr}\mathcal{C}^c$, we have
	\begin{align*}
		&\sigma_t^{tr}\pi_t^{\alpha,\star}(F)-\sigma_t^{tr}\pi_t^{\alpha,\star}(0)-
		\mbox{Proj}_{\sigma_t^{tr}\mathcal{C}^c}\left(H_t^{0}(F)\right)\\
		=&\mbox{Proj}_{\sigma_t^{tr}\mathcal{C}^c}\left(Z_t^{\alpha}(F)+\frac{\theta_t}{\alpha}\right)
		-\mbox{Proj}_{\sigma_t^{tr}\mathcal{C}^c}\left(Z_t^{\alpha}(0)+\frac{\theta_t}{\alpha}\right)
		-\mbox{Proj}_{\sigma_t^{tr}\mathcal{C}^c}\left(H_t^{0}(F)\right)\\
		=&\mbox{Proj}_{\sigma_t^{tr}\mathcal{C}^c}\left(H_t^{\alpha}(F)-H_t^{0}(F)\right).
	\end{align*}
	Since $H^{\alpha}(F)\rightarrow H^0(F)$ in $L^2[0,T]$ (see the proof
	of Theorem \ref{thm:asy0}), it follows that
	$$\lim_{\alpha\rightarrow 0}\int_0^T|\sigma_t^{tr}\pi_t^{\alpha,\star}(F)-\sigma_t^{tr}\pi_t^{\alpha,\star}(0)-
	\mbox{Proj}_{\sigma_t^{tr}\mathcal{C}^c}\left(H_t^{0}(F)\right)|^2dt=0,
	\ \mbox{$\mathbb{P}$-a.s.}$$
\end{remark}

Based on Theorem \ref{thm:asy0}, we can further represent the
solution to (\ref{linear_BSDE}) in terms of the expected value of
the payoff $F$ under some equivalent probability measure. To see
this, we define $\mathbb{Q}^Z$ as
\begin{equation}\label{density}
	\frac{d\mathbb{Q}^{Z}}{d\mathbb{P}}:=\mathcal{E}_T(\int_0^{\cdot}
	\left(Z^1_t(0)-\mbox{Proj}_{\sigma_t^{tr}\mathcal{C}^c}(Z_t^{1}(0)+\theta_t)\right)^{tr}dB_t).
\end{equation}
Since $\int_0^{\cdot}(Z^1_t(0))^{tr}dB_t$ is a BMO martingale (see
Lemma 12 in \cite{him} or Remark \ref{remark_BMO}), the stochastic exponential in
(\ref{density}) is indeed a uniformly integrable martingale, and
$\mathbb{Q}^{Z}$ is therefore well defined. We then obtain the
asymptotic representation of $C^{\alpha}(F)$ when $\alpha\rightarrow
0$ as
\begin{equation}\label{lowerlimit}
	\lim_{\alpha\rightarrow 0}C_0^{\alpha}(F)= C_0^{0}(F)=
	E^{\mathbb{Q}^{Z}}[F].
\end{equation}

In Chapter \ref{section:example}, we shall show that the probability
measure $\mathbb{Q}^Z$ will reduce to the minimal entropy martingale
measure if the portfolio constraint is a subspace of $\mathbb{R}^d$.


To prove Theorem \ref{thm:asy0}, we start by proving some estimates
of the generator $g^{\alpha}$.

\begin{lemma}\label{lemma1} The generator $g^{\alpha}(t,h,Z_t^{1}(0))$ has the following
	properties:
	
	(i) $g^{\alpha}(t,h,Z_t^{1}(0))$ is nondecreasing in $\alpha$;
	
	(ii) For $\alpha\in(0,1]$, $g^{\alpha}(t,h,Z_t^{1}(0))$ has the
	upper and lower bounds, both of which are independent of $\alpha$,
	\begin{equation}\label{bound_generator}
		h^{tr}\left((Z_t^{1}(0)+\theta_t)-\mbox{Proj}_{\sigma_t^{tr}\mathcal{C}^c}(
		Z_t^{1}(0)+\theta_t)\right)\leq g^{\alpha}(t,h,Z_t^{1}(0)) \leq
		\frac{|h|^2}{2}+h^{tr}m_t,
	\end{equation}
	for some predictable process $m$ satisfying $|m_t|\leq
	|Z_t^1(0)|+|\theta_t|$, for $t\in[0,T]$.
\end{lemma}

\begin{proof}
	
	(i) To prove $g^{\alpha}$ is nondecreasing in $\alpha$, we recall
	that $\sigma_t^{tr}\mathcal{C}^{c}$ is convex, so
	$\mbox{dist}^2_{\sigma_t^{tr}\mathcal{C}^c}(\cdot)$ is convex. It
	then follows that
	$$g^{\alpha}(t, h, Z_t^{1}(0))=
	\frac{\mbox{dist}^2_{\sigma_t^{tr}\mathcal{C}^c}(\alpha
		h+Z_t^1(0)+\theta_t)-\mbox{dist}^2_{\sigma_t^{tr}\mathcal{C}^c}(Z_t^1(0)+\theta_t)}{2\alpha}$$
	is nondecreasing in $\alpha$.\\
	
	(ii) According to (i), we know $g^{\alpha}(t,h,Z_t^{1}(0))\geq
	\lim_{\alpha\rightarrow 0}g^{\alpha}(t,h,Z_t^{1}(0))$. Next, we
	calculate the limit of $g^{\alpha}$ when $\alpha\rightarrow 0$ as
	\begin{align}\label{limit}
		\lim_{\alpha\rightarrow
			0}g^{\alpha}(t,h,Z_t^{1}(0))&=\frac12\left(\frac{\partial}{\partial{\alpha}}\mbox{dist}^2_{\sigma_t^{tr}\mathcal{C}^c}(\alpha
		h+Z_t^1(0)+\theta_t)\right)|_{\alpha=0}\notag\\
		&=h^{tr}\left((\alpha
		h+Z_t^{1}(0)+\theta_t)-\mbox{Proj}_{\sigma_t^{tr}\mathcal{C}^c}(\alpha
		h+Z_t^{1}(0)+\theta_t)\right)|_{\alpha=0}\notag\\
		&=h^{tr}\left((Z_t^{1}(0)+\theta_t)-\mbox{Proj}_{\sigma_t^{tr}\mathcal{C}^c}(
		Z_t^{1}(0)+\theta_t)\right).
	\end{align}
The second equality can be derived as we did in Lemma \ref{lemma_convex} for the calculation of the derivative of
	quadratic distance functions.
	
	On the other hand, using an elementary equality
	$a^2-b^2=(a-b)^2+2b(a-b)$, we rewrite the generator $g^{\alpha}$ as
	\begin{align*}
		&g^{\alpha}(t,h,Z_t^1(0))\\
		=&\frac{1}{2\alpha}|\mbox{dist}_{\sigma_t^{tr}\mathcal{C}^c}(\alpha
		h+Z_t^1(0)+\theta_t)-\mbox{dist}_{\sigma_t^{tr}\mathcal{C}^c}(Z_t^1(0)+\theta_t)|^2\\
		&+\frac{\mbox{dist}_{\sigma_t^{tr}\mathcal{C}^c}(Z_t^1(0)+\theta_t)\left(\mbox{dist}_{\sigma_t^{tr}\mathcal{C}^c}(\alpha
			h+Z_t^1(0)+\theta_t)-\mbox{dist}_{\sigma_t^{tr}\mathcal{C}^c}(Z_t^1(0)+\theta_t)\right)}{\alpha}.
	\end{align*}
	
	By the Lipschitiz continuity of distance functions, $g^{\alpha}$ is
	dominated by
	$$g^{\alpha}(t,h,Z_t^1(0))\leq \frac{\alpha}{2}|h|^2+h^{tr}m_t\leq \frac{1}{2}|h|^2+h^{tr}m_t,$$
	where
	$$m_t:=\frac{\mbox{dist}_{\sigma_t^{tr}\mathcal{C}^c}(Z_t^1(0)+\theta_t)\left(\mbox{dist}_{\sigma_t^{tr}\mathcal{C}^c}(\alpha
		h+Z_t^1(0)+\theta_t)-\mbox{dist}_{\sigma_t^{tr}\mathcal{C}^c}(Z_t^1(0)+\theta_t)\right)}{\alpha|h|^2}h.$$
	Thus, we conclude by noting that $|m_t|\leq |Z_t^1(0)|+|\theta_t|$,
	for $t\in[0,T]$.
\end{proof}

\textbf{Proof of Theorem \ref{thm:asy0}.} We apply the stability
property of quadratic BSDE with bounded terminal data in Lemma \ref{lemma_stability} to study the
limit of $C^{\alpha}(F)$ when $\alpha\rightarrow 0$.

To this end, we consider BSDE$(F,g^{\alpha})$ for $\alpha\in(0,1]$
(cf. (\ref{BSDE_indifference_price})):
\begin{align}\label{truncated_BSDE_11}
	C_t^{\alpha}(F)=&\
	F+\int_t^Tg^{\alpha}(u,H_u^{\alpha}(F),Z_u^{1}(0))du\notag\\
	&\
	-\int_t^T(H^{\alpha}_u(F))^{tr}\theta_udu-\int_t^T(H^{\alpha}_u(F))^{tr}dB_u,\
	t\in[0,T].
\end{align}

Since the generator $g^{\alpha}$ satisfies (\ref{bound_generator})
and the terminal condition $F$ satisfies Assumption
\ref{assumption3}, the comparison theorem for quadratic BSDE with
unbounded terminal data (see section 3 of \cite{BH2}) implies
$$\underline{C}_t^0\leq C_t^{\alpha}(F)\leq C_t^{\alpha^{\prime}}(F)\leq  \overline{C}_t^1,$$
for $0<\alpha\leq \alpha^{\prime}\leq 1$, where $\underline{C}^{0}$
solves BSDE
\begin{align}\label{linear_BSDE_1}
	\underline{C}_t^{0}=&\
	-k+\int_t^T(\underline{H}_u^{0})^{tr}\left(Z_u^{1}(0)-\mbox{Proj}_{\sigma_u^{tr}\mathcal{C}^{c}}(Z_u^{1}(0)+\theta_u)\right)du\notag\\
	&\ -\int_t^{T}(\underline{H}^{0}_u)^{tr}dB_u,\quad t\in [0,T],
\end{align}
and $\overline{C}^1$ solves BSDE
\begin{align}\label{quadratic_BSDE}
	\overline{C}^1_t=&\ F^++\int_t^T\frac12\left|\overline{H}^1_u\right|^2du\notag\\
	&-\int_t^T\left(\overline{H}^1_u\right)^{tr}(dB_u-m_udu),\quad t\in
	[0,T].
\end{align}
It is routine to check that both $\underline{C}^0$ and
$\overline{C}^1$ have the explicit expressions
$$\underline{C}^0_t=-k;\quad \overline{C}^1_t=\ln E^{\mathbb{Q}^m}\left[e^{F^+}|\mathcal{F}_t\right],$$
where $\mathbb{Q}^m$ is defined as
$\frac{d\mathbb{Q}^m}{d\mathbb{P}}:=L^{m}_T=\mathcal{E}_T(\int_0^{\cdot}m_u^{tr}dB_u)$.

We claim that $\overline{C}_t^{1}<+\infty$. Indeed, since
$\int_0^{\cdot}(Z_u^1(0))^{tr}dB_u$ is a BMO martingale (see Lemma
12 in \cite{him} or Remark \ref{remark_BMO}), $\int_0^{\cdot}m_u^{tr}dB_u$ is also a BMO
martingale. It then follows from reverse H\"older's inequality (see \cite{Kaz}) that
there exists some $p>1$ such that
$$E\left[\left.\left(\frac{L_T^m}{L_t^m}\right)^p\right|\mathcal{F}_t\right]\leq C$$
for some constant $C>0$. In turn, H\"older's inequality implies
\begin{align*}
	\ln E^{\mathbb{Q}^m}\left[e^{F^+}|\mathcal{F}_t\right]&=\ln
	E\left[\frac{L_T^m}{L_t^m}e^{F^+}|\mathcal{F}_t\right]\\
	&\leq\frac{1}{p}\ln
	E\left[\left.\left(\frac{L_T^m}{L_t^m}\right)^p\right|\mathcal{F}_t\right]+(1-\frac{1}{p})\ln
	E[e^{\frac{p}{p-1}F^+}|\mathcal{F}_t]<+\infty.
\end{align*}

Next, we pass to the limit in (\ref{truncated_BSDE_11}). But since
the upper bound of $g^{\alpha}$ involves $m$ (cf.
(\ref{bound_generator})), which is typically unbounded due to the
unboundedness of $Z^1(0)$, and the terminal condition $F$ is
unbounded, the stability property does not apply directly. To
overcome this difficulty, we apply the localization argument as in
Theorem \ref{theorem:BSDE1} and define the stopping time
$$\tau_j=T\wedge\inf\left\{t\in [0,T]: \max\left\{\int_0^{t}|Z^1_s(0)|^2ds,\overline{C}_t^{1}\right\} >j\right\},$$
for integer $j\geq 1$.

Then $(C^{\alpha}_j(t),H_j^{\alpha}(t)):=(C^{\alpha}_{t\wedge
	\tau_j}(F),H_t^{\alpha}(F) 1_{\{t\le\tau_j\}})$, $t\in[0,T]$,
satisfies
\begin{align}\label{BSDE_indifference_price_local}
	C_j^{\alpha}(t)=&\
	F_j^{\alpha}+\int_t^T1_{\{u\leq \tau_j\}}g^{\alpha}(u,H_j^{\alpha}(u),Z_u^1(0))du\notag\\
	&\
	-\int_t^T(H^{\alpha}_j(u))^{tr}\theta_udu-\int_t^T(H^{\alpha}_j(u))^{tr}dB_u,
\end{align}
where $F_j^{\alpha}=C_j^{\alpha}(T)=C^{\alpha}_{\tau_j}(F)$.

For fixed $j$, we observe that $C_j^{\alpha}(\cdot)$ is bounded, and
the generator of BSDE (\ref{BSDE_indifference_price_local})
satisfies
$$\left|1_{\{u\leq
	\tau_j\}}g^{\alpha}(u,h,Z_u^1(0))\right|\leq 1_{\{u\leq
	\tau_j\}}|Z_u^1(0)|^2+|\theta_u|^2+|h|^2,$$ for $u\in[0,T]$, so
$\int_0^{T}1_{\{u\leq \tau_j\}}|Z_u^{1}(0)|^2du\leq j$. Moreover,
according to (\ref{limit}),
\begin{align*}
	&1_{\{u\leq
		\tau_j\}}g^{\alpha}(u,h,Z_u^1(0))\\
	\rightarrow&\ 1_{\{u\leq
		\tau_j\}}h^{tr}\left((Z_u^{1}(0)+\theta_u)-\mbox{Proj}_{\sigma_u^{tr}\mathcal{C}^c}(
	Z_u^{1}(0)+\theta_u)\right),
\end{align*}
as $\alpha\rightarrow 0$. Consequently, the quadratic BSDE
(\ref{BSDE_indifference_price_local}) satisfies the stability
property conditions in section 2.2 of \cite{Morlais}.


Hence, setting $C^0_{j}(t)=\inf_{\alpha}C_{j}^{\alpha}(t)$, it
follows from the stability property of quadratic BSDE with bounded
terminal data (see Lemma \ref{lemma_stability}) that there exists
$H_{j}^0(\cdot)\in M^2$ such that $\lim_{\alpha \rightarrow
	0}H_{j}^{\alpha}(\cdot)=H_{j}^0(\cdot)$ in $M^2$, and
$\left(C_{j}^0(\cdot),H_{j}^0(\cdot)\right)$ satisfies
\begin{align}\label{BSDE_indifference_price_local2}
	C_j^{0}(t)=&\
	F_j+\int_t^{\tau_j}(H_j^{0}(u))^{tr}\left(Z_u^{1}(0)-\mbox{Proj}_{\sigma_u^{tr}\mathcal{C}^{c}}(Z_u^{1}(0)+\theta_u)\right)du\notag\\
	&\ -\int_t^{\tau_j}(H_j^{0}(u))^{tr}dB_u,
\end{align}
where $F_{j}=C_{j}^0(T)=\inf_{\alpha}C_{\tau_{j}}^{\alpha}(F)$. The
linear BSDE (\ref{linear_BSDE}) then follows by sending
$j\rightarrow\infty$ in (\ref{BSDE_indifference_price_local2}).

\subsection{Asymptotics as $\alpha\rightarrow\infty$}

We complete the asymptotic analysis of $C^{\alpha}(F)$ by
considering the situation $\alpha\rightarrow\infty$ in the next
theorem.  It is based on
Peng's monotonic limit theorem (see Theorem 2.4 in Peng  \cite{Peng}), which is recalled below.

\begin{lemma}[\cite{Peng}]\label{PengMonotone}
Let $(Y^n,Z^n)\in\mathcal{S}^2\times M^2$ be the solution of the BSDE
\begin{equation}\label{g_BSDE}
Y_t^n=F+\int_t^TdA^n_u-\int_t^Tg(u,Y^n_u,Z^n_u)du-\int_t^T(Z^n_u)^{tr}dB_u
\end{equation}
where $F$ is square integrable, the generator $g$ satisfies the Lipschitz continuous condition with a Lipschitiz constant $C_g$:
$$|g(t,y,z)-g(t,\bar{y},\bar{z})|\leq C_g(|y-\bar{y}|+|z-\bar{z}|)$$
for $(t,y,z)\in[0,T]\times\mathbb{R}\times\mathbb{R}^m$ and $(t,\bar{y},\bar{z})\in[0,T]\times\mathbb{R}\times\mathbb{R}^m$, and $A^{n}\in\mathcal{S}^2$ is an increasing process. If $Y^n$ is increasing and converges to $Y\in\mathcal{S}^2$. then there exist $(Z,A)\in M^2\times\mathcal{S}^2$ such that $A$ is increasing, $A^n\rightarrow A$ weakly in $\mathcal{S}^2$ and $Z^n\rightarrow Z$ weakly in $M^2$, strongly in $M^p$ for $p<2$ with $(Y,Z,A)$ satisfies (\ref{g_BSDE}).
\end{lemma}

\begin{theorem}\label{thm:asy} Suppose that Assumption \ref{assumption3} holds, and the admissible set is
	$\mathcal{A}_D^{cone}$ as in Definition \ref{admiss4}.
	Moreover, suppose that the following constrained BSDE
	\begin{align}\label{constraint_BSDE}
		C_t^{\infty}(F)&=
		F+\int_t^TdA^{\infty}_u(F)-\int_t^T(H_u^{\infty}(F))^{tr}\theta_udu-\int_t^{T}(H^{\infty}_u(F))^{tr}dB_u,\notag\\
		&\ \mbox{subject\ to}\
		H_u^{\infty}(F)\in\sigma_u^{tr}\mathcal{C}^c,\ \mbox{for}\ a.e.\
		u\in[0,T],
	\end{align}
	admits at least one solution
	$(\overline{C}^{\infty}(F),\overline{H}^{\infty}(F),\overline{A}^{\infty}(F))\in\mathcal{S}^2\times
	{M}^2\times\mathcal{S}^2$ with $\overline{A}^{\infty}(F)$ being an
	increasing process.
	
	Then,
	$\lim_{\alpha\rightarrow\infty}C_0^{\alpha}(F)=C_0^{\infty}(F)$,
	where $(C^{\infty}(F),H^{\infty}(F),A^{\infty}(F))$ is the minimal
	solution to the constrained BSDE (\ref{constraint_BSDE}). Herein,
	the minimal solution means if
	$(\overline{C}^{\infty}(F),\overline{H}^{\infty}(F),\overline{A}^{\infty}(F))$
	is a solution to (\ref{constraint_BSDE}), then $
	{C}_t^{\infty}(F)\leq \overline{C}_t^{\infty}(F)$,
	$\mathbb{P}$-a.s., for $t\in[0,T]$.
\end{theorem}


\begin{remark}
	
	The assumption on the constrained BSDE (\ref{constraint_BSDE}) means
	that the payoff $F$ can be superreplicated by using a trading
	strategy constrained in the set $\sigma_u^{tr}\mathcal{C}^c$, for
	a.e. $u\in[0,T]$.
	
	Note that if $F$ is bounded, the above assumption is indeed
	satisfied by taking $\overline{C}^{\infty}(F)\equiv||F||_{\infty}$,
	the essential supremum of $F$, $\overline{H}^{\infty}(F)\equiv0$,
	and
	$\overline{A}_u^{\infty}(F)=1_{\{u=T\}}(||F||_{\infty}-F)+1_{\{u<T\}}0$.
	
	Furthermore, if $F$ is bounded, and $\mathcal{C}^c$ is a subspace of
	$\mathbb{R}^d$ as in Remark \ref{remark}, we also have the
	convergence of the optimal trading strategy. To see this, since
	\begin{equation*}
		\sigma_t^{tr}\pi_t^{\alpha,\star}(F)-\sigma_t^{tr}\pi_t^{\alpha,\star}(0)-
		H_t^{\infty}(F)
		=\mbox{Proj}_{\sigma_t^{tr}\mathcal{C}^c}\left(H_t^{\alpha}(F)-H_t^{\infty}(F)\right),
	\end{equation*}
	and $H^{\alpha}(F)\rightarrow H^{\infty}(F)$ weakly in $M^2$,
	strongly in $M^p$ for $p<2$ (see the proof of Theorem
	\ref{thm:asy}), we obtain
	$$\lim_{\alpha\rightarrow\infty}E\left[\int_0^T|\sigma_t^{tr}\pi_t^{\alpha,\star}(F)-\sigma_t^{tr}\pi_t^{\alpha,\star}(0)-
	H_t^{\infty}(F)|^pdt\right]=0.$$
\end{remark}

Let us recall that, according to Cvitanic et al.~\cite{CKS}, the minimal solution
$C^{\infty}(F)$ actually admits a stochastic control representation.
To see this, we introduce the admissible set
\begin{equation*}
	\mathcal{A}_{\Gamma}^{\star}=\cup_{m\geq 0}\left\{ v\in
	L^2[0,T]:\ |v_t|\leq m\ \mbox{and}\ v_t\ \mbox{is\ valued\ in}\
	\Gamma_t,\ t\in[0,T]\right\} \text{.}
\end{equation*}%

The domain $\Gamma_t$ is define as follows: For $t\in[0,T]$, given
the closed and convex cone $\sigma_t^{tr}\mathcal{C}^c$, we define
its support function
$\delta^{\star}_{\sigma_t^{tr}\mathcal{C}^c}(\cdot)$ as the convex
dual of the characteristic function
$\delta_{\sigma_t^{tr}\mathcal{C}^c}(\cdot)$ of
$\sigma_t^{tr}\mathcal{C}^c$,
$$\delta^{\star}_{\sigma_t^{tr}\mathcal{C}^c}(v)=\sup_{z\in\mathbb{R}^m}\left(z^{tr}v-\delta_{\sigma_t^{tr}\mathcal{C}^c}(z)\right)$$
for $(t,v)\in[0,T]\times\mathbb{R}^{m}$.

Then, $\delta^{\star}_{\sigma_t^{tr}\mathcal{C}^c}(\cdot)$ is valued
in $\mathbb{R}\cup\{+\infty\}$, and is bounded on compact subsets of
the barrier cone
$$\Gamma_t=\left\{v\in\mathbb{R}^{m}:\delta^{\star}_{\sigma_t^{tr}\mathcal{C}^c}(v)<+\infty\right\}.$$

In our case, since $\sigma_t^{tr}\mathcal{C}^c$ is a closed and
convex cone, it follows that $\Gamma_t=\{v\in\mathbb{R}^{m}:
z^{tr}v\leq 0\ \text{for}\ z\in\sigma_t^{tr}\mathcal{C}^{c}\}$, and
$\delta^{\star}_{\sigma_t^{tr}\mathcal{C}^c}\equiv 0$ on $\Gamma_t$.

It then follows from Theorem 5.1 and Corollary 5.1 of \cite{CKS}
(with $\mu=-\theta$) that
\begin{align}\label{upperlimit}
	\lim_{\alpha\rightarrow \infty}C_0^{\alpha}(F)&=C_0^{\infty}(F)\notag\\
	&=
	\sup_{v\in\mathcal{A}_{\Gamma}^{\star}}E^{\mathbb{Q}^{v}}\left[F-\int_0^T\delta^{\star}_{\sigma^{tr}_u\mathcal{C}^{c}}(v_u)du\right]\notag\\
	&=
	\sup_{v\in\mathcal{A}_{\Gamma}^{\star}}E^{\mathbb{Q}^{v}}\left[F\right],
\end{align}
where
\begin{equation}\label{density2}
	\frac{d\mathbb{Q}^{v}}{d\mathbb{P}}:=\mathcal{E}_T(\int_0^{\cdot}
	\left(v_t-\theta_t\right)^{tr}dB_t).
\end{equation}
In Chapter \ref{section:example}, we shall show that
$C_0^{\infty}(F)$ is nothing but the superreplication price of $F$
under MLMM when the portfolio constraint is a subspace.


\textbf{Proof of Theorem \ref{thm:asy}.} We start with
bounded $F$, and proceed to general $F$ by an approximation
procedure.

(i) \textbf{The case that $F$ is bounded.} We first rewrite
(\ref{BSDE_indifference_price}) as
\begin{align}\label{BSDE_indifference_price2}
	\overline{C}_t^{\alpha}(F)=&\
	\left(F-\int_0^T\frac{\mbox{dist}^2_{\sigma_u^{tr}\mathcal{C}^c}
		\left(Z_u^1(0)+\theta_u\right)}{2\alpha}du\right)+\int_t^TdA_u^{\alpha}(F)\notag\\
	&\
	-\int_t^T(H^{\alpha}_u(F))^{tr}\theta_udu-\int_t^T(H^{\alpha}_u(F))^{tr}dB_u,
\end{align}
where
\begin{equation}\label{auxilaryprocess}
	\overline{C}_t^{\alpha}(F):=C_t^{\alpha}(F)-\int_0^t\frac{\mbox{dist}^2_{\sigma_u^{tr}\mathcal{C}^c}
		\left(Z_u^1(0)+\theta_u\right)}{2\alpha}du,
\end{equation}
and $A^{\alpha}$ is an adapted, continuous and increasing process
defined as
\begin{align}\label{processA}
	A_t^{\alpha}(F):&=\int_0^t\frac{1}{2\alpha}\mbox{dist}^2_{\sigma_u^{tr}\mathcal{C}^c}
	\left(\alpha
	H_u^{\alpha}(F)+Z_u^1(0)+\theta_u\right)du\notag\\
	&=\frac{\alpha}{2}\int_0^t\mbox{dist}^2_{\sigma_u^{tr}\mathcal{C}^c}
	\left(H_u^{\alpha}(F)+\frac{Z_u^1(0)+\theta_u}{\alpha}\right)du.
\end{align}
We regard (\ref{BSDE_indifference_price2}) as a penalized equation
for the constrained BSDE (\ref{constraint_BSDE}), and we shall prove
$(\overline{C}^{\alpha}(F),A^{\alpha}(F),H^{\alpha}(F))\rightarrow
(C^{\infty}(F),A^{\infty}(F),H^{\infty}(F))$ as
$\alpha\rightarrow\infty$.

To this end, we note that since $F$ is bounded and $g^{\alpha}$ is
increasing in $\alpha$ (see Lemma \ref{lemma1}), the comparison
theorem for quadratic BSDE with bounded terminal data then implies
$C^{\alpha}(F)$ is increasing in $\alpha$. Consequently,
$\overline{C}^{\alpha}(F)$ is also increasing in $\alpha$.

We claim that for any bounded solution
$(\overline{C}^{\infty}(F),\overline{H}^{\infty}(F),\overline{A}^{\infty}(F))$
to the constrained BSDE (\ref{constraint_BSDE}) (the bounded
solution means $\overline{C}^{\infty}(F)$ is bounded), it holds that
\begin{equation}\label{boundforprice}
	{C}_t^{\alpha}(F)\leq \overline{C}_t^{\infty}(F),\
	\mbox{$\mathbb{P}$-a.s.},\ \mbox{for}\ t\in[0,T].
\end{equation}
We defer the proof of the above inequality to Lemma \ref{lemma00}.
Then (\ref{boundforprice}) will further imply
$\overline{C}_t^{\alpha}(F)\leq \overline{C}_t^{\infty}(F)$ by the
definition of $\overline{C}^{\alpha}(F)$ in (\ref{auxilaryprocess})
and the non-negativity of quadratic distance functions.

Thus, there exists $C^{\infty}(F)\in\mathcal{S}^{\infty}$ such that
$\overline{C}^{\alpha}(F)\rightarrow C^{{\infty}}(F)$ in
$\mathcal{S}^{\infty}$, and $C_t^{{\infty}}(F)\leq
\overline{C}^{\infty}_t(F)$, for $t\in[0,T]$.

It then follows from Lemma 2.5 in \cite{Peng} that
\begin{equation}\label{estimate2}
	\sup_{\alpha\geq 1}E[|A_T^{\alpha}(F)|^2]\leq C;\quad
	\sup_{\alpha\geq 1}E\left[\int_0^T|H^{\alpha}_t(F)|^2dt\right]\leq
	C,
\end{equation}
for some constant $C>0$.

Applying Lemma \ref{PengMonotone}, we
obtain that there exist $(H^{\infty}(F),A^{\infty}(F))\in
M^2\times\mathcal{S}^2$ such that $A^{\infty}(F)$ is increasing,
$A^{\alpha}(F)\rightarrow A^{\infty}(F)$ weakly in $\mathcal{S}^2$,
and $H^{\alpha}(F)\rightarrow H^{\infty}(F)$ weakly in ${M}^2$,
strongly in ${M}^{p}$ for $p<2$, with
$(C^{\infty}(F),H^{\infty}(F),A^{\infty}(F))$ satisfying the
equation in (\ref{constraint_BSDE}).

We next show that $H^{\infty}_u(F)\in\sigma_u^{tr}\mathcal{C}^{c}$,
for a.e. $u\in[0,T]$. Indeed, by the definition of $A^{\alpha}(F)$
in (\ref{processA}) and the first estimate in (\ref{estimate2}), we
have
\begin{align*}
	&E\left[\int_0^T\mbox{dist}^2_{\sigma_u^{tr}\mathcal{C}^c}
	\left(H_u^{\alpha}(F)+\frac{Z_u^1(0)+\theta_u}{\alpha}\right)du\right]\\
	=&\
	\frac{2E[A_T^{\alpha}(F)]}{\alpha}\leq\frac{2E[|A_T^{\alpha}(F)|^2]^{\frac12}}{\alpha}\leq
	\frac{2C^{\frac12}}{\alpha}\rightarrow 0,\ \text{as}\
	\alpha\rightarrow \infty,
\end{align*}
which forces that $H_u^{\infty}(F)\in\sigma_u^{tr}\mathcal{C}^{c}$,
for a.e. $u\in[0,T]$, i.e. the constraint condition in
(\ref{constraint_BSDE}) holds.

(ii) \textbf{The case that $F$ satisfies Assumption
	\ref{assumption3}.} In general, we approximate $F$ from below by
introducing $F^n:=F\wedge n$. It then follows from the comparison
theorem and the stability property for quadratic BSDE with bounded
terminal data that $C^{\alpha}(F^n)$ is increasing in $n$, and
$$\lim_{n\rightarrow\infty}C_0^{\alpha}(F^n)=\sup_{n}C_0^{\alpha}(F^n)=C_0^{\alpha}(F).$$ On the
other hand, we have also proved in (i) that $C^{\alpha}(F^n)$ is
increasing in $\alpha$, and
$$\lim_{\alpha\rightarrow\infty}C_0^{\alpha}(F^n)=\sup_{\alpha}C_0^{\alpha}(F^n)=C_0^{\infty}(F^n).$$
Thus, by interchanging the above two limiting procedures, we obtain
\begin{align*}
	\lim_{\alpha\rightarrow\infty}C_0^{\alpha}(F)&=\sup_{\alpha}\sup_{n}C_0^{\alpha}(F^n)\\
	&=\sup_{n}\sup_{\alpha}C_0^{\alpha}(F^n)\\
	&=\lim_{n\rightarrow\infty}C_0^{\infty}(F^n).
\end{align*}

According to (\ref{upperlimit}), $C_0^{\infty}(F^n)$ admits the
stochastic control representation
$$C_0^{\infty}(F^n)=\sup_{v\in\mathcal{A}_{\Gamma}^{\star}}E^{\mathbb{Q}^{v}}\left[F^n\right].$$
Since $F^n\rightarrow F$ from below, by monotone convergence
theorem, we have
$$\lim_{n\rightarrow\infty}C_0^{\infty}(F^n)=
\sup_{v\in\mathcal{A}_{\Gamma}^{\star}}
E^{\mathbb{Q}^{v}}\left[F\right].$$ Using (\ref{upperlimit}) again,
we know the right hand side of the above equality is nothing but the
stochastic control representation of $C_0^{\infty}(F)$, i.e. (the
first component of) the minimal solution to the constrained BSDE
(\ref{constraint_BSDE}). Hence,
$\lim_{\alpha\rightarrow\infty}C_0^{\alpha}(F)=C_0^{\infty}(F)$.

\begin{lemma}\label{lemma00} Suppose that the payoff $F$ is bounded, and let $({C}^{\alpha}(F), H^{\alpha}(F))$ be the unique solution to
	BSDE (\ref{BSDE_indifference_price}). Then, for any bounded solution
	$$(\overline{C}^{\infty}(F),\overline{H}^{\infty}(F),\overline{A}^{\infty}(F))$$
	to the constrained BSDE (\ref{constraint_BSDE}), we have
	${C}_t^{\alpha}(F)\leq \overline{C}_t^{\infty}(F)$,
	$\mathbb{P}$-a.s., for $t\in[0,T]$.
\end{lemma}

\begin{proof}
	If
	$(\overline{C}^{\infty}(F),\overline{H}^{\infty}(F),\overline{A}^{\infty}(F))$
	is a bounded solution to the constrained BSDE
	(\ref{constraint_BSDE}), it follows from
	(\ref{BSDE_indifference_price}) and (\ref{constraint_BSDE}) that
	\begin{align*}
		C^{\alpha}_t(F)-\overline{C}_t^{\infty}(F)=&\
		\int_t^Tg^{\alpha}(u,H_u^{\alpha}(F),Z_u^{1}(0))du-\int_t^Td\overline{A}_u^{\infty}(F)\notag\\
		&\
		-\int_t^T(H^{\alpha}_u(F)-\overline{H}^{\infty}_u(F))^{tr}(\theta_udu+dB_u)\notag\\
		\leq &\
		\int_t^Tg^{\alpha}(u,\overline{H}_u^{\infty}(F),Z_u^{1}(0))du\notag\\
		&\
		-\int_t^T(H^{\alpha}_u(F)-\overline{H}^{\infty}_u(F))^{tr}\left((\beta_u+\theta_u)du+dB_u\right),
	\end{align*}
	where $\beta$ is defined as
	$$\beta_u=\frac{g^{\alpha}(u,\overline{H}_u^{\infty}(F),Z_u^{1}(0))-g^\alpha(u,H_u^{\alpha}(F),Z_u^{1}(0))}{|H_u^{\alpha}(F)-\overline{H}_u^{\infty}(F)|^2}\left(H_u^{\alpha}(F)-\overline{H}^{\infty}_u(F)\right),$$
	for $u\in[0,T]$. Using the similar arguments as in the proof of the
	upper bound of $g^{\alpha}$ in Lemma \ref{lemma1}, we deduce
	\begin{equation*}
		|\beta_u|\leq
		|Z_u^{1}(0)|+|\theta_u|+\frac{\alpha}{2}(|\overline{H}^{\infty}_u(F)|+|H^{\alpha}_u(F)|).
	\end{equation*}
	
	But
	$\int_0^{\cdot}(Z^1_u(0)+\overline{H}_u^{\infty}(F)+H_u^{\alpha}(F))^{tr}dB_u$
	is a BMO martingale, which follows along the similar arguments as in
	Lemma 12 of \cite{him} by noting that $F$ is bounded, so
	$\int_0^{\cdot}\beta_u^{tr}dB_u$ is also a BMO martingale. In turn,
	defining $\mathbb{Q}^{\beta}$ as
	$$\frac{d\mathbb{Q}^{\beta}}{d\mathbb{P}}:=\mathcal{E}_T(\int_0^{\cdot}(\beta_u+\theta_u)^{tr}dB_u),$$
	we have
	\begin{equation}\label{estimate}
		C_t^{\alpha}(F)-\overline{C}^{\infty}_t(F)\leq
		E^{\mathbb{Q}^{\beta}}\left[\int_t^Tg^{\alpha}(u,\overline{H}_u^{\infty}(F),Z_u^{1}(0))du|\mathcal{F}_t\right].
	\end{equation}
	
	Recall from Lemma \ref{lemma1} that
	$g^{\alpha}(u,\overline{H}_u^{\infty}(F),Z_u^{1}(0))$ is
	increasing in $\alpha$. Setting $\beta:=1/\alpha$, similar to
	(\ref{limit}), we have
	\begin{align*}
		&\lim_{\alpha\rightarrow\infty}g^{\alpha}(u,\overline{H}_u^{\infty}(F),Z_u^{1}(0))\\
		=&\lim_{\alpha\rightarrow\infty}\frac{1}{2\alpha}\mbox{dist}^2_{\sigma_u^{tr}\mathcal{C}^c}\left(\alpha
		\overline{H}^{\infty}_u(F)+Z_u^{1}(0)+\theta_u\right)\\
		=&\lim_{\beta\rightarrow 0}\frac{1}{2\beta}\left(
		\mbox{dist}^2_{\sigma_u^{tr}\mathcal{C}^c}\left(
		\overline{H}^{\infty}_u(F)+\beta(Z_u^{1}(0)+\theta_u)\right)-\mbox{dist}^2_{\sigma_u^{tr}\mathcal{C}^c}\left(
		\overline{H}^{\infty}_u(F)\right)\right)\\
		=&\frac12\left(\frac{\partial}{\partial\beta}\mbox{dist}^2_{\sigma_u^{tr}\mathcal{C}^c}\left(
		\overline{H}^{\infty}_u(F)+\beta(Z_u^{1}(0)+\theta_u)\right)\right)|_{\beta=0},
	\end{align*}
	where we used the constraint condition
	$\overline{H}^{\infty}_u(F)\in\sigma_u^{tr}\mathcal{C}^c$ in the
	last but one equality. We then calculate the derivative of the above
	quadratic distance function
	$\mbox{dist}^2_{\sigma_u^{tr}\mathcal{C}^c}(\cdot)$ with respect to
	$\beta$ as
	\begin{equation*}
		(Z_u^1(0)+\theta_u)^{tr}\left(
		\overline{H}^{\infty}_u(F)+\beta(Z_u^{1}(0)+\theta_u)-
		\mbox{Proj}_{\sigma_u^{tr}\mathcal{C}^c}(
		\overline{H}^{\infty}_u(F)+\beta(Z_u^{1}(0)+\theta_u)) \right),
	\end{equation*}
	which is $0$ at $\beta=0$. Thus,
	$g^{\alpha}(u,\overline{H}_u^{\infty}(F),Z_u^{1}(0))\leq 0$, from
	which we conclude $C_t^{\alpha}(F)\leq C_t^{\infty}(F)$.
\end{proof}

\section{A specific case with subspace portfolio constraints}\label{section:example}

We show that the results herein cover the existing literature of
indifference valuation (\cite{dgr}, \cite{El} and \cite{MS}), if the
trading strategies stay in a subspace of $\mathbb{R}^d$ rather than
a general set.

To facilitate the discussion, we only consider an example with a
single traded asset. The more general case follows along similar
arguments. Consider a market with a single stock whose coefficients
depend on a single stochastic factor driven by a $2$-dimensional
Brownian motion, namely, $m=2$, $d=1$, and
$$
dS_{t}=b\left( V_{t}\right) S_{t}dt+\sigma \left( V_{t}\right)
S_{t}dB_{1,t},
$$
$$
dV_{t}=\eta \left( V_{t}\right) dt+\kappa _{1}dB_{1,t}+\kappa
_{2}dB_{2,t}.
$$
The payoff $F$ has the form $F=F(V_{\cdot})$, which may depend on
the whole path of the stochastic factor process $V$.

We assume that the two positive constants $\kappa _{1},\kappa _{2}$
satisfy
$|\kappa _{1}|^{2}+|\kappa _{2}|^{2}=1$, the functions $%
b(\cdot )$, $\sigma (\cdot )$ and $\eta(\cdot)$ are uniformly
bounded, and $\sigma(\cdot)>0$. Then, the wealth equation becomes
$dX_{t}^{\pi}=\pi _{t}\sigma(V_t)\left( \theta
(V_{t})dt+dB_{1,t}\right) ,$ where $\theta(V_{t}) =\frac{b( V_{t})
}{\sigma ( V_{t}) }$. We also choose $\mathcal{C}=\mathbb{R}$.

\begin{remark}
{In the Markovian setting, we may impose an additional dissipative condition on the function $\eta(v)$. Specifically, we require that $\eta(v)$ satisfies:
\[
(\eta(v) - \eta(\bar{v}))(v - \bar{v}) \leq -C |v - \bar{v}|^2,
\]
for some constant $C > 0$. Consider a payoff function $F = f(V_T)$, where $f$ is a function with linear growth, such as a call option payoff of the form $\pm (V_T - K)^+$. Under these conditions, $F$ has exponential moments. Consequently, the exponential integrability condition on the payoff $F$ as stated in Assumption \ref{assumption1} is satisfied.}
\end{remark}

In this case, for $z=(z_{1},z_{2}),$ the generator $f(t,z)$ in
(\ref{driver}) reduces to
\begin{equation}\label{example_driver}
	f\left(t,(z_{1},z_{2})\right) =-\theta \left( v\right) z_{1}-\frac{1}{2\alpha}|\theta \left( v\right) |^{2}+%
	\frac{\alpha}{2}|z_{2}|^{2},
\end{equation}
and its convex dual $f^{\star}(t,q)$ with $q=(q_1,q_2)$ in
(\ref{dual}) becomes
\begin{equation}\label{example_driver_dual}
	f^{\star}(t,(q_1,q_2))=1_{\{q_1+\theta(v)=0\}}\frac{|q_2|^2+|\theta(v)|^2}{2\alpha}+
	1_{\{q_1+\theta(v)\neq 0\}}\times(+\infty).
\end{equation}

According to Theorem \ref{theorem_dual}, the convex dual
representation of $Y$ in (\ref{dual_formula}) becomes
\begin{equation*}
	Y_t=\esssup_{\substack{q\in\mathcal{A}^{\star}_D\\q_{1,\cdot}=-\theta(V_{\cdot})}}E^{\mathbb{Q}^q}
	\left[\left.F-\int_t^{T}\frac{|q_{2,s}|^2+|\theta(V_s)|^2}{2\alpha}ds\right\vert\mathcal{F}_t\right].
\end{equation*}
Note that the first component $q_1$ of the density process $q$ must
be the negative market price of risk $-\theta(V_s)$, $s\in[0,T]$. Thus,
under $\mathbb{Q}^q$, the stock price process $S$ follows
\begin{equation}\label{stockprice }
	dS_t=\sigma(V_t)S_tdB_{1,t}^{q},
\end{equation} where
$$B_{1,t}^{q}:=B_{1,t}-\int_0^{t}q_{1,s}ds=B_{1,t}+\int_0^{t}\theta(V_s)ds,\
t\in[0,T],$$ is a one-dimensional Brownian motion under
$\mathbb{Q}^{q}$, so $\mathbb{Q}^q$ is an MLMM.

We also note that
\begin{equation}\label{relation}
	E\left[\frac{d\mathbb{Q}^q}{d\mathbb{P}}\ln\frac{d\mathbb{Q}^q}{d\mathbb{P}}\right]=
	E^{\mathbb{Q}^q}\left[\int_0^T\frac{|q_{2,s}|^2+|q_{1,s}|^2}{2}ds\right].
\end{equation}
Hence, we can rewrite the convex dual representation of $Y$ as the
following minimal entropy representation
\begin{equation}\label{relation2}
	Y_0=\sup_{\substack{q\in\mathcal{A}^{\star}_D\\q_{1,\cdot}=-\theta(V_{\cdot})}}\left(E^{\mathbb{Q}^q}
	\left[F\right]-\frac{1}{\alpha}
	E\left[\frac{d\mathbb{Q}^q}{d\mathbb{P}}\ln\frac{d\mathbb{Q}^q}{d\mathbb{P}}\right]\right)
	.
\end{equation}

Consequently,
\begin{align}\label{formula2}
	C_0(F)=&\
	\sup_{\substack{q\in\mathcal{A}_D^{\star}\\q_{1,\cdot}=-\theta(V_{\cdot})}}\left(E^{\mathbb{Q}^q}
	\left[F\right]-\frac{1}{\alpha}
	E\left[\frac{d\mathbb{Q}^q}{d\mathbb{P}}\ln\frac{d\mathbb{Q}^q}{d\mathbb{P}}\right]\right)
	\notag\\
	&\
	-\sup_{\substack{q\in\mathcal{A}_D^{\star}\\q_{1,\cdot}=-\theta(V_{\cdot})}}\left(
	-\frac{1}{\alpha}
	E\left[\frac{d\mathbb{Q}^q}{d\mathbb{P}}\ln\frac{d\mathbb{Q}^q}{d\mathbb{P}}\right]\right),
\end{align}
which is precisely the minimal entropy representation of the utility
indifference price obtained in \cite{dgr}, \cite{El} and \cite{MS}.

Next, we investigate the asymptotic behavior of the utility
indifference price in this example. One could, of course, work on
(\ref{formula2}) to obtain the asymptotic results as in \cite{El}
and \cite{MS}. However, we shall apply the asymptotic results
obtained in Chapter \ref{section:asymptotic} directly, in particular
(\ref{lowerlimit}) and (\ref{upperlimit}), and compare them with the
asymptotic results established in \cite{El} and \cite{MS}.

We follow the notations in Chapter \ref{section:asymptotic}. We
first show that when $\alpha\rightarrow 0$, the probability measure
$\mathbb{Q}^{Z}$ introduced in (\ref{density}) is the minimal
entropy martingale measure.

Denote the optimal density process by
$q^{\star}=(q_1^{\star},q_2^{\star})$, i.e. $q^{\star}$ is the
maximizer of
\begin{equation*}
	Y_0^{\alpha}(0)=\sup_{\substack{q\in\mathcal{A}_D^{\star}\\q_{1,\cdot}=-\theta(V_{\cdot})}}\left(
	-\frac{1}{\alpha}
	E\left[\frac{d\mathbb{Q}^q}{d\mathbb{P}}\ln\frac{d\mathbb{Q}^q}{d\mathbb{P}}\right]\right)=
	E^{\mathbb{Q}^{q^{\star}}}\left[\int_0^T\frac{|q_{2,s}^{\star}|^2+|\theta(V_s)|^2}{-2\alpha}ds\right],
\end{equation*}
where we used (\ref{relation}) in the second equality. Similar to
the proof of Theorem \ref{theorem_dual}, the martingale
representation theorem implies
\begin{equation}\label{BSDE_qqq}
	dY_t^{\alpha}(0)=\frac{|q_{2,t}^{\star}|^2+|\theta(V_t)|^2}{2\alpha}dt
	+\left(Z_{1,t}^{q^{\star}}dB_{1,t}^{q^{\star}}+Z_{2,t}^{q^{\star}}dB_{2,t}^{q^{\star}}\right),
\end{equation}
for some $\mathbb{R}^2$-valued predictable process
$Z^{q^{\star}}=(Z_1^{q^{\star}},Z_2^{q^{\star}})\in L^2[0,T]$, where
\begin{align*}
	B^{q^{\star}}_{1,t}:=&B_{1,t}-\int_0^{t}q_{1,s}^{\star}ds=B_{1,t}+\int_0^{t}\theta(V_s)ds;\\
	B^{q^{\star}}_{2,t}:=&B_{2,t}-\int_0^{t}q_{2,s}^{\star}ds,\quad
	t\in[0,T],
\end{align*}
is a two-dimensional Brownian motion under $\mathbb{Q}^{q^{\star}}$.

On the other hand, according to the primal BSDE$(0,f)$ with the
generator $f$ given in (\ref{example_driver}), we also have
\begin{align}\label{BSDE_qqqq}
	dY_t^{\alpha}(0)=&
	\left(\theta \left( V_{t}\right)Z_{1,t}^{\alpha}(0)-\frac{\alpha}{2}%
	|Z_{2,t}^{\alpha}(0)|^{2}+\frac{|\theta(V_t)|^2}{2\alpha}\right)
	dt\notag\\
	&+\left(Z_{1,t}^{\alpha}(0)dB_{1,t}+Z_{2,t}^{\alpha}(0)dB_{2,t}\right)\notag\\
	=&\frac{2\alpha Z_{2,t}^{\alpha}(0)q_{2,t}^{\star}-|\alpha
		Z_{2,t}^{\alpha}(0)|^2+|\theta(V_t)|^2}{2\alpha}dt\notag\\
	&+\left(Z_{1,t}^{\alpha}(0)dB_{1,t}^{q^{\star}}+Z_{2,t}^{\alpha}(0)dB_{2,t}^{q^{\star}}\right).
\end{align}

Comparing (\ref{BSDE_qqq}) and (\ref{BSDE_qqqq}) gives
$(Z_{1,t}^{q^{\star}},Z_{2,t}^{q^{\star}})=(Z_{1,t}^{\alpha}(0),Z_{2,t}^{\alpha}(0))$,
and moreover,
$$q_{2,t}^{\star}=\alpha Z_{2,t}^{\alpha}(0)=Z_{2,t}^{1}(0),\quad t\in[0,T],$$
by the scaling property (\ref{invariant2}). Thus, the optimal
density process is
$$(q_{1,t}^{\star},q_{2,t}^{\star})=(-\theta(V_t),Z^{1}_{2,t}(0)),\quad t\in[0,T].$$

Note that the above density is nothing but the density of the
probability measure $\mathbb{Q}^{Z}$ in (\ref{density}):
\begin{equation*}
	\frac{d\mathbb{Q}^{Z}}{d\mathbb{P}}=\mathcal{E}_T
	(\int_0^{\cdot}-\theta(V_t)dB_{1,t}+Z_{2,t}^{1}(0)dB_{2,t}),
\end{equation*}
which means $\mathbb{Q}^{Z}$ is the minimal entropy martingale
measure.

To conclude the paper, we show that $C^{\infty}_0(F)$ in
(\ref{upperlimit}) is the superreplication price of $F$ under MLMM
in this example. To this end, observe that the barrier cone
$\Gamma_t$ reduces to $\Gamma_t=\{v\in\mathbb{R}^{2}: v_1\equiv
0\}$, and the support function
$\delta^{\star}_{\sigma_t^{tr}\mathcal{C}^{c}}\equiv 0$ on
$\Gamma_t$. Thus, (\ref{upperlimit}) becomes
$$C^{\infty}_0(F)=\sup_{v\in\mathcal{A}_{\Gamma}^{\star}}E^{\mathbb{Q}^{v}}[F],$$
where
$$\frac{d\mathbb{Q}^v}{d\mathbb{P}}=\mathcal{E}_T(\int_0^{\cdot}-\theta(V_t)dB_{1,t}+v_{2,t}dB_{2,t}),$$
so the stock price $S$ follows (\ref{stockprice }), and
$\mathbb{Q}^{v}$ is an MLMM, i.e. $C_0^{\infty}(F)$ is the
superreplication price of $F$ under MLMM $\mathbb{Q}^v$.

%
%
%


\chapter{Regime switching markets}\label{Chapter6}

In this chapter, we address the utility maximization problem with random endowment in regime switching markets by employing an \textit{exponential quadratic BSDE} framework.
Let $(\Omega,\mathcal{F},\mathbb{F},\mathbb{P})$ be the filtered
probability space introduced in Chapter \ref{Chapter2}. Assume
the probability space also supports a finite-state Markov chain $\Lambda$
with its transition rate matrix $\mathcal{Q}=\{\lambda^{i,k}\}_{i,k\in I}$ satisfying (i) $\sum_{k\in I}\lambda^{i,k}=0$; (ii)
$\lambda^{i,k}\geq 0$ for $i\neq k$, and admits the
representation
$$d \Lambda_t=\sum_{k,k'\in
I}(k-k')\chi_{\{H_{t-}=k'\}}d{N}_t^{k',k},$$
where $N=(N^{k^{\prime},k})_{k^{\prime},k\in I}$ are independent
Poisson processes each with intensity $\lambda^{k^{\prime},k}$ and augmented filtration $\mathbb{H}=\{\mathcal{H}_t\}_{t\geq 0}$. Let $T_0=0$ and $T_1,T_2,\dots$ be the jump times of the Markov
chain $\Lambda$, and $(\Lambda^{j})_{j\geq 1}$ be a sequence of
$\mathcal{H}_{T_{j}}$-measurable random variables representing the
position of $\Lambda$ in the time interval $[T_{j-1},T_{j})$. Hence, $$\Lambda_t=\sum_{j\geq 1}\Lambda^{j-1}\chi_{[T_{j-1},T_{j})}(t).$$
Without
loss of generality, assume that $\Lambda^0=i\in I$. Denote the
smallest filtration generated by $\mathbb{F}$ and $\mathbb{H}$ as
$\mathbb{G}=\{\mathcal{G}_t\}_{t\geq 0}$, i.e.
$\mathcal{G}_t=\mathcal{F}_t\vee\mathcal{H}_t$.
Each state $k\in I$ of the Markov chain $\Lambda$ represents a market
regime.

In this regime switching market, assume that all the coefficients are $\mathbb{G}$-predictable and bounded stochastic
processes. In particular,  the corresponding $\mathbb{R}^{d\times m}$-valued volatility matrix $\sigma$ and $\mathbb{R}^m$-valued market price of risk $\theta$ admit the representation
\begin{align*}
\sigma_t=\sum_{j\geq 1}\sigma_t^{\Lambda^{j-1}}\chi_{(T_{j-1},T_{j}]}(t), \quad  \theta_t=\sum_{j\geq 1}\theta_t^{\Lambda^{j-1}}\chi_{(T_{j-1},T_{j}]}(t),
\end{align*}
where $\sigma^{k}$ and $\theta^{k}$, for $k\in I$, are $\mathbb{F}$-predictable bounded stochastic
processes. Trading strategies will be chosen from the admissible set $\mathcal{A}_D^{\mathbb{G}}$ as defined below:

\begin{definition} \label{admiss5}
	[{\bf Admissible strategies with constraints $\mathcal{A}_D^{\mathbb{G}}$}]

Let $\mathcal{C}$ be a closed and convex set in $\mathbb{R}^d$ satisfying $0\in\mathcal{C}$.
The set of admissible trading strategies $\mathcal{A}_D^{\mathbb{G}}$ consists of $\mathbb{R}^d$-valued $\mathbb{G}$-predictable processes $\pi$ such that $$\pi
_{t}=\pi_0^i\chi_{\{0\}}(t)+\sum_{j\geq
1}\pi_t^{\Lambda^{j-1}}\chi_{(T_{j-1},T_{j}]}(t),\quad t\in[0,T],$$
where, for $k\in I$, $\pi^k\in L^2[0,T]$ is an $\mathbb{R}^d$-valued $\mathbb{F}$-predictable
processes satisfying $\pi^k_t\in\mathcal{C}$, $\mathbb{P}$-a.s., for $t\in[0,T]$. Moreover, the following class (D) condition holds:
$$\{ V(\tau,X_{\tau}^{\pi};\Lambda_{\tau}):\ \tau\ \mbox{is a $\mathbb{G}$-stopping time
		taking values in}\ [0,T]\}
$$ is a uniformly integrable family, where the conditional value process $V(\cdot,\cdot;\Lambda)$ is defined by
\begin{equation}\label{expoopt_222_regime}
	V\left(t,\xi;\Lambda_t\right):= \esssup_{\substack{\{\pi_u,\ u\in[t,
			T]\}\\\text{admissible}}}E \left[\left. -\exp
	\left(-\alpha\left( \xi + \int_t^T \pi_u^{tr}
	\sigma_u (dB_u + \theta_u du) -F
	\right) \right)\right|\mathcal{G}_t\right],
\end{equation}
for any random variable $\xi\in\mathcal{F}_t$ such that $\xi$ satisfying Assumption \ref{assumption1}.
\end{definition}

Herein, we use the superscript $^{\mathbb{G}}$ to emphasize the dependency on $\mathbb{G}$. Note that, for $\pi\in\mathcal{A}_D^{\mathbb{G}}$, the corresponding wealth equation (\ref{wealth}) becomes
\begin{align}\label{wealth_regime_switching}
	X^{\pi}_t&=x + \int_0^t \pi_u^{tr}
	\sigma_u (dB_u + \theta_u du)\notag\\
 &=x+\int_0^{t}\sum_{j\geq 1}(\pi_u^{\Lambda^{j-1}})^{tr}\sigma_u^{\Lambda^{j-1}}(dB_u+\theta_u^{\Lambda^{j-1}}du)\chi_{(T_{j-1},T_{j}]}(u).
\end{align}

Our utility maximization problem is therefore formulated as follows:
\emph{Solving the optimization problem
\begin{equation}\label{expoopt_regime_switching}
V(0,x;i):=\sup_{\pi\in\mathcal{A}_{D}^{\mathbb{G}}} E \left[ -\exp \left(-\alpha\left(x + \int_0^T \pi_u^{tr}
	\sigma_u (dB_u + \theta_u du) -F \right) \right)  \right],
\end{equation}
with the payoff $F$ satisfying Assumption \ref{assumption1}.}

\begin{remark}
In Chapter \ref{section:regime}, we will show that the conditional value process has the form $V(t,\xi;\Lambda_t)=-e^{-\alpha(\xi-Y_t)}$ with $Y$ solving an exponential quadratic BSDE (see (\ref{BSDEfg})), which \emph{in prior} has nothing to do with the optimization problem (\ref{expoopt_regime_switching}), so there is no circular dependency in the definition of the admissible set $\mathcal{A}_D^{\mathbb{G}}$.
\end{remark}

\section{Multidimensional quadratic BSDE with unbounded terminal data}

To solve the optimization problem (\ref{expoopt_regime_switching}), we utilize a BSDE driven by the Brownian motion $B$ and Markov chain $\Lambda$. For an $|I|\times|I|$ matrix $\lambda=(\lambda^{i,k})_{{i,k}\in I}$, denote its $i$th row as $\lambda^{i}$. The operation of multiplying the $i$th row of matrix $\lambda$ with the $i$th row of another matrix $U$ is represented as $$\lambda^{i}\cdot U^{i}=\sum_{k\in I}\lambda^{i,k} U^{i,k}.$$
The BSDE is formulated as follows:
\begin{equation}\label{BSDEfg}
Y_t = F + \int_t^T\left[ f(s,Z_s)ds +  \lambda^{\Lambda_{s-}}\cdot g(U_s^{\Lambda_{s-}})\right]ds - \int_t^TZ_s^{tr}dB_s - \int_t^TU_{s}^{\Lambda_{s-}}\cdot d\tilde{N}_s^{\Lambda{s-}},
\end{equation}
where $\tilde{N}=(\tilde{N}^{i,k})_{i,k\in I}$, with $\tilde{N}_t^{i,k} = N_t^{i,k} - \lambda^{i,k}t$, $t\in [0,T]$, represent the compensated Poisson martingales of $N$ under the filtration $\mathbb{G}=\mathbb{F}\vee\mathbb{H},$
and the generator $g$ is given by
\begin{equation*}
g(u): = \frac{1}{\alpha}\left(e^{\alpha u}-1-\alpha u\right)\in\mathbb{R},\quad u\in\mathbb{R}.
\end{equation*}
Note that the notation $g(U_s^{\Lambda_{s-}})$ for a row vector process $U_s^{\lambda_{s-}}$ means that the function $g$ is evaluated on each component of $U_s^{\lambda_{s-}}$ without any confusion arising.
Finally, the generator $f$ is given in (\ref{driver}), i.e.
\begin{align*}
f(s,z)&=\frac{\alpha}{2}\mbox{dist}_{\sigma_s^{tr}\mathcal{C}}^2\Big(z+\frac{1}{\alpha}\theta_s\Big)
          - z^{tr}\theta_s - \frac{1}{2 \alpha}
		|\theta_s|^2 \\
  &=\sum_{j\geq 1} \left(\frac{\alpha}{2}\mbox{dist}_{(\sigma_s^{\Lambda^{j-1}})^{tr}\mathcal{C}}^2\Big(z+\frac{1}{\alpha}\theta_s^{\Lambda^{j-1}}\Big)
          - z^{tr}\theta_s^{\Lambda^{j-1}} - \frac{1}{2 \alpha}
		|\theta_s^{\Lambda^{j-1}}|^2\right)\chi_{(T_{j-1},T_{j}]}(s)\\
  &=:\sum_{j\geq 1}f^{\Lambda^{j-1}}(s,z)\chi_{(T_{j-1},T_{j}]}(s), \quad (s,z)\in[0,T]\times \mathbb{R}^{m}
  \end{align*}
Since the generator $f$ has quadratic growth (cf. Lemma \ref{lemma_generator}), and the generator $g$ has an exponential form, BSDE (\ref{BSDEfg}) is referred to as \emph{an exponential quadratic BSDE}, as first introduced in \cite{Matoussi2016} and \cite{Matoussi2019}. See also \cite{Lin2023} for its recent development.

By a solution to BSDE (\ref{BSDEfg}), we mean a triplet of $\mathbb{G}$-progressively measurable processes $(Y,Z,U)=(Y_t,Z_t,U_t)_{t\in [0,T]}$, with values in $\mathbb
R\times\mathbb R^m\times\mathbb{R}^{I\times I}$ such that $\mathbb P$-a.s., $t\mapsto Y_t$ is
right continuous, $t\mapsto Z_t$ belongs to $L^2(0,T)$, i.e.
$\int_0^T|Z_t|^2dt<+\infty$, $$t\mapsto f(t,Z_t)+\lambda^{\Lambda_{t-}}\cdot g(U_t^{\Lambda_{t-}})$$ belongs to
$L^1(0,T)$, and $(Y,Z,U)$ satisfies (\ref{BSDEfg}).

The existence of solutions will be established via solving an associated multidimensional quadratic BSDE system
\begin{align}\label{BSDEf_multi}
Y_t^i&=F+\int_t^T\left[f^i(s,Z_s^i)+\sum_{k\in I}\lambda^{i,k}[g(Y^k_s-Y^i_s)+Y^k_s-Y^i_s]\right]ds-\int_t^T(Z_s^i)^{tr}dB_s\\
&=F+\int_t^T\left[f^i(s,Z_s^i)+\sum_{k\in I}\frac{\lambda^{i,k}}{\alpha}(e^{\alpha(Y^k_s-Y^i_s)}-1)\right]ds-\int_t^T(Z_s^i)^{tr}dB_s\notag,
\end{align}
for $t\in [0,T]$ and $i\in I$.
A solution to the BSDE system (\ref{BSDEf_multi}) is denoted by $(Y^i,Z^i)_{i\in I}$. The solutions of (\ref{BSDEfg}) and (\ref{BSDEf_multi}) are connected via
\begin{align*}
Y_t=\sum_{j\geq 1} Y_t^{\Lambda^{j-1}}\chi_{[T_{j-1},T_{j})}(t),\quad
Z_t=\sum_{j\geq 1} Z_t^{\Lambda^{j-1}}\chi_{(T_{j-1},T_{j}]}(t),
\end{align*}
and $U_t^{i,k}=Y_t^{k}-Y_t^{i}.$
The above connection can be established by applying
It\^o's formula for the Markov chain $\Lambda$. We recall it in
the following lemma, which will be frequently used in the rest of
the paper. Its proof is a straightforward extension of Chapter 9.1.2 in Bremaud \cite{Bremaud} and is thus omitted here.

\begin{lemma}\label{lemma_Ito}
For $i\in I$, let $F^i_t$, $t\geq 0$, be a family of
$\mathbb{F}$-progressively measurable and continuous stochastic
processes. Then,
\begin{align*}
&\sum_{j\geq
1}\left[F^{\Lambda_{T_j}}_{T_j}-F^{\Lambda_{T_{j}-}}_{T_{j}-}\right]\chi_{\{T_{j}\leq
t\}}\\
=&\int_{0}^{t}\sum_{k\in
I}\lambda^{\Lambda_{s-}k}[F^k_s-F^{\Lambda_{s-}}_s]ds+\int_{0}^{t}\sum_{k,k^{\prime}\in
I}[F_s^k-F_s^{k'}]\chi_{\{\Lambda_{s-}=k^{\prime}\}}d\tilde{N}_s^{k^{\prime}k},
\end{align*}
\end{lemma}

The rest of this section is devoted to solving the multidimensional quadratic BSDE system (\ref{BSDEf_multi}), which in turn provides a solution to (\ref{BSDEfg}). The uniqueness of the solution to (\ref{BSDEfg}) will be established in Chapter \ref{section6.3} as a corollary of the convex dual representation of
the solution component $Y$.

The proof of existence relies on the multidimensional comparison theorem
for systems of BSDE established in \cite{HuPeng}. The proof of the following version can be found in Lemma 2.2 in \cite{HLT}.

\begin{lemma}\label{comparsion_lemma} For $i\in I$, consider a system of
BSDE$(F^i,f^i+g^i)$ with the terminal data $F^i$ and the generator
$(f^i+g^i)$, namely,
$$Y^{i}_t=F^i+\int_t^T\left[f^i(s,Z_s^i)+g^i(Y_s^i,Y_s^{-i})\right]ds-\int_t^T(Z_s^i)^{tr}dB_s,\quad t\in [0,T], $$
where
$Y_s^{-i}:=(Y_s^{1},\dots,Y_s^{i-1},Y_s^{i+1},\dots,Y_s^{|I|})$. Let
$(\bar{Y}^i,\bar{Z}^i)$ be a solution of the system of
BSDE $(\bar{F}^i,\bar{f}^i,\bar{g}^i)$ with the
terminal data $\bar{F}^i$ and the generator $(\bar{f}^i+\bar{g}^i)$.
Suppose that

(i) both $F^i$ and $\bar{F}^i$ are square integrable and
satisfying $F^i\leq \bar{F}^i$ for $i\in I$;

(ii) there exist constants $C_f$ and $C_g$ such that, for $i\in
I$ and $z,\bar{z}\in\mathbb{R}^m$, $y=(y^i,y^{-i}),
\bar{y}=(\bar{y}^i,\bar{y}^{-i})\in\mathbb{R}^{|I|}$,
\begin{align}
|f^i(s,z)-f^i(s,\bar{z})|&\leq C_f|z-\bar{z}|\label{Lip_F},\\
|g^i(y^i,y^{-i})-g^i(\bar{y}^i,\bar{y}^{-i})|&\leq
C_g|y-\bar{y}|;\label{Lip_G}
\end{align}

(iii) the generator $g^{i}(y^i,{y}^{-i})$ is increasing in all
of its components other than $y^i$, i.e. it is increasing in
$y^{k}$, for $k\neq i$;

(iv) the following inequalities hold,
\begin{align}
f^{i}(s,\bar{Z}_s^i)&\leq \bar{f}^{i}(s,\bar{Z}_s^i),\label{compare_F}\\
g^i(s,\bar{Y}_s^i,\bar{Y}_s^{-i})&\leq
\bar{g}^i(s,\bar{Y}_s^i,\bar{Y}_s^{-i})\label{compare_G}.
\end{align}
Then, $Y_t^{i}\leq \bar{Y}_t^i$ for $t\in[0,T]$ and $i\in I$.
\end{lemma}

\begin{theorem}\label{theorem:BSDE1_regime}
	Suppose that Assumption \ref{assumption1} holds. Then the BSDE system (\ref{BSDEf_multi}) admits a solution $(Y^i,Z^i)_{i\in I}$, where $e^{\alpha
		(Y^i)^+}\in\mathcal{S}^{p}$, $e^{\varepsilon (Y^i)^-}\in\mathcal{S}^{1}$,
	and $Z^i\in M^2$, i.e.
	\begin{equation*}
		E\left[e^{p\alpha (Y^i)^{+}_{\star}}+e^{\varepsilon
			(Y^i)^{-}_{\star}}+\int_0^T|Z_s^i|^2ds\right]<+\infty,
	\end{equation*}
	where $Y^i_{\star}=\sup_{t\in[0,T]}Y^i_t$ is the running maximum of a
	stochastic process $Y$.
	\end{theorem}

\begin{proof}
We first truncate the terminal data $F$ as in the proof of Theorem \ref{theorem:BSDE1} to define $F^{n,k}:=F^+\wedge n-F^-\wedge k$. Since there is no comparison theorem for general quadratic BSDE systems, we next apply the inf convolution to the generator $f^i$ in the BSDE system (\ref{BSDEf_multi}) as in Remark \ref{remark_approx}, i.e.,
$$f^{i,n}(t,z):=\inf_{q}\{f^i(t,q)+n|q-z|\},$$
which is Lipschitz continuous in $z$ and has the growth property
$$-|\theta_t||z|-\frac{1}{2\alpha}|\theta_t|^2\leq f^{i,n}(t,z)\leq f^i(t,z)\leq \frac{\alpha}{2}|z|^2$$
with $|\theta_t|:=\max_{i\in I}|\theta_t^i|$. 

Consider the following BSDE system
\begin{align}\label{BSDE_system_truncate}
Y_t^{i,(n,k)} &= F^{n,k} - \int_t^T (Z_s^{i,(n,k)})^{\text{tr}} dB_s \\
&\quad + \int_t^T \bigg[ f^{i,n}(s,Z_s^{i,(n,k)}) 
+ \sum_{j\in I} \frac{\lambda^{i,j}}{\alpha} \Big( e^{\alpha\big(p(Y^{j,(n,k)}_s) - p(Y^{i,(n,k)}_s)\big)}-1\Big) \bigg] ds, \notag
\end{align}
where the truncation function \( p(\cdot) \) will be chosen to ensure that \( p(y) \) is bounded and Lipschitz continuous. As a result, the term \( e^{\alpha (p(y^j) - p(y^i))} \) will also be bounded and Lipschitz continuous.
 Note that the above BSDE system has bounded terminal data $F^{n,k}$ and its generator $f^{i,n}(t,z)$ is Lipschitz continuous in $z$. Moreover, the generator
\begin{align*}
\sum_{j\in I}\frac{\lambda^{i,j}}{\alpha}\left(e^{\alpha (p(y^{j})-p(y^{i}))}-1\right)
\end{align*}
is also Lipschitz continuous in $y=(y^{i},y^{-i})$ due to the introduction of the truncation function $p(y)$, and satisfies the monotonic condition (iii) in Lemma \ref{comparsion_lemma}. Therefore, the comparison theorem in Lemma \ref{comparsion_lemma} applies to (\ref{BSDE_system_truncate}), from which we deduce that $Y^{i,(n,k)}$ is increasing in $n$ and decreasing in $k$.

Let $(\underline{Y},\underline{Z})\in\mathcal{S}^{\infty}\times M^2$ solve
BSDE 
\begin{equation}\label{MBSDEl*}
	\underline{Y}_t=F^{n,k}+\int_t^T
	(-|\theta_s||\underline{Z}_s|-\frac{1}{2\alpha}|\theta_s|^2
	)ds-\int_t^T \underline{Z}_s^{tr}dB_s,
\end{equation}
and $(\overline{Y},\overline{Z})\in\mathcal{S}^{\infty}\times M^2$ solve BSDE
\begin{equation}\label{MBSDE2*}
	\overline{Y}_t=F^{n,k}+\int_t^T
	\frac{\alpha}{2}|\overline{Z}_s|^2ds-\int_t^T
	\overline{Z}_s^{tr}dB_s.
\end{equation}
Comparing the generators in (\ref{BSDE_system_truncate}) and (\ref{MBSDEl*}) at $(\underline{Y},\underline{Z})$, we have
$$f^{i,n}(s,\underline{Z}_s)+\sum_{j\in I}\frac{\lambda^{i,j}}{\alpha}\left(e^{\alpha(p(\underline{Y}_s)-p(\underline{Y}_s))}-1\right)\geq -|\theta_s||\underline{Z}_s|-\frac{1}{2\alpha}|\theta_s|^2.$$
Likewise, comparing the generators in (\ref{BSDE_system_truncate}) and (\ref{MBSDE2*}) at $(\overline{Y},\overline{Z})$, we have
$$f^{i,n}(s,\overline{Z}_s)+\sum_{j\in I}\frac{\lambda^{i,j}}{\alpha}\left(e^{\alpha(p(\overline{Y}_s)-p(\overline{Y}_s))}-1\right)\leq \frac{\alpha}{2}|\overline{Z}_s|^2.$$
By the comparison theorem in Lemma \ref{comparsion_lemma}, we deduce that
$\underline{Y}_t\leq Y_t^{i,(n,k)}\leq \overline{Y}_t$, and we may choose the truncation function $$p(y)=\max\{-M(n,k),\min\{y,M(n,k)\}\},$$  where the constant $M(n,k)$, depending on $n$ and $k$, is such that 
$\max\{|\underline{Y}^{-}_t|,|\overline{Y}_t^+|\}\leq M$ for $t\in[0,T]$. Therefore, the truncation function $p(y)$ does not play a role in the BSDE system (\ref{BSDE_system_truncate}), and $(Y^{i,(n,k)},Z^{i,(n,k)})_{i\in I}$ is also a solution to (\ref{BSDE_system_truncate}) without applying the truncation function $p$.    

Furthermore, we obtain from (\ref{MBSDEl*}) and (\ref{MBSDE2*}) that $\underline{Y}$ further bounded from below by the solution component $P$ of the following BSDE
\begin{equation*}
	{P}_t=-F^{-}+\int_t^T
	(-|\theta_s||{Q}_s|-\frac{1}{2\alpha}|\theta_s|^2
	)ds-\int_t^T {Q}_s^{tr}dB_s,
\end{equation*}
and $\overline{Y}_t$ bounded from above by 
\begin{align*}
	\overline{Y}_t&=\frac{1}{\alpha}\ln E[e^{\alpha F^{n,k}}|\mathcal{F}_t]\leq \frac{1}{\alpha}\ln E[e^{\alpha F^+}|\mathcal{F}_t],
\end{align*}

The remainder of the proof proceeds using the analogous localization argument to each equation in the BSDE system, as presented in the proof of Theorem \ref{theorem:BSDE1}. Specifically, we employ the localization sequence \( (\tau_j)_{j \geq 1} \), defined by
\[
\tau_j = T \wedge \inf\left\{t \in [0, T] : \max\left\{\frac{1}{\alpha}\ln E[e^{\alpha F^+}|\mathcal{F}_t], -P_t\right\} > j \right\}.
\]
Since the steps are analogous to the previous proof for Theorem \ref{theorem:BSDE1}, we omit the detailed presentation here.
\end{proof}

\section{Exponential utility maximization in regime switching markets}\label{section:regime}

We are now ready to provide one of the main results in this section, which is the characterization of the conditional value function and the corresponding optimal trading strategy for the optimization problem (\ref{expoopt_regime_switching})
under Assumption \ref{assumption1}.

\begin{theorem} \label{thm:optimal_regime}
	Suppose that Assumption 1 holds. Let $(Y,Z,U)$ be the unique solution
	to BSDE (\ref{BSDEfg}). Then, the value function of
	the optimization problem (\ref{expoopt_regime_switching}) with admissible set
	$\mathcal{A}_D^{\mathbb{G}}$ is given by
	\begin{equation}\label{valuefunction_regime}
		V(0, x;i)= -\exp(-\alpha(x-Y_0)),
	\end{equation}
	and there exists an optimal trading strategy
 $\pi^{\star}\in \mathcal{A}_D^{\mathbb{G}}$ such that
	\begin{equation}  \label{pi_regime}
\sigma_t^{tr}\pi^{\star}_t=\mbox{Proj}_{\sigma_t^{tr}\mathcal{C}}\left(Z_t+\frac{\theta_t}{\alpha}\right),\
		\mbox{$\mathbb{P}$-a.s.},\ \text{for}\ t \in [0,T].
	\end{equation}
\end{theorem}

The proof of Theorem \ref{thm:optimal_regime} follows similar arguments to those of Theorem \ref{thm:opti2}, and thus, we only highlight their major differences.

\begin{lemma}\label{lemma:regime} Let $X^{\pi}$ solve the wealth equation (\ref{wealth_regime_switching}) and $(Y,Z,U)$ solve BSDE (\ref{BSDEfg}). Then,
$$-e^{-\alpha(X_{t}^{\pi}-Y_{t})}=-e^{-\alpha(x-Y_0)}A_t^{\pi}L_t^{\pi},$$
	where
	$$A_t^{\pi}=\exp\left(\alpha\int_0^t
	\left(\frac{\alpha}{2}|\sigma_u^{tr}\pi_u-(Z_u +
	\frac{\theta_u}{\alpha})|^2 - Z_u^{tr}\theta_u - \frac{|\theta_u|^2}{2
		\alpha}-f(u,Z_u)\right)du
	\right),$$ and
	\begin{equation*}
		L_{t}^{\pi}=\mathcal{E}_t\left(\alpha\int_0^{\cdot}(Z_u^{tr}-\pi_u^{tr}\sigma_u)dB_u+\int_0^{\cdot}
(e^{\alpha U_{u}^{\Lambda_{u-}}}-\mathbf{1}^{tr})\cdot d\tilde{N}^{\Lambda_{u-}}_u\right),
	\end{equation*}
with $\mathbf{1}$ representing the column vector of ones.
\end{lemma}

\begin{proof} We denote the continuous part and pure jump part of $-e^{-\alpha(X_{t}^{\pi}-Y_{t})}$ by $V^{c}(t,X_t^{\pi})$ and $V^{d}(t,X_t^{\pi})$, respectively. Observe that on $[T_{j-1},T_j)$, $j\geq 1$, $Y_t=Y_t^{\Lambda^{j-1}}$ and $X_t^{\pi}=(X_t^{\pi})^{\Lambda^{j-1}}$, and they follow
\begin{align*}
dY_t^i&=-\left\{f^i(t,Z_t^i)+{\lambda}^{i}\cdot [g(U_t^{i})+U_t^{i}]\right\}dt+(Z_t^i)^{tr}dB_t,\\
d(X_t^{\pi})^i&=(\pi_t^{i})^{tr}\sigma_t^i(dB_t+\theta_t^{i}dt),
\end{align*}
on the event $\{\Lambda^{j-1}=i\}$. Hence, It\^o's formula yields that
\begin{align*}
&d(-e^{-\alpha((X_{t}^{\pi})^i-Y_{t}^i)})\\
=&\ -\alpha e^{-\alpha((X_{t}^{\pi})^i-Y_{t}^i)}\left\{\frac{\alpha}{2}\left|(\sigma_t^i)^{tr}\pi^i_t-(Z_t^i +
	\frac{\theta^i_t}{\alpha})\right|^2 - (Z^i_t)^{tr}\theta^i_t - \frac{|\theta^i_t|^2}{2
		\alpha}-f^i(t,Z_t^i)\right\}dt\\
&\ -\alpha e^{-\alpha((X_{t}^{\pi})^i-Y_{t}^i)}\left((Z_t^i)^{tr}-(\pi_t^i)^{tr}\sigma_t^i\right)dB_t\\
&\ + e^{-\alpha((X_{t}^{\pi})^i-Y_{t}^i)}\lambda^{i}\cdot(e^{\alpha U_t^{i}}-\mathbf{1}^{tr})dt.
\end{align*}
By summing up $-e^{-\alpha((X_{t}^{\pi})^i-Y_{t}^i)}$ over the interval $[T_{j-1},T_j)$ for $j\geq 1$, we obtain
\begin{align*}
&V^{c}(t,X_t^{\pi})\ =\ \ -e^{\alpha(x-Y_0)}\\
&-\int_0^t \!\!\! \alpha\   e^{-\alpha(X_{u-}^{\pi}-Y_{u-})} \Bigl[\frac{\alpha}{2}\Bigm|(\sigma_u)^{tr}\pi_u-(Z_u +
	\frac{\theta_u}{\alpha})\Bigm|^2\!\! - (Z_u)^{tr}\theta_u - \frac{|\theta_u|^2}{2
		\alpha}-f(u,Z_u)\Bigr]du\\
&\quad\quad+\int_0^t \lambda^{\Lambda_{u-}} e^{-\alpha(X_{u-}^{\pi}-Y_{u-})} (e^{\alpha U_u^{\Lambda_{u-}}}-\mathbf{1}^{tr}) du\\
&\quad\quad -\int_0^t \alpha \, e^{-\alpha(X_{u-}^{\pi}-Y_{u-})}\left(Z_u^{tr}-\pi_u^{tr}\sigma_u\right)dB_u
\end{align*}
On the other hand, the pure jump part $V^{d}(t,X_t^{\pi})$ can be expressed as
\begin{align*}
V^{d}(t,X_t^{\pi})&=\sum_{j\geq 1}\left(-e^{-\alpha(X_{T_{j}}^{\pi}-Y_{T_{j}})}+e^{-\alpha(X_{T_{j}-}^{\pi}-Y_{T_{j}-})}\right)\chi_{T_j\leq t}\\
&=-\int_0^te^{-\alpha(X_{u-}^{\pi}-Y_{u-})}(e^{\alpha U_{u}^{\Lambda_{u-}}}-1)\cdot dN_u^{\Lambda_{u-}}.
\end{align*}
The result then follows by summing up $V^{c}(t,X_t^{\pi})$ and $V^{d}(t,X_t^{\pi})$.
\end{proof}

\textbf{Proof of Theorem \ref{thm:optimal_regime}.} 
We apply the martingale optimality principle. The supermartingale property of $-e^{-\alpha(X_{t}^{\pi}-Y_{t})}$, $t\in[0,T]$, follows from Lemma \ref{lemma:regime} and the definition of the admissible set $\mathcal{A}_{D}^{\mathbb{G}}$. To prove the martingale property of $-e^{-\alpha(X_{t}^{\pi^{\star}}-Y_{t})}$, $t\in[0,T]$,
it is sufficient to demonstrate that, with $\pi^{\star}$ given in (\ref{pi_regime}), the optimal density process $L^{\pi^{\star}}$ is in Class (D). This will further imply that $L^{\pi^{\star}}$ is a uniformly integrable martingale. Consequently, $-e^{-\alpha(X_{t}^{\pi^{\star}}-Y_{t})}$, $t\in[0,T]$, also belongs to Class (D).

We observe that
the convex dual of the generator $g$ is given by $g^{\star}$ with
$$g^{\star}(v)=\sup_{u\in\mathbb{R}}(u v-g(u))=\frac{1}{\alpha}((v+{1})\ln(v+{1})-v)\geq 0,\quad v>-1.$$
On other hand, $\alpha(Z_t-\sigma_t^{tr}\pi_t^{\star})\in\partial_z f(t,Z_t)$ by (\ref{subdifferential}) in Lemma \ref{lemma_convex}. Therefore, in order to show that $L^{\pi^{\star}}_T$ has a finite entropy, we consider
\begin{equation*}		L_{t}^{q^{\star},v^{\star}}=\mathcal{E}_t\left(\int_0^{\cdot}(q_u^{\star})^{tr}dB_u
+v_u^{\star,\Lambda_{u-}}\cdot d\tilde{N}_u^{\Lambda_{u-}}\right),\quad t\in[0,T].
	\end{equation*}
for any $q^{\star}_t\in\partial_z f(t,Z_t)$ and $v_t^{\star}=(v_t^{\star,i,k})_{i,k\in I}$ with $v^{\star,i, k}_t=\partial_ug(U_t^{i,k})>-{1}$, and show that $L_{T}^{q^{\star},v^{\star}}$ has a finite entropy in the following lemma, which will then conclude the proof.

\begin{lemma}\label{lemma_switching}
	The optimal density process $L_T^{q^{\star},v^{\star}}$ has a finite entropy. Hence, by De la Vall\'ee-Poussin theorem, $L^{q^{\star},v^{\star}}$ is in
	Class (D) and, therefore, it is a uniformly integrable martingale. 
 \end{lemma}

\begin{proof}

To show that
	$L^{q^{\star},v^{\star}}_T$
	has finite entropy, we introduce
	the following stopping times: for any integer $j\ge 1$,
	$$\sigma_j:=T\wedge \inf\left\{t\in [0,T]: \max\left\{\int_0^t|q^{\star}_s|^2ds,|v_t^{\star}|,\int_0^t|Z_s|^2ds \right\} >j\right\},$$
	so that $L^{q^{\star},v^{\star}}_{\cdot\wedge\sigma_j}$ is a uniformly
	integrable martingale under $\mathbb{P}$.
	
	We can then define a probability measure $\mathbb{Q}^{q^{\star},v^{\star}}$ on
	$\mathcal{F}_{\sigma_j}$ by
	$\frac{d\mathbb{Q}^{q^{\star},v^{\star}}}{d\mathbb{P}}:=L_{\sigma_j}^{q^{\star},v^{\star}}$,
	together with the Brownian motion
	$B^{q^{\star}}:=B-\int_0^{\cdot}q_u^{\star}du$ and the compensated Poisson martingales $\tilde{N}^{v^{\star}}=(\tilde{N}^{v^{\star},i,k})_{i,k\in I}$ with  $\tilde{N}^{v^{\star},i,k}_t:=\tilde{N}^{i,k}_t-\int_0^{t}v_u^{\star,i,k} \lambda^{i,k}du$, for
	$t\in[0,\sigma_j]$.
	
Applying the Fenchel inequality (\ref{Young_inequ})
	to $L_{\sigma_{j}}^{q^{\star},v^{\star}}(\alpha Y_{\sigma_j})$ gives
	\begin{equation}\label{inequality6}
		E[L_{\sigma_{j}}^{q^{\star},v^{\star}}\alpha Y_{\sigma_j}]\leq
		\frac{E[L_{\sigma_j}^{q^{\star},v^{\star}}\ln L_{\sigma_j}^{q^{\star},v^{\star}}]}{p} -
		\frac{E[L_{\sigma_j}^{q^{\star},v^{\star}}]\ln p}{p}+E[e^{p\alpha
			Y^+_{\star}}]
	\end{equation}
	with $p>1$ given in (\ref{condition}).

We use BSDE
	(\ref{BSDEfg}) to estimate the LHS of the inequality (\ref{inequality6}).
 Under $\mathbb{Q}^{q^{\star},v^{\star}}$, we have
	{\small \begin{align}\label{inequality555}
		E[L_{\sigma_{j}}^{q^{\star}, v^{\star}}\alpha
		Y_{\sigma_j}] & =E^{\mathbb{Q}^{q^{\star},v^{\star}}}[\alpha
		Y_{\sigma_{j}}]\notag\\
	&=E^{\mathbb{Q}^{q^{\star},v^{\star}}}\left[\alpha
		Y_0+\alpha\int_0^{\sigma_{j}}
(-f(u,Z_u)+Z_u^{tr}q_u^{\star})du+\alpha\int_0^{\sigma_j}(Z_u)^{tr}dB_u^{q^{\star}}\right.\notag\\
     &\quad +\alpha\int_0^{\sigma_{j}}\sum_{k\in I} \lambda^{\Lambda_{u-},k}(-g(U_u^{\Lambda_{u-},k})+U_u^{\Lambda_{u-},k}v_u^{\star,\Lambda_{u-},k})du\notag\\
     &\quad \left.+\ \alpha\int_0^{\sigma_{j}}U_u^{\Lambda_{u-}}\cdot d\tilde{N}_u^{v^{\star},\Lambda_{u-}}\right]\notag\\
        &=\alpha Y_0+E^{\mathbb{Q}^{q^{\star},v^{\star}}}\left[\alpha\int_0^{\sigma_j}(f^{\star}(u,q_u^{\star})+\lambda^{\Lambda_{u-}}\cdot g^{\star}(v_u^{\star,\Lambda_{u-}}))du\right]\notag\\
		&\geq \alpha Y_0+E^{\mathbb{Q}^{q^{\star},v^{\star}}}\left[
		\int_0^{\sigma_{j}}\left(\frac12|q_u^{\star}|^2+\alpha\lambda^{\Lambda_{u-}}\cdot g^{\star}(v_u^{\star,\Lambda_{u-}})\right)du\right],
	\end{align}}
	where we used the definitions of convex dual $f^{\star}$ and $g^{\star}$ in the last equality, and the lower bound of $f^{\star}$ (cf.
	(\ref{lowerbound})) in the last inequality.
	
Next, we calculate $E[L_{\sigma_j}^{q^{\star},v^{\star}}\ln
L_{\sigma_j}^{q^{\star},v^{\star}}]$. Denote the continuous part and pure jump part of $L^{q^{\star},v^{\star}}$ by $L^{c}$ and $L^{d}$, respectively. Then, since
$$dL_{t}^{q^{\star},v^{\star}}=L_{t-}^{q^{\star},v^{\star}}
\left[(q^{\star}_t)^{tr}dB_t-\lambda^{\Lambda_{t-}}\cdot v_t^{\star,\Lambda_{t-}}dt+v_t^{\star,\Lambda_{t-}}\cdot dN_t^{\Lambda_{t-}}\right],$$
we have
\begin{align*}
\ln L_t^c&=\int_0^t(q_u^{\star})^{tr}dB_u-\int_0^t\left(\lambda^{\Lambda_{u-}}\cdot v_u^{\star,\Lambda_{u-}}+\frac{1}{2}|q^{\star}_u|^2\right)du\\
&=\int_0^t(q_u^{\star})^{tr}dB_u^{q^{\star}}-\int_0^t\left(\lambda^{\Lambda_{u-}}\cdot v_u^{\star,\Lambda_{u-}}-\frac{1}{2}|q^{\star}_u|^2\right)du
\end{align*}
and
\begin{align*}
\ln L_t^d&=\sum_{0<s\leq t}\ln(L_{s}^{q^{\star},v^{\star}})-\ln(L_{s-}^{q^{\star},v^{\star}})\\
&=\sum_{0<s\leq t}\ln(L_{s-}^{q^{\star},v^{\star}}+L_{s-}^{q^{\star},v^{\star}}\sum_{k\in I}v_s^{\star,\Lambda_{s-},k}\Delta N_s^{\Lambda_{s-},k})-\ln(L_{s-}^{q^{\star},v^{\star}})\\
&=\sum_{0<s\leq t}\sum_{k\in I}\ln(v_s^{\star,\Lambda_{s-},k}+1)\Delta N_s^{\Lambda_{s-},k}\\
&=\int_0^t\sum_{k\in I}\ln(v_s^{\star,\Lambda_{s-},k}+1)dN_s^{\Lambda_{s-},k}\\
&=\int_0^t\sum_{k\in I}\ln(v_s^{\star,\Lambda_{s-},k}+1)\left(d\tilde{N}_s^{v^{\star},\Lambda_{s-},k}+
\lambda^{\Lambda_{s-},k}(v_s^{\star,\Lambda_{s-},k}+1)ds\right).
\end{align*}
Hence, by recalling the expression of the convex dual $g^{\star}$, we obtain
$$E\left [L_{\sigma_j}^{q^{\star},v^{\star}}\ln
L_{\sigma_j}^{q^{\star},v^{\star}}\right]=E^{\mathbb{Q}^{q^{\star},v^{\star}}}\left[\int_0^{\sigma_{j}}\left(\frac12|q_u^{\star}|^2+\alpha\lambda^{\Lambda_{u-}}\cdot g^{\star}(v_u^{\star,\Lambda_{u-}})\right)du\right].$$

Combining the above expression together with (\ref{inequality4}) and (\ref{inequality5}),
	we obtain
	$$(1-\frac{1}{p})E\left [L_{\sigma_j}^{q^{\star},v^{\star}
    }\ln
	L_{\sigma_j}^{q^{\star},v^{\star}}\right]\leq -\alpha Y_0-\frac{\ln
		p}{p}+E\left [e^{p\alpha Y^+_{\star}}\right]<+\infty.$$ The assertion then
	follows by sending $j\rightarrow\infty$ in the above inequality, and
	using Fatou's lemma.
\end{proof}

\section[A convex dual representation]{A convex dual representation of the regime switching model}\label{section6.3}

We provide a convex dual representation of
the solution component $Y$ of BSDE (\ref{BSDEfg}) in this section, which completes the
proof of Theorem \ref{theorem:BSDE1} for its solution uniqueness.

We first introduce the admissible set of the convex dual problem. For
an $\mathbb{R}^m$-valued $\mathbb{G}$-predictable process $q$ and an $\mathbb{R}^{|I|\times |I|}$-valued
$\mathbb{G}$-predicable process $v$,
we
define its stochastic exponential as $$
L^{q,v}:=\mathcal{E}\left(\int_0^{\cdot}q_u^{tr}dB_u+v_u^{\Lambda_{u-}}\cdot d\tilde{N}_u^{\Lambda_{u-}}\right).$$
If $L_T^{q,v}$ has a finite entropy, i.e. $E[L_T^{q,v}\ln L_T^{q,v}]<+\infty$, then De la
Vall\'ee-Poussin theorem implies that $L^{q,v}$ is in Class (D) and is
therefore a uniformly integrable martingale. We can then define a
probability measure $\mathbb{Q}^{q,v}$ on $\mathcal{G}_T$ by
$d\mathbb{Q}^{q,v}:=L_T^{q,v} d\mathbb{P}$, and introduce the
admissible set
\begin{align*}
	\mathcal{A}_D^{\mathbb{G},\star}&=\biggl\{(q,v)\in L^2(0,T)\times L^{g^{\star}}[0,T]:\
	L_T^{q,v}\ \mbox{has a finite
		entropy such}\\
	&\ \text{that}\ E^{\mathbb{Q}^{q,v}}\left[|F|+\int_0^T|f^{\star}(s,q_s)|+|g^{\star}(v_s)|ds\right]<+\infty\  \mbox{with}\ d\mathbb{Q}^{q,v}:=L_T^{q,v}d\mathbb{P}\ \biggr\}
	\text{,}
\end{align*}%
where $L^{g^{\star}}[0,T]$ is the space of $\mathbb{G}$-predictable processes $v=(v^{i,k})_{i,k\in I}$ such that $$\int_0^Tg^{\star}(v^{i,k}_t)dt=\int_0^T\frac{1}{\alpha}[(v^{i,k}_t+1)\ln(v^{i,k}_t+1)-v_t^{i,k}]dt<\infty, \quad \mathbb{P}-a.s..$$

\begin{theorem}\label{theorem_dual_regime}
	Suppose that Assumption~\ref{assumption1} holds. Then, the solution component
	$Y$ to BSDE (\ref{BSDEfg}) admits the following
	convex dual representation
	\begin{equation}\label{dual_formula_regime}
		Y_t=\esssup_{(q,v)\in\mathcal{A}_D^{\mathbb{G},\star}}E^{\mathbb{Q}^{q,v}}
		\left[F-\int_t^{T}[{f}^{\star}(s,q_s)+\lambda^{\lambda_{s-}}\cdot g^{\star}(v_s^{\Lambda_{s-}})]ds\Bigm\vert\mathcal{G}_t\right].
	\end{equation}
Moreover, there exists an
	optimal pair of density processes $(q^{\star},v^{\star})\in\mathcal{A}_D^{\mathbb{G},\star}$
	such that
	\begin{equation}\label{dual_formula_2_regime}
		Y_t=E^{\mathbb{Q}^{q^{\star},v^{\star}}}
		\left[F-\int_t^{T}[{f}^{\star}(s,q_s^{\star})+\lambda^{\Lambda_{s-}}\cdot g^{\star}(v_s^{\Lambda_{s-}})]ds\Bigm \vert\mathcal{G}_t\right].
	\end{equation}
\end{theorem}

\begin{proof}

The proof is similar to that of Theorem \ref{theorem_dual}, but with the addition of new jump components, and we will only highlight the main differences.
For any $(q,v)\in\mathcal{A}_D^{\mathbb{G},\star}$, we define
\begin{equation*}
	Y_t^{q,v}:=E^{\mathbb{Q}^{q,v}}
	\left[F-\int_t^{T}[{f}^{\star}(s,q_s)+\lambda^{\Lambda_{s-}}\cdot g^{\star}(v_s^{\Lambda_{s-}})]ds\Bigm \vert\mathcal{G}_t\right],
\end{equation*}
which is finite due to the integrability condition in the admissible
set $\mathcal{A}_D^{\mathbb{G},\star}$. Then, 
according to Chapter 3 of \cite{Bremaud_0},
the
martingale representation theorem (under the filtration
$\{\mathcal{G}_t\}_{t\geq 0}$) implies that
\begin{align}\label{BSDE_qq_regime}
	Y_t^{q,v}=\ &F-\int_t^{T}\left[{f}^{\star}(s,q_s)+\lambda^{\Lambda_{s-}}\cdot g^{\star}(v_s^{\Lambda_{s-}})\right]ds\notag\\
	&-\int_t^{T}(Z_s^{q,v})^{tr}dB_s^{q}-\int_t^TU_s^{q,v,\Lambda_{s-}}\cdot d\tilde{N}_s^{v,\Lambda_{s-}}
\end{align}
for some $\mathbb{G}$-predictable density processes $(Z^{q,v},U^{q,v})\in\mathbb{R}^{m}\times\mathbb{R}^{|I|\times|I|}$, where $$B^{q}_t:=B_t-\int_0^{t}q_sds, \quad t\in[0,T]$$
 is an
$m$-dimensional Brownian motion, and 
$\tilde{N}^{v}=(\tilde{N}^{v,i,k})_{i,k\in I}$ with  
$$\tilde{N}^{v,i,k}_t:=\tilde{N}^{i,k}_t-\int_0^{t}v_u^{i,k} \lambda^{i,k}du, \quad t\in[0,T]$$
are the compensated Poisson martingales
under $\mathbb{Q}^{q,v}$.

On the other hand, BSDE (\ref{BSDEfg}) also reads
\begin{align}\label{BSDE_q_fg_regime}
Y_{t}=&\ F+\int_{t}^{T}\left(f(s,Z_{s})-Z_s^{tr}q_s\right)ds-\int_{t}^{T}Z_{s}^{tr}dB_{s}^q\notag \\
&\ +\int_t^{T}\sum_{k\in I} \lambda^{\Lambda_{s-},k}(g(U_s^{\Lambda_{s-},k})-U_s^{\Lambda_{s-},k}v_s^{\Lambda_{s-},k})ds\\
&-\int_t^{T}U_s^{\Lambda_{s-}}\cdot d\tilde{N}_s^{v,\Lambda_{s-}}.\notag
\end{align}

For integer $j\ge 1$, we introduce the following stopping time
$$\sigma_j=T\wedge \inf\left\{t\in [0,T]: \max\left\{\int_0^t|Z_s|^2ds, \int_0^t|Z_s^{q,v}|^2ds,|U_t|,|U_t^{q,v}|\right\} >j\right\}, $$
and therefore,  all the processes 
$$\int_0^{\cdot\wedge\sigma_j}Z_s^{tr}dB_s^{q}, \quad 
\int_0^{\cdot\wedge\sigma_j}(Z_s^{q,v})^{tr}dB_s^q, \quad  \int_0^{\cdot\wedge\sigma_j}U_s^{\Lambda_{s-}}\cdot d\tilde{N}_s^{v,\Lambda_{s-}}$$
 and 
 $$
 \int_0^{\cdot\wedge\sigma_j}U_s^{q,v,\Lambda_{s-}}\cdot d\tilde{N}_s^{v,\Lambda_{s-}}
 $$ 
 are $\mathbb{Q}^{q,v}$-martingales. 

Combining (\ref{BSDE_qq_regime}) and (\ref{BSDE_q_fg_regime}) and taking the
expectation conditioned on $\mathcal{G}_t$ give
\begin{align*}
	Y_{t}-Y_t^{q,v}=&\ E^{\mathbb{Q}^{q,v}}\left[Y_{\sigma_{j}}-Y_{\sigma_{j}}^{q,v}+ \int_t^{\sigma_j}\left(f(s,Z_{s})-Z_s^{tr}q_s+{f}^{\star}(s,q_s)\right)ds\right.\\
&\left.\left.+ \int_t^{\sigma_j}\sum_{k\in I}\lambda^{\Lambda_{s-},k}[ g(U_s^{\Lambda_{s-},k})-U_s^{\Lambda_{s-},k}v_s^{\Lambda_{s-},k}+g^{\star}(v_s^{\Lambda_{s-},k})]
ds\right|\mathcal{G}_t\right].
\end{align*}
By the Fenchel-Moreau theorem, we then deduce that, for any
$(q,v)\in\mathcal{A}_D^{\mathbb{G},\star}$,
\begin{align*}
	f(s,Z_{s})-Z_s^{tr}q_s+{f}^{\star}(s,q_s)&\geq 0,\\
g(U_s^{\Lambda_{s-},k})-U_s^{\Lambda_{s-},k}v_s^{\Lambda_{s-},k}+g^{\star}(v_s^{\Lambda_{s-},k})&\geq 0
\end{align*}
and the equalities hold
for $q^{\star}_s\in\partial_z f(s,Z_s)$ and $v_s^{\star}=(v_s^{\star,i,k})_{i,k\in I}$ with $v^{\star,i, k}_t=\partial_ug(U_s^{i,k})$. 
Hence,
$$Y_t-Y_t^{q,v}\geq E^{\mathbb{Q}^{q,v}}\left [Y_{\sigma_j}-Y_{\sigma_j}^{q,v}\bigm |\mathcal{G}_t\right ], \quad  \forall (q,v)\in\mathcal{A}_D^{\mathbb{G},\star}
$$
 and $$Y_t-Y_t^{q^{\star},v^{\star}}= E^{\mathbb{Q}^{q^{\star},v^{\star}}}\left[Y_{\sigma_j}-Y_{\sigma_j}^{q^{\star},v^{\star}}\bigm |\mathcal{G}_t\right].$$
Note both $Y$ and $Y^{q,v}$ are uniformly integrable under $\mathbb{Q}^{q,v}$.
Sending $j\rightarrow \infty$, we obtain that $Y_t\geq Y_t^{q,v}$. We conclude the proof by asserting that $Y_t=Y_t^{q^{\star},v^{\star}}$.
To this end, it
remains to prove $(q^{\star},v^{\star})\in\mathcal{A}_D^{\mathbb{G},\star}$, which is
established in the following lemma.
\end{proof}

\begin{lemma}\label{lemma000_regime} The optimal density process
	$(q^{\star},v^{\star})$ belongs to $\mathcal{A}_D^{\mathbb{G},\star}$.
\end{lemma}

\begin{proof}
	Note that $q^{\star} \in L^2(0,T)$ has been proved in Lemma \ref{lemma000}. To show that $v^{\star} \in L^{g^{\star}}[0,T]$, it is sufficient to observe that for $v^{\star,i,k}_t = \partial_ug(U_s^{i,k}) = e^{\alpha U_t^{i,k}} - 1$,
\begin{align*}
\int_0^Tg^{\star}(v^{\star,i,k}_t)dt &= \int_0^T\frac{1}{\alpha} \left[(v^{\star,i,k}_t+1)\ln(v^{\star,i,k}_t+1)-v_t^{\star,i,k}\right]dt \\
&= \int_0^T\frac{1}{\alpha} \left[e^{\alpha U^{i,k}_t}(\alpha U^{i,k}_t-1)+1 \right]dt<\infty.
\end{align*}

Next, Lemma \ref{lemma_switching} implies that
$L^{q^{\star},v^{\star}}_T =
\mathcal{E}_T\left(\int_0^{\cdot}(q_u^{\star})^{tr}dB_u+v_u^{\star,\Lambda{u-}}\cdot d\tilde{N}_u^{\Lambda{u-}}\right)$
has a finite entropy.

Finally, the integrability of $F$ under $\mathbb{Q}^{q^{\star},v^{\star}}$ has been proved in Lemma \ref{lemma000}. We conclude the proof by verifying that $f^{\star} + g^{\star}$ is integrable under $\mathbb{Q}^{q^{\star},v^{\star}}$. Since both $f^{\star}$ and $g^{\star}$ are nonnegative, it suffices to prove that their integrals are bounded from above, which follows from BSDE (\ref{BSDE_qq_regime}):
\begin{align*}
E^{\mathbb{Q}^{q^{\star},v^{\star}}}\left[\int_0^{\sigma_j}\left(f^{\star}(u,q_u^{\star})+\lambda_{\Lambda_{s-}}\cdot g^{\star}(v^{\Lambda_{s-}})\right)du\right] =\ & E^{\mathbb{Q}^{q^{\star},v^{\star}}}\left[Y_{\sigma_{j}}^{q^{\star}}-Y_0^{q^{\star}}\right]\\
 \leq\ & E^{\mathbb{Q}^{q^{\star},v^{\star}}}\left[Y_{\star}\right]-Y_0 < \infty.
\end{align*}

\end{proof}

%
%
%


\chapter{The presence of consumption}\label{Chapter7}
In this chapter, our objective is to enhance the utility maximization problem by incorporating a consumption process into our model. We consider the filtered probability space $(\Omega,\mathcal{F},\mathbb{F},\mathbb{P})$ introduced in Chapter \ref{Chapter2}. The investor aims to consume a stream of goods with a rate of $\beta\geq 0$, which modifies the wealth equation as follows:
\begin{equation}\label{wealth_cons}
	X^{\pi,C}_t=x + \int_0^t \pi_u^{tr}
	\sigma_u (dB_u + \theta_u du)-\int_0^t\beta C_udu, \quad x\in\mathbb{R},
\end{equation}
where $C \in L^1(0,T)$ is an $\mathbb{F}$-predictable consumption process to be determined, i.e. $\int_0^T|C_t|dt<\infty$, $\mathbb{P}$-a.s. The investor’s intertemporal utility of consumption is given by
$$U(c)=-e^{-\alpha c},\quad c\in\mathbb{R}.$$
Note that negative consumption means the infusion of funds (or income), which also generates utility.

The investor~\textit{chooses an admissible pair of optimal trading strategy $\pi^{\star}$ and optimal consumption process $C^{\star}$ so as to maximize the expected utility of the terminal wealth and accumulated consumption:}
{\small
\begin{equation}\label{expoopt_cons} V(0,x):=\sup_{\substack{\{\pi_u,C_u,\ u\in[0,T]\}\\{\text{admissible}}}} E \left[ -\exp \left(-\alpha\left( X_T^{\pi,C} -F \right)\right)-\int_0^T\lambda\exp\left(-\alpha C_u\right)du  \right],
\end{equation}}
\textit{where $\lambda \geq 0$ denotes the weight of utility of consumption}. 

It is important to note that when $\beta = \lambda = 0$, we revert to the utility maximization problem (\ref{expoopt}) discussed in Chapter \ref{section:main}. For our current analysis, we assume that $(\pi,C) \in \mathcal{A}_D^{C}$, where we use the superscript $^C$ to highlight the dependency on the consumption $C$.

\begin{definition}\label{admiss6}
	[{\bf Admissible strategies with constraints $\mathcal{A}_D^C$}]

Let $\mathcal{C}$ be a closed and convex set in $\mathbb R^{d}$ and $0\in\mathcal{C}$. 	
The admissible set $\mathcal{A}_D^C$ consists of all trading strategies $\pi$ and consumption processes $C$ such that $\pi\in L^2[0,T]$ is predictable, satisfying $\pi_t\in\mathcal{C}$, $\mathbb{P}$-a.s., {for} $t \in
		[0,T]$; and $C\in L^1[0,T]$ is predictable. Moreover, the following integrability condition holds:  
  \begin{equation}\label{consumption_integrability_condtion}
 E[\int_0^T\exp\left(-\alpha C_u\right)du]<\infty, 
  \end{equation}
  and the class (D) condition holds:
\begin{equation}\label{class_D_condition}
\{ V(\tau,X_{\tau}^{\pi,C}):\ \tau\ \mbox{is an $\mathbb{F}$-stopping time
		taking values in}\ [0,T]\}
\end{equation} is a uniformly integrable family, where the conditional value process $V(\cdot,\cdot)$ is defined as 
{\small
\begin{align}\label{expoopt_222_con}
	V\left(t,\xi\right):=& \esssup_{\substack{\{\pi_u,C_u\ u\in[t,
			T]\}\\\text{admissible}}}E \Biggl[ -\int_t^T\lambda\exp(-\alpha C_u)du\\
&
-\exp
	\left(-\alpha\left( \xi + \int_t^T \left[\pi_u^{tr}
	\sigma_u (dB_u + \theta_u du)-\beta C_u du\right] -F
	\right) \right)
\Bigm |\mathcal{F}_t\Biggr],\notag
\end{align}
}
for any random variable $\xi\in\mathcal{F}_t$ satisfying Assumption \ref{assumption1}. 
\end{definition}

\begin{remark}
In Theorem \ref{thm:optimal_cons}, we will demonstrate that the conditional valueprocess takes the form of
$V(t,\xi)=-e^{-\alpha(h(t)\xi-Y_t)}$, with $Y$ solving BSDE$(F,f)$ with $f$ given by (\ref{driver_cons}).
Hence, the integrability condition is equivalent to stating that $e^{-\alpha h(t)X_{t}^{\pi,C}}e^{\alpha Y_{t}}$, $t\in[0,T]$, belongs to Class (D). It is important to note that this condition is independent of the optimization problem and its derivation. Therefore, there is no circular dependency in the definition of the admissible set.
\end{remark}

\section[Quadratic BSDE ]{Quadratic BSDE with unbounded terminal data and generator linear in $y$}

 To solve the optimization problem (\ref{expoopt_cons}), we consider the quadratic BSDE$(F,f)$ with the generator $f$ given by
\begin{align}\label{driver_cons}
		f(t,y,z) =& \frac{\alpha}{2}\min_{\pi\in {\mathcal C}}
		\left|h(t)\sigma_t^{tr}\pi-(z +
		\frac{1}{\alpha} \theta_t)\right|^2 - z^{tr}\theta_t - \frac{1}{2 \alpha}
		|\theta_t|^2\\
&\quad -\beta h(t)y-\frac{\beta h(t)}{\alpha}(\ln\frac{\beta h(t)}{\lambda}-1),\notag
	\end{align}
where $h$ solves the (backward) differential Riccati equation:
$$h(t)=1-\int_t^T\beta h^2(s)ds, \quad t\in [0,T]$$
which has a unique solution
$$h(t)=\frac{1}{\beta(1+T-t)}
\in\left[\frac{1}{\beta(1+T)},\frac{1}{\beta}\right], \quad t\in [0,T].$$
Note that when $\beta=\lambda=0$, the above BSDE will reduce to BSDE (\ref{BSDEexpo}) studied in Chapter \ref{Chapter3}. On the other hand, when $F$ is bounded, the aforementioned BSDE, along with the associated optimal consumption and investment problem, has been studied in \cite{CH}. Its extension to a regime switching setting has been investigated in \cite{HSX2022} recently.

\begin{theorem}\label{theorem:BSDE_cons}
	Suppose that Assumption \ref{assumption1} holds.
Then BSDE$(F,f)$ with the generator $f$ given in (\ref{driver_cons}) admits a unique solution $(Y,Z)$, where $e^{\alpha
		Y^+}\in\mathcal{S}^{p}$, $e^{\varepsilon Y^-}\in\mathcal{S}^{1}$,
	and $Z\in M^2$, i.e.
	\begin{equation*}
		E\left[e^{p\alpha Y^{+}_{\star}}+e^{\varepsilon
			Y^{-}_{\star}}+\int_0^T|Z_s|^2ds\right]<+\infty,
	\end{equation*}
	where $Y_{\star}=\sup_{t\in[0,T]}Y_t$ is the running maximum of a
	stochastic process $Y$.
\end{theorem}

\begin{proof} We only prove the solution existence and leave its solution uniqueness in Chapter \ref{Section_convex_consump}. Define the positive constant
$$h^{\star}=\sup_{t\in[0,T]}\left|\frac{\beta h(t)}{\alpha}(\ln\frac{\beta h(t)}{\lambda}-1)\right|.$$
Then, the generator $f$ satisfies
$$-z^{tr}\theta_t-\frac{1}{2\alpha}|\theta_t|^2-\beta h(t)y-h^{\star}\le
f(t,y,z)\leq \frac{\alpha}{2}|z|^2-\beta h(t) y+h^{\star}.$$
By truncating the terminal data $F$ as in Theorem \ref{theorem:BSDE1} by approximating it with $F^{n,k}$, and applying the comparison theorem for quadratic BSDE with bounded terminal data to BSDE$(F^{n,k},f)$, we obtain 
that $Y^{n,k}$ is increasing in $n$ and decreasing in $k$, with upper and lower bounds as
$\underline{Y}_t\leq Y_t^{n,k} \leq \overline{Y}_t$,
where $\underline{Y}$ solves
BSDE$(-F^{-},-z^{tr}\theta_t-\frac{1}{2\alpha}|\theta_t|^2-\beta h(t)y-h^{\star})$, namely,
\begin{equation}\label{BSDElconsump*}
	\underline{Y}_t=-F^-+\int_t^T
	(-\underline{Z}_s^{tr}\theta_s-\frac{1}{2\alpha}|\theta_s|^2-\beta h(s)\underline{Y}_s-h^{\star}
	)ds-\int_t^T \underline{Z}_s^{tr}dB_s,\ t\in[0,T],
\end{equation}
and $\overline{Y}$ solves BSDE$(F^{+},\frac{\alpha}{2}|z|^2-\beta h(t)y+h^{\star})$, namely,
\begin{equation}\label{BSDE2consump*}
	\overline{Y}_t=F^{+}+\int_t^T
	(\frac{\alpha}{2}|\overline{Z}_s|^2-\beta h(s)\overline{Y}_s+ h^{\star})ds-\int_t^T
	\overline{Z}_s^{tr}dB_s,\ t\in[0,T].
\end{equation}

It is routine to check that both $\underline{Y}$ and $\overline{Y}$ can be further bounded as

\begin{align*}
	\underline{Y}_t&\geq E^{\mathbb{Q}^{\theta}}\left[\left.-e^{-\int_t^T\beta h(s)ds}F^{-}-\int_t^Te^{-\int_t^s\beta h(s)ds}
	(\frac{1}{2\alpha}|\theta_s|^2+h^{\star})ds\right|{\mathcal F}_t\right]\\
&\geq -E^{\mathbb{Q}^{\theta}}\left[F^{-}|{\mathcal F}_t\right]+C;\\
	\overline{Y}_t&\leq \frac{1}{\alpha}\ln E\left[\left.\exp\left(\alpha e^{-\int_t^T\beta h(s)ds} F^{+}+\alpha \int_t^Te^{-\int_t^T\beta h(s)ds}h^{\star}ds\right)\right|\mathcal{F}_t\right]\\
&\leq \frac{1}{\alpha}\ln E[e^{\alpha F^+}|\mathcal{F}_t]+C,
\end{align*}
for some constant $C>0$ independent of $n$ and $k$, where $\mathbb{Q}^{\theta}$ is the minimal local martingale measure given in (\ref{EMM}).
The rest of the proof proceeds using the analogous localization argument as presented in the proof of Theorem \ref{theorem:BSDE1}. Specifically, we employ the localization sequence \( (\tau_j)_{j \geq 1} \), defined by
\[
\tau_j = T \wedge \inf\left\{t \in [0, T] : \max\left\{\frac{1}{\alpha}\ln E[e^{\alpha F^+}|\mathcal{F}_t]+C, E^{\mathbb{Q}^{\theta}}\left[F^{-}|{\mathcal F}_t\right]-C\right\} > j \right\},
\]
As the remaining steps closely mirror those in the proof of Theorem \ref{theorem:BSDE1}, we omit the detailed exposition for brevity.
\end{proof}

\begin{remark} Since the generator $f$ is decreasing and Lipschitz continuous in $y$, as shown from the proof, including $y$ in the generator has no impact on the exponential integrability of the solution component $Y$.
\end{remark}

\section{Optimal consumption and investment with unbounded payoffs}

We are in a position to solve the optimization problem (\ref{expoopt_cons}) via the quadratic BSDE$(F,f)$ with generator (\ref{driver_cons}).

\begin{theorem}\label{thm:optimal_cons}
	Suppose that Assumption 1 holds. Let $(Y,Z)$ be the unique solution
	to BSDE$(F,f)$ with generator (\ref{driver_cons}). Then, the value function of
	the optimization problem (\ref{expoopt_cons}) with admissible set
	$\mathcal{A}_D^{C}$ is given by
	\begin{equation}\label{valuefunction_cons}
		V(0, x)= -\exp(-\alpha(h(0)x-Y_0)),
	\end{equation}
	and there exists an optimal trading strategy $\pi^{\star}$ and an optimal consumption process $C_t^{\star}$ (with the feedback form depending on the wealth level $X^{\pi}_t$), both in the admissible set $\mathcal{A}_D^C$, given by
\begin{align}\label{pi_cons}
h(t)\sigma_t^{tr}\pi^{\star}_t&=\mbox{Proj}_{h(t)\sigma_t^{tr}\mathcal{C}}\left(Z_t+\frac{\theta_t}{\alpha}\right);\\
C_t^{\star}&=h(t)X_t^{\pi, C^{\star}}-Y_t-\frac{1}{\alpha}\ln\frac{\beta h(t)}{\lambda},\ \hbox{$\mathbb{P}$-a.s.,}\ \text{for}\ t \in	[0,T].\notag
\end{align}
\end{theorem}

\begin{proof} We construct the conditional valueprocess as
 \begin{equation}\label{value_process_consump}
 V(t,X_t^{\pi,C})=-e^{-\alpha (h(t)X_t^{\pi,C}-Y_t)},
 \end{equation}
for $t\in[0,T]$, and verify the martingale optimality principle, i.e. $V(t,X_t^{\pi,C})-\int_0^t\lambda e^{-\alpha C_s}ds$, $t\in[0,T]$,
is
	a supermartingale for any admissible $(\pi, C)$, and is a
	martingale for $(\pi^*,C^{\star})$ given in (\ref{pi_cons}). Consequently,
\begin{align*}
-\exp(-\alpha(h(0)x-Y_0))&=\sup_{(\pi,C)\in\mathcal{A}_D^C}E\left[-\exp(-\alpha(X_T^{\pi,C}-F))-\int_0^T\lambda\exp(-\alpha C_u)du\right]\\
&=E\left[-\exp(-\alpha(X_T^{\pi^{\star},C^{\star}}-F))-\int_0^T\lambda\exp(-\alpha C_u^{\star})du\right].
\end{align*}

	To verify the supermartingale and martingale properties, an application of
	It\^o's formula to $-e^{-\alpha(h(t)X_{t}^{\pi,C}-Y_{t})}$, for $(\pi,C)\in\mathcal{A}_{D}^C$, gives
\begin{align}\label{value_process_consump_2}	
&V(t,X_t^{\pi,C})-\int_0^t\lambda e^{-\alpha C_s}ds\\
=&\ -e^{-\alpha(h(0)x-Y_0)}+\int_0^t\alpha e^{-\alpha (h(u)X_u^{\pi,C}-Y_u)}\left[A_u^{\pi,C}du+(h(u)\pi_u^{tr}\sigma_u-Z_u^{tr})dB_u\right],\notag
\end{align}
	where
\begin{align*}	
A_u^{\pi,C}=& -\frac{\alpha}{2}\left|h(u)\sigma_u^{tr}\pi_u-(Z_u +
	\frac{\theta_u}{\alpha})\right|^2 + Z_u^{tr}\theta_u + \frac{|\theta_u|^2}{2
		\alpha}\\
&-\beta h(u)C_u-\frac{\lambda}{\alpha}e^{\alpha (h(u)X_u^{\pi,C}-Y_u)}e^{-\alpha C_u}+\beta h^2(u)X_u^{\pi,C}+ f(u,Y_u,Z_u).
\end{align*}

We require that $A_u^{\pi,C}\leq 0$ for any admissible $(\pi,C)\in\mathcal{A}_D^{C}$ and $A_u^{\pi^{\star},C^{\star}}=0$ for $(\pi^{\star},C^{\star})$ given by (\ref{pi_cons}). Hence, the generator $f$ must take the form
\begin{align*}	
f(u,Y_u,Z_u)=& \min_{\pi_u\in\mathcal{C}}\left\{\frac{\alpha}{2}|h(u)\sigma_u^{tr}\pi_u-(Z_u +
	\frac{\theta_u}{\alpha})|^2 - Z_u^{tr}\theta_u - \frac{|\theta_u|^2}{2
		\alpha}\right\}\\
&+\min_{C_u\in\mathbb{R}}\left\{\beta h(u)C_u+\frac{{\lambda}}{\alpha}e^{\alpha (h(s)X_u^{\pi,C}-Y_u)}e^{-\alpha C_u}-\beta h^2(u)X_u^{\pi,C}\right\}.
\end{align*}
Note that the second term can be rewritten as
\begin{align*}
&-\frac{{\lambda}}{\alpha}e^{\alpha (h(s)X_u^{\pi,C}-Y_s)}\max_{C_u\in\mathbb{R}}\left\{-C_u\underbrace{\frac{\alpha\beta h(u)}{\lambda}e^{-\alpha (h(s)X_u^{\pi,C}-Y_u)}}_{P_u}-e^{-\alpha C_u}\right\}-\beta h^2(u)X_u^{\pi,C}\\
&\quad =-\frac{{\lambda}}{\alpha}e^{\alpha (h(s)X_u^{\pi,C}-Y_s)}\frac{P_u}{\alpha}(\ln\frac{P_u}{\alpha}-1)-\beta h^2(u)X_u^{\pi,C}\\
&\quad= -\beta h(u)Y_u-\frac{\beta h(u)}{\alpha}(\ln\frac{\beta h(u)}{\lambda}-1),
\end{align*}
where we utilized the convex duality between $e^{\alpha x}$ and $\frac{y}{\alpha}(\ln\frac{y}{\alpha}-1)$ in the first equality. This lead to the generator $f$ as defined in (\ref{driver_cons}).

We next show that $V(t,X_t^{\pi^{\star},C^{\star}})$, $t\in[0,T]$, satisfies the Class (D) condition.
From the expression (\ref{value_process_consump}), it is sufficient to verify that $e^{-\alpha (h(t)X_t^{\pi^{\star},C^{\star}}-Y_t)},$ $t\in[0,T]$, satisfies the Class (D) condition. To this end, note that (\ref{value_process_consump_2}) yields that, with $(\pi^{\star},C^{\star})$ given in (\ref{pi_cons}),
\begin{align*}
de^{-\alpha (h(t)X_t^{\pi^{\star},C^{\star}}-Y_t)}=&\ e^{-\alpha (h(t)X_t^{\pi^{\star},C^{\star}}-Y_t)}\\
&\ \times\left[\alpha(Z_t^{tr}-h(t)(\pi_t^ {\star})^{tr}\sigma_t)dB_t-\lambda e^{\alpha (h(t)X_t^{\pi^{\star},C^{\star}}-Y_t)}e^{-\alpha C_t^{\star}}dt\right].
\end{align*}
In turn,
\begin{equation}\label{value_function_equ}
e^{-\alpha (h(t)X_t^{\pi^{\star},C^{\star}}-Y_t)}=e^{-\alpha (h(0)x-Y_0)}e^{-\int_0^t\lambda e^{\alpha (h(u)X_u^{\pi^{\star},C^{\star}}-Y_u)}e^{-\alpha C_u^{\star}}du}L_t^{\pi^{\star}},
\end{equation}
where
\begin{equation*} L_{t}^{\pi}=\mathcal{E}_t\left(\alpha\int_0^{\cdot}(Z_u^{tr}-h(u)\pi_u^{tr}\sigma_u)dB_u\right).
\end{equation*}
Note that
\begin{equation}\label{inequality_for_consumption}
e^{-\int_0^t\lambda e^{\alpha (h(u)X_u^{\pi^{\star},C^{\star}}-Y_u)}e^{-\alpha C_u^{\star}}du}\leq 1,
\end{equation}
so it is sufficient to prove that the optimal density process $L^{\pi^{\star}}$ is in Class (D). Similar to (\ref{subdifferential}) in Lemma \ref{lemma_convex}, it can be shown that  $$\alpha(Z_t-h(t)\sigma_t^{tr}\pi_t^{\star})\in\partial_z f(t,Y_t,Z_t).$$  Hence,  it is sufficient to prove that $L_T^{q^{\star}}$ has finite entropy for any $q^{\star}_t\in\partial_z f(t,Y_t,Z_t)$, where
\begin{equation*}		L_{t}^{q^{\star}}=\mathcal{E}_t\left(\int_0^{\cdot}q_u^{\star}dB_u\right),\quad t\in[0,T],
	\end{equation*}
which will be shown in Lemma \ref{lemma_consumption} after the proof. 

We are left to show the integrability condition (\ref{consumption_integrability_condtion}), so that 
$(\pi^{\star},C^{\star})\in\mathcal{A}_D^{C}$. Indeed, by the expression of $C^{\star}$ in (\ref{pi_cons}) and the equality (\ref{value_function_equ}) for the conditional value process, we have
\begin{align*}
E[\exp\left(-\alpha C_t^{\star}\right)]&=
E\left[e^{\left(-\alpha (h(t)X_t^{\pi^{\star}, C^{\star}}-Y_t)\right)}\right]\frac{\beta h(t)}{\lambda}\\
&=E\left[e^{-\int_0^t\lambda e^{\alpha (h(u)X_u^{\pi^{\star},C^{\star}}-Y_u)}e^{-\alpha C_u^{\star}}du}L_t^{\pi^{\star}}\right]e^{-\alpha (h(0)x-Y_0)}\frac{\beta h(t)}{\lambda}.
\end{align*}
Note that the above expectation is less than $1$ thanks to (\ref{inequality_for_consumption}) and the Class (D) property of $L^{\pi^{\star}}$. Hence,
\begin{align*}
 E[\int_0^T\exp\left(-\alpha C_u^{\star}\right)du]&\leq \int_0^Te^{-\alpha (h(0)x-Y_0)}\frac{\beta h(t)}{\lambda}dt\\
 &=e^{-\alpha (h(0)x-Y_0)}\frac{T(\ln(1+T))}{\lambda},
\end{align*}
from which we conclude the proof. 
\end{proof}

\begin{lemma}\label{lemma_consumption}
	The optimal density process $L_T^{q^{\star}}$ has finite entropy.
Hence, by De la Vall\'ee-Poussin theorem, $L^{q^{\star}}$ is in
	Class (D) and, therefore, it is a uniformly integrable martingale. 
 \end{lemma}

\begin{proof} The proof is similar to Lemma \ref{lemma}, so we only highlight the main differences. With a localization sequence $(\sigma_j)_{j\geq 1}$, we have
	\begin{equation}\label{inequality_consump}
		E[L_{\sigma_{j}}^{q^{\star}}\alpha Y_{\sigma_j}]\leq
		\frac{E[L_{\sigma_j}^{q^{\star}}\ln L_{\sigma_j}^{q^{\star}}]}{p} -
		\frac{E[L_{\sigma_j}^{q^{\star}}]\ln p}{p}+E[e^{p\alpha
			Y^+_{\star}}]
	\end{equation}
by the Fenchel inequality.
	
	Rewriting BSDE$(F,f)$ under $\mathbb{Q}^{q^{\star}}$
gives
$$dY_t=-(f(t,Y_t,Z_t)-Z_t^{tr}q_t^{\star})dt+Z_t^{tr}dB_t^{q^{\star}},$$
and moreover,
$$d\left(e^{-\int_0^t\beta h(s)ds}Y_t\right)=-e^{-\int_0^t\beta h(s)ds}(f(t,Y_t,Z_t)+\beta h(t)Y_t-Z_t^{tr}q_t^{\star})dt+Z_t^{tr}dB_t^{\star}.$$
Hence,
	\begin{align*}
E[L_{\sigma_{j}}^{q^{\star}}\alpha
		Y_{\sigma_j}] &= E^{\mathbb{Q}^{q^{\star}}}[\alpha
		Y_{\sigma_{j}}]\\		
&=E^{\mathbb{Q}^{q^{\star}}}\left[e^{\int_0^{\sigma_j}\beta h(s)ds}\alpha
		Y_0\right.\\
&\quad\left. +\alpha\int_0^{\sigma_{j}}e^{\int_u^{\sigma^{j}}\beta h(s)ds}
(-f(u,Y_u,Z_u)-\beta h(u)Y_u+Z_u^{tr}q_u^{\star})du\right]\notag\\
	\end{align*}
Define $f^{\star}$ as the convex dual of $f$:
\begin{equation}\label{dual_consump}
	{f}^{\star}(t,y,q):=\sup_{z\in\mathbb{R}^m}\left({z}^{tr}q-f(t,y,z)\right),
\end{equation}
which is linear in $y$ and has lower bound
$$f^{\star}(t,y,q)\geq \frac{1}{2\alpha}|q|^2+
 \beta h(t)y+\frac{\beta h(t)}{\alpha}(\ln\frac{\beta h(t)}{\lambda}-1).$$
Then,
\begin{align}\label{BSDE_inequality_consump}
&E[L_{\sigma_{j}}^{q^{\star}}\alpha
		Y_{\sigma_j}] \notag\\
		=\ &E^{\mathbb{Q}^{q^{\star}}}\left[e^{\int_0^{\sigma_j}\beta h(s)ds}\alpha
		Y_0+\alpha\int_0^{\sigma_{j}}e^{\int_u^{\sigma^{j}}\beta h(s)ds}
(f^{\star}(u,Y_u,q^{\star}_u)-\beta h(u)Y_u)du\right]\notag\\
		&\geq -e^{T}\alpha |Y_0|\\
		&+E^{\mathbb{Q}^{q^{\star}}}\left[
	\alpha\int_0^{\sigma_{j}}e^{\int_u^{\sigma_j}\beta h(s)ds}\left(\frac{1}{2\alpha}|q_u^{\star}|^2+\frac{\beta h(u)}{\alpha}(\ln\frac{\beta h(u)}{\lambda}-1)\right)du\right],\notag\\
&\geq -e^{T}\alpha|Y_0|+E^{\mathbb{Q}^{q^{\star}}}\left[\int_0^{\sigma_{j}}\frac{1}{2}|q_u^{\star}|^2du\right]-e^{T}\alpha h^{\star}T.\notag
	\end{align}

	Finally, combining (\ref{inequality_consump}) and (\ref{BSDE_inequality_consump}), and
	observing
$$	
E[L_{\sigma_j}^{q^{\star}}\ln
L_{\sigma_j}^{q^{\star}}]=E^{\mathbb{Q}^{q^{\star}}}\left[\int_0^{\sigma_{j}}\frac12|q_u^{\star}|^2du\right],
$$
	we obtain
	$$(1-\frac{1}{p})E[L_{\sigma_j}^{q^{\star}}\ln
	L_{\sigma_j}^{q^{\star}}]\leq e^T\alpha |Y_0|+e^T\alpha h^{\star}T-\frac{\ln
		p}{p}+E[e^{p\alpha Y^+_{\star}}]<+\infty.$$ The assertion then
	follows by sending $j\rightarrow\infty$ in the above inequality, and
	using Fatou's lemma.
\end{proof}

\section[A convex dual representation]{A convex dual representation of the consumption and investment model}\label{Section_convex_consump}

In this section, we provide a convex dual representation of
the solution component $Y$ of BSDE$(F,f)$ with the generator $f(t,y,z)$ given by (\ref{driver_cons}) in this section, which completes the
proof of Theorem \ref{theorem:BSDE_cons} for its solution uniqueness.

We first introduce the admissible set of the convex dual problem. For
an $\mathbb{R}^m$-valued $\mathbb{F}$-predictable process $q$,
we
define the stochastic exponential$$
L^{q}:=\mathcal{E}_t\left(\int_0^{\cdot}q_u^{tr}dB_u\right). $$
If $L_T^{q}$ has a finite entropy,  i.e. $E[L_T^{q}\ln L_T^{q}]<+\infty$, then De la
Vall\'ee-Poussin theorem implies that $L^{q}$ is in Class (D) and is
therefore a uniformly integrable martingale. We then define the
probability measure $\mathbb{Q}^{q}$ on $\mathcal{F}_T$ by
$d\mathbb{Q}^{q}:=L_T^{q}d\mathbb{P}$, and introduce the
admissible set
\begin{align*}
	\mathcal{A}_D^{C,\star}\ :=\ \Biggl\{&q\in L^2[0,T]:\
	L_T^{q}\ \mbox{has a finite
		entropy such}\\
	&\ \text{that }  E^{\mathbb{Q}^q}\left[|F|+\int_0^T|f^{\star}(s,0,q_s)|ds\right]<\infty\ \mbox{with}\ {d\mathbb{Q}^{q}}:=L_T^{q}{d\mathbb{P}}\ \Biggr\}
	\text{.}
\end{align*}%

\begin{theorem}\label{theorem_dual_consump}
	Suppose that Assumption \ref{assumption1} holds. Then, the solution component
	$Y$ to BSDE$(F,f)$ with $f(t,y,z)$ given by (\ref{driver_cons}) admits the following
	convex dual representation
	\begin{equation}\label{dual_formula_consump}
	Y_t=\esssup_{q\in\mathcal{A}_D^{C,\star}}Y_t^{q}, \quad t\in [0,T], 
	\end{equation}
where $Y^q$ is the unique solution to BSDE$(F,-f^{\star})$ with the generator $f^{\star}$ given by (\ref{dual_consump}):
\begin{equation}\label{dual_BSDE_consump}
Y_t^q=F-\int_t^Tf^{\star}(s,Y_s^q,q_s)ds-\int_t^T(Z_s^{q})^{tr}dB_s^q, \quad t\in [0,T].
\end{equation}
Moreover, there exists an
	optimal density process $q^{\star}\in\mathcal{A}_D^{C,\star}$
	such that $Y=Y^{q^{\star}}$.
\end{theorem}

\begin{proof} Since the term involving \( y \) in the generator \( -f^*(s, y, q_s) \) is \( -\beta h(s) y \), which is linear and decreasing in \( y \), this ensures that \(\text{BSDE}(F, -f^*)\) has a unique solution \( (Y^q, Z^q) \) for any \( q \in \mathcal{A}_D^{C,\star} \). On the other hand, BSDE$(F,f)$ reads
\begin{equation}\label{dual_BSDE_consump_2}
Y_t=F+\int_t^T\left[f(s,Y_s,Z_s)-Z_s^{tr}q_s\right]ds-\int_t^TZ_s^{tr}dB_s^{q}, \quad t\in [0,T],
\end{equation}
for any \( q \in \mathcal{A}_D^{C,\star} \), where \( B^q_t := B_t - \int_0^{t} q_u du \), \( t \in [0,T] \), is an \( m \)-dimensional Brownian motion under \( \mathbb{Q}^q \) defined at the beginning of this subsection.

With a localization sequence $
(\sigma_{j})_{j\geq 1}$,
combining (\ref{dual_BSDE_consump}) and (\ref{dual_BSDE_consump_2}) yields
{
\begin{align*}
	&Y_{t}-Y_t^{q}\\
=&\ Y_{\sigma_{j}}-Y_{\sigma_{j}}^{q}+
\int_t^{\sigma_j}\left[f(s,Y_s,Z_{s})-Z_s^{tr}q_s+{f}^{\star}(s,Y_s^{q},q_s)\right]ds\\
&\ -\int_t^{\sigma_j}(Z_s-Z_s^q)^{tr}dB_s^q\\[2mm]
=&\ e^{-\int_t^{\sigma_j}\beta h(u)du}(Y_{\sigma_{j}}-Y_{\sigma_{j}}^{q})\\
&\ +\int_t^{\sigma_j}e^{-\int_t^{s}\beta h(u)du}\left[f(s,0,Z_{s})-Z_s^{tr}q_s+{f}^{\star}(s,0,q_s)\right]ds\\
&\ -\int_t^{\sigma_j}e^{-\int_t^{s}\beta h(u)du}(Z_s-Z_s^q)^{tr}dB_s^{q}\\[2mm]
=&\ E^{\mathbb{Q}^q}\biggl[e^{-\int_t^{\sigma_j}\beta h(u)du}(Y_{\sigma_{j}}-Y_{\sigma_{j}}^{q})\\
&\ +\int_t^{\sigma_j}e^{-\int_t^{s}\beta h(u)du}\left[f(s,0,Z_{s})-Z_s^{tr}q_s+{f}^{\star}(s,0,q_s)\right]ds\Bigm|\mathcal{F}_t\biggr].
\end{align*}
}

By the Fenchel-Moreau theorem, we then deduce that
\begin{align*}
	f(s,0,Z_{s})-Z_s^{tr}q_s+{f}^{\star}(s,0,q_s)&\geq 0, \quad \forall q\in\mathcal{A}_D^{C,\star}
\end{align*}
and the equality holds
for $q^{\star}_s\in\partial_z f(s,0,Z_s)$.
Hence,
$$
Y_t-Y_t^{q}\geq E^{\mathbb{Q}^{q}}\left[e^{-\int_t^{\sigma_j}\beta h(u)du}(Y_{\sigma_j}-Y_{\sigma_j}^{q})\Bigm| \mathcal{F}_t\right], 
\quad \forall q\in\mathcal{A}_D^{C,\star}
$$ 
and 
$$Y_t-Y_t^{q^{\star}}= E^{\mathbb{Q}^{q^{\star}}}\left[e^{-\int_t^{\sigma_j}\beta h(u)du}(Y_{\sigma_j}-Y_{\sigma_j}^{q^{\star}})\Bigm|\mathcal{F}_t\right].$$
Note that $e^{-\int_t^{\sigma_j}\beta h(u)du}\leq 1$, and both $Y$ and $Y^{q}$ are uniformly integrable under $\mathbb{Q}^q$. Sending $j\rightarrow \infty$, we obtain that $Y_t\geq Y_t^{q}$ for any $q\in\mathcal{A}_D^{C,\star}$.

On the other hand, $Y_t=Y_t^{q^{\star}}$. To achieve this, it
remains to prove $q^{\star}\in\mathcal{A}_D^{C,\star}$.
Indeed, the finite entropy of $L_{T}^{q^{\star}}$ has already been established in Lemma \ref{lemma_consumption}. Moreover, we can prove  $q^{\star}\in L^{2}[0,T]$ and
$$E^{\mathbb{Q}^{q^{\star}}}\left[|F|+\int_0^T|f^{\star}(s,0,q^{\star}_s)|ds\right]<+\infty$$
with a similar approach to the proof in Lemma \ref{lemma000}, and as such, the details are omitted.
\end{proof}

%
%
%


\chapter[Epstein-Zin recursive utility maximization]{Epstein-Zin recursive utility maximization in an unbounded market}\label{Chapter8}
In this chapter, we examine an optimal investment and consumption problem within the framework of an investor with Epstein-Zin recursive utility. A notable departure from earlier models is the unbounded nature of the market coefficients, particularly in scenarios where the market price of risk $\theta$ is unbounded.

\section{Epsetin-Zin recursive utility}

The Epstein-Zin recursive utility, initially introduced in \cite{EpsteinZin}, represents a generalization of the traditional expected power utility. It allows for the separation of time aggregation and risk aggregation, thereby permitting the capturing of intertemporal substitution and risk aversion within a more sophisticated and general framework. The Epstein-Zin recursive utility is governed by the Epstein-Zin aggregator, represented as
\begin{equation}\label{EZ_aggregator}
f^{EZ}(C,V)=\delta\frac{(1-\gamma)V}{1-1/\psi}
\left[\left(\frac{C}{[(1-\gamma)V]^{\frac{1}{1-\gamma}}}\right)^{1-1/\psi}-1\right]
\end{equation}
for a given consumption level $C\geq 0$ and a utility $V\in\mathbb{R}$ such that $(1-\gamma)V\geq 0$\footnote{We follow the notation convention in the Epstein-Zin literature, where capital letters $C$ and $V$ are used as variables to represent consumption and utility, respectively.},
where $0<\gamma\neq 1$ is the relative risk aversion, $0<\psi\neq 1$ is the elasticity of intertemporal substitution (EIS) and $\delta>0$ is the discounting rate.

A special case arises when $\gamma=1/\psi$ resulting in the intertwining of risk aversion and intertemporal substitution. Consequently, the Epstein-Zin aggregator reduces to the conventional power utility aggregator:
$$f^{EZ}(C,V)=\delta\frac{C^{1-\gamma}}{1-\gamma}-\delta V.$$

Given the Epstein-Zin aggregator $f^{EZ}$, the Epstein-Zin recursive utility of consumption $C$ is represented as
$$V_t^{C}=E\left[\left.\int_t^Tf^{EZ}(C_s,V_s^C)ds+U(C_T)\right|\mathcal{F}_t\right],$$
where $U(C)=\frac{C^{1-\gamma}}{1-\gamma}$ represents the terminal utility of consumption. When $\gamma=1/\psi$, then $V_t^C$ reduces to the traditional expected power utility of consumption $C$:
$$V_t^C=E\left[\left.\int_t^T\delta e^{-\delta(s-t)}\frac{C_s^{1-\gamma}}{1-\gamma}ds+e^{-\delta(T-t)}U(C_T)\right|\mathcal{F}_t\right].$$

It is more convenient to introduce an additional parameter that measures the ratio between risk aversion and intertemporal substitution:
\begin{equation}\label{beta}
\beta=\frac{1-\gamma}{1-1/\psi},
\end{equation}
and the Epstein-Zin aggregator can be rewritten as
\begin{equation}\label{EZ_aggregator}
f^{EZ}(C,V)=\delta\frac{C^{1-\frac{1}{\psi}}}{1-\frac{1}{\psi}}[(1-\gamma)V]^{1-\frac{1}{\beta}}-\delta\beta V.
\end{equation}

We consider an optimal investment and consumption problem utilizing the Epstein-Zin recursive utility within the incomplete financial market introduced in Chapter \ref{Chapter2}. One of the key distinctions lies in allowing the market price of risk $\theta$ to be \emph{unbounded} in this chapter. As we shall demonstrate, this optimization problem leads to quadratic BSDEs with unbounded coefficients. Similar optimization problems have received considerable attention recently, particularly in \cite{Kraft2017,Xing2017,Feng2023}. The latter two are motivated by the quadratic BSDE representation for the expected power utility derived in Cheridito and Hu \cite{CH}. Additionally, see \cite{Matoussi2018} for a convex duality representation, and \cite{Huang2023,Herdegen2023(1),Herdegen2023(2)} for extensions to an infinite horizon framework. Notably, the aforementioned works do not cover portfolio constraints, which constitutes one of the primary contributions of our work.

Let \( u_t \) represent the \textit{proportion} of wealth invested in risky assets at time \( t \), and \( c_t \) represent the \textit{proportion} of wealth allocated for consumption at \( t \). The wealth equation with consumption in (\ref{wealth_cons}) is modified as
\begin{equation}\label{wealth_2}
X_t^{u,c} = x + \int_0^t X_s^{u,c}\left[u_s^{\text{tr}}\sigma_s\left(dB_s + \theta_s ds\right) - c_s ds\right], \quad x > 0,
\end{equation}
with initial wealth \( x \). Note that \( \pi_t = X_t^{u,c} u_t \) and \( C_t = X_t^{u,c} c_t \), where \( \pi_t \) is the amount of money invested in risky assets, and \( C_t \) is the amount of money allocated for consumption. Here, \( \theta \) is the \textit{unbounded} market price of risk, which is assumed to satisfy the Novikov's condition 
\begin{equation}\label{Novikov}
E\left[\exp\left(\frac{1}{2} \int_0^T |\theta_t|^2 dt\right)\right] < \infty.
\end{equation}
As in (\ref{EMM}), we can therefore define an equivalent minimal local martingale measure (MLMM) \( \mathbb{Q}^{\theta} \) on \( \mathcal{F}_T \) by
\begin{equation}\label{EMM2}
    \frac{d\mathbb{Q}^{\theta}}{d\mathbb{P}} := L_T^{\theta} = \mathcal{E}_T\left( -\int_0^{\cdot} \theta_t^{\text{tr}} dB_t \right),
\end{equation}
where \( \mathcal{E}(\cdot) \) denotes the stochastic exponential. Under the MLMM \( \mathbb{Q}^{\theta} \), the wealth process becomes a nonnegative local supermartingale. Therefore, it follows that \( X_t^{u,c} \) is a supermartingale. 
 
We consider the Epstein-Zin recursive utility of consumption and wealth governed by the Epstein-Zin aggregator $f^{EZ}$ and the terminal power utility $U$ of bequest $X_T^{u, c}$, which plays the role of reward functional: 
\begin{equation}\label{EZ_BSDE}
V_t^{\xi;u,c}=U(X_T^{\xi;u,c})+\int_t^Tf^{EZ}(X_s^{\xi;u,c}c_s,V_s^{\xi;u,c})ds-\int_t^T(Z_s^{\xi;u,c})^{tr}dB_s,
\end{equation}
where 
\begin{equation*}
X_r^{\xi;u,c}=\xi+\int_t^rX_s^{\xi;u,c}\left[u_s^{tr}\sigma_s(dB_s+\theta_sds)-c_sds\right],\quad r\in[t,T],
\end{equation*}
for any random variable $\xi\in\mathcal{F}_t$ such that $U(\xi)$ is integrable.

We present the well-posedness result for the Epstein-Zin recursive utility. The main difficulty is that the Epstein-Zin aggregator \(f^{EZ}(C, V)\) is unbounded, polynomial in \(V\), and singular at \(V=0\). For technical reasons, we only consider EIS \(\psi > 1\). See the discussions in Remark \ref{remark_difficulty} for the reasons.

\begin{lemma}\label{lemma_EZ} 
Let $(V_t^{u,c},Z_t^{u,c}):=(V^{X_t^{u,,c},u,c}_t,Z^{X_t^{u,c},u,c}_t)$. 
Suppose that $\psi>1$, and that $(u,c)$ satisfy the following integrability conditions: 

(i) for $\gamma>1$, 
\begin{equation}\label{integrability_condition1}  
E[(X^{u,c}_T)^{1-\gamma}]<\infty,\quad E[\int_0^TC_t^{1-\frac{1}{\psi}}dt]<\infty;
\end{equation}

(ii) for $\gamma<1$, 
\begin{equation}\label{integrability_condition2}  
E[(X^{u,c}_T)^{1-\gamma\wedge\frac{1}{\psi}}]<\infty,\quad E[\int_0^TC_t^{1-\gamma\wedge\frac{1}{\psi}}dt]<\infty,
\end{equation}
where $X_T^{u,c}$ is given by (\ref{wealth_2}) and $C_t=X_t^{u,c}c_t$.

Then, the Epstein-Zin BSDE (\ref{EZ_BSDE}) admits a unique solution\footnote{For ${\gamma<\frac{1}{\psi}}$, by the unique solution, we mean its minimal solution. Recall $(V^{u,c},Z^{u,c})$ is the minimal solution to  (\ref{EZ_BSDE}), if for any other solution $(\bar{V}^{u,c},\bar{Z}^{u,c})$ satisfying  (\ref{EZ_BSDE}), it holds that $V^{u,c}\leq \bar{V}^{u,c}$. Also refer to Theorem \ref{thm:asy} for the definition of the minimal solution.} $(V^{u,c},Z^{u,c})$ such that $V^{u,c}$ is continuous, $(1-\gamma)V^{u.c}> 0$ and of Class (D), and $Z^{u,c}\in L^2(0,T)$. Moreover, the comparison theorem holds in the sense that if ${U}(X_T^{\bar{u},\bar{c}})\geq U(X_T^{u,c})$,  $\mathbb{P}$-a.s., then their corresponding solutions satisfy ${V}^{\bar{u},\bar{c}}_t\geq V^{u,c}_t$, $\mathbb{P}$-a.s., for $t\in[0,T]$.
\end{lemma}

\begin{proof} Case (i) that $\gamma>1$ has been proved in Proposition 2.2 in Xing \cite{Xing2017}. Hence, we focus on Case (ii) that $\gamma<1$. As in \cite{Xing2017}, we first make the transformation $(P^{u,c},Q^{u,c}):=e^{-\delta \beta t}(1-\gamma) (V^{u,c},Z^{u,c})$. Then
(\ref{EZ_BSDE}) is transformed to 
\begin{equation}\label{EZ_BSDE_1}
    P_t^{u,c}=\eta+\int_t^TF^{EZ}(s,X_s^{u,c}c_s, P_s^{u,c})ds-\int_t^T(Q_s^{u,c})^{tr}dB_s,
\end{equation}
where $\eta:=e^{\delta\beta T}(X_T^{u,c})^{1-\gamma}$, and 
$$F^{EZ}(t, C, V):=\delta\beta e^{-\delta t}C^{1-\frac{1}{\psi}}V^{1-\frac{1}{\beta}},\quad t\in[0,T], C\geq 0, V\geq 0.$$
Recall $C_s=X_s^{u,c}c_s$, so we will use $C_s$ in the rest of proof. Also for notational simplicity, we omit $(u,c)$ in the rest of the proof. 

We consider $\gamma<\frac{1}{\psi}$ and $\gamma>\frac{1}{\psi}$ separately. 

For $\underline{\gamma<\frac{1}{\psi}}$, so $\beta>1$ and $1-\frac{1}{\beta}\in(0,1)$, $F^{EZ}(t,C,V)$ is increasing in $C$ and $V$, and locally Lipschitz continuous in $V$ with its Lipschitz constant exploding as $V\rightarrow 0$. Since
\[
y^{1 - \frac{1}{\beta}} \leq 1_{y \leq 1} + 1_{y > 1} y \leq 1 + y \quad \text{for } y \geq 0,
\]
we may approximate it by inf convolution as in Remark \ref{remark_approx}: For integers $n\geq 1$,
\begin{equation*}
f^n(y)=\inf_{q}\{q^{1-\frac{1}{\beta}}+n|q-y|\}.
\end{equation*}
Then, $f^n(y)$ is well defined and globally Lipschitz continuous with constant $n$. 
Moreover, $f^n(y)$ is increasing and converges pointwisely to $y^{1-\frac{1}{\beta}}$. In turn, Dini's theorem implies that convergence is also uniform on compact sets. 
We consider the truncated BSDE$(\eta\wedge n, F^{EZ}_{n})$:  \begin{equation}\label{EZ_BSDE11}
    P_t^{n}=\eta\wedge n +\int_t^TF^{EZ}_{n}(s,C_s, P_s^{n})ds-\int_t^T(Q_s^{n})^{tr}dB_s,
\end{equation}
where the generator
\begin{equation*}
F^{EZ}_{n}(t,C,V)=\delta\beta e^{-\delta t}(C\wedge n)^{1-\frac{1}{\psi}}f^n(V).
\end{equation*}
is globally Lipschitz continuous in $V$:
$$|F^{EZ}_{n}(t,C,V)-F^{EZ}_{n}(t,C,V')|\leq \delta\beta n^{2-\frac{1}{\psi}}|V-V'|.$$
Moreover,
we have 
\[
0 \leq F^{EZ}_{n}(t,C,V) \leq \delta \beta n^{1 - \frac{1}{\psi}} (1 + V).
\]
Note that since $\psi>1$, the above Lipschitz constant and the bound will explode when $n\rightarrow\infty$.

The truncated BSDE (\ref{EZ_BSDE11}) admits a unique solution $(P^{n},Q^{n})$. Furthermore, since $F^{EZ}_{n}(t,C,V)$ is increasing in $n$, the comparison theorem for Lipschitz BSDE applied to (\ref{EZ_BSDE11}) implies that $P^{n}$ is also increasing in $n$. 

Note that the component \(P^{n}\) is nonnegative, which follows from a comparison between the truncated BSDE\((\eta\wedge n, F^{EZ}_{n})\) and the BSDE\((0, F^{EZ}_{n})\), the latter clearly having the unique solution \((0,0)\). On the other hand, using $ab\leq \frac{a^p}{p}+\frac{b^q}{q}$ for $p,q>1$ and $\frac{1}{p}+\frac{1}{q}=1$, we have 
$$C^{1-\frac{1}{\psi}}V^{1-\frac{1}{\beta}}\leq \frac{1}{\beta}C^{1-\gamma}+(1-\frac{1}{\beta})V,$$
so that $$F^{EZ}_{n}(t,C,V)\leq \delta\beta\left(\frac{1}{\beta}C^{1-\gamma}+(1-\frac{1}{\beta})V\right).$$
Hence, the comparison theorem for 
$(\ref{EZ_BSDE11})$ yields that 
$
P^n_t \leq \bar{P}_t,
$
where 
\[
\bar{P}_t = e^{\delta\beta T}X_T^{1-\gamma}  + \int_t^T \left[\delta  C_s^{1 - \gamma} + \delta (\beta-1)\bar{P}_s)\right] ds - \int_t^T (\bar{Q}_s)^{\text{tr}} dB_s,
\]
which has a unique solution given by 
\[
\bar{P}_t = E\left[e^{\delta(\beta-1)(T-t)}e^{\eta\beta T}X_T^{1-\gamma} + \int_t^T e^{\delta(\beta-1)(r-t)} \delta C_r^{1 - \gamma} \, dr \,\Bigg| \, \mathcal{F}_t\right]<\infty
\]
by the integrability condition (\ref{integrability_condition2}).

Hence, we can apply the localization argument analogous to those in 
Theorem \ref{theorem:BSDE1} by introducing the stopping time sequence \((\tau_j)_{j \geq 1}\):
\[
\tau_j = T \wedge \inf\{t \in [0,T] : \bar{P}_t > j\}.
\]
Then, the stability property of the BSDE \((\ref{EZ_BSDE11})\) (see, for example, Lemma \ref{lemma_stability}) yields that there exists a limiting pair \((P, Q)\) satisfying \((\ref{EZ_BSDE_1})\) with \(0 \leq P_t \leq \bar{P}_t\) and \(Q \in L^2\)(0,T).

It is clear that \(P\) is continuous and of Class (D). Moreover, since $\eta>0$ and $F^{EZ}\geq 0$, it follows that 
$$P_t=E\left[\left.\eta+\int_t^TF^{EZ}(s,C_s, P_s)ds\right|\mathcal{F}_t\right]\geq E[\eta|\mathcal{F}_t]>0.$$
Since we approximate \(P\) from below by an increasing sequence \(P^n\) satisfying (\ref{EZ_BSDE11}), \(P\) is the minimal solution to (\ref{EZ_BSDE_1}). 
Indeed, for any other solution \((\bar{P}, \bar{Q})\) of (\ref{EZ_BSDE_1}), by the comparison theorem for (\ref{EZ_BSDE11}), we have 
$
\bar{P}_t \geq P_t^{n},
$
which implies that 
$
\bar{P}_t \geq \lim_{n \rightarrow \infty} P_t^n = P_t.
$
Finally, the comparison theorem for (\ref{EZ_BSDE_1}) follows from the comparison theorem for (\ref{EZ_BSDE11}), which gives 
$P^{n;\bar{u},\bar{c}}_t \geq P^{n;u,c}_t,
$
where we use \((u,c)\) to emphasize the dependence of \(P^{n;u,c}\) on \((u,c)\).

For $\underline{\frac{1}{\psi}<\gamma<1}$, so
\(\beta \in (0,1)\) and \(1 - \frac{1}{\beta} < 0\), \(F^{EZ}(C,V)\) is increasing in \(C\) and decreasing in \(V\), but it becomes singular at \(V=0\), with its derivative exploding as \(V\) approaches \(0\). Hence, we need to truncate \(F^{EZ}(C,V)\) near \(V=0\).
We consider the truncated BSDE\((\eta\wedge m, F^{EZ}_{m})\): 
 \begin{equation}\label{EZ_BSDE12}
    P_t^{m}=\eta\wedge m+\int_t^TF^{EZ}_{m}(C_s, P_s^{m})ds-\int_t^T(Q_s^{m})^{tr}dB_s,
\end{equation}
where the generator
\begin{equation*}
F^{EZ}_{m}(t,C,V)=\delta\beta e^{-\delta t}(C\wedge m)^{1-\frac{1}{\psi}}(V\vee \frac{1}{m})^{1-\frac{1}{\beta}}
\end{equation*}
is globally Lipschitz continuous in $V$:
$$|F^{EZ}_{m}(t,C,V)-F^{EZ}_{m}(t,C,V')|\leq \delta (\beta-1)m^{1-\frac{1}{\psi}+\frac{1}{\beta}}|V-V'|.$$
Moreover,
we have 
\[
0 \leq F^{EZ}_{m}(t,C,V) \leq \delta \beta m^{\frac{1}{\beta}- \frac{1}{\psi}}.
\]
Note that since $\beta<1<\psi$, the above Lipschitz constant and the bound will explode when $m\rightarrow\infty$.

The truncated BSDE (\ref{EZ_BSDE12}) admits a unique solution $(P^{m},Q^{m})$. Furthermore, since $F^{EZ}_{m}(t,C,V)$ is increasing in $m$, the comparison theorem for Lipschitz BSDE applied to (\ref{EZ_BSDE12}) implies that $P^{m}$ is also increasing in $m$. In a manner similar to the case of \(\gamma<\frac{1}{\psi}\), one can show that the component \(P^{m}\) is nonnegative. 

Now, we need to find an upper bound for \(P^{m}\) that is independent of \(m\). This will allow us to take the limit as \(m \to \infty\) for \(P^{m}\).
Applying Itô's formula to \(\beta (P_t^{m})^{\frac{1}{\beta}}\) yields
\begin{align}\label{Itoformula}
d\left(\beta (P_t^{m})^{\frac{1}{\beta}}\right) =& -\left\{ 
\delta \beta e^{-\delta t} (C_t \wedge m)^{1 - \frac{1}{\psi}} 
(P^{m}_t)^{\frac{1}{\beta} - 1} \left(P^{m}_t \vee \frac{1}{m}\right)^{1 - \frac{1}{\beta}} \right. \\
& \quad + \left. \frac{1}{2} (1 - \frac{1}{\beta})(P_t^{m})^{\frac{1}{\beta} - 2} |Q_t^{m}|^2 \right\} \notag \\
& + (P_t^{m})^{\frac{1}{\beta} - 1} (Q_t^{m})^{tr} dB_t. \notag
\end{align}
Hence, by taking the conditional expectation with respect to \(\mathcal{F}_t\), we have
\begin{align*}
\beta (P_t^{m})^{\frac{1}{\beta}} = & E\left[\beta (\eta\wedge m)^{\frac{1}{\beta}} + \int_t^T \delta\beta e^{-\delta s} (C_s\wedge)^{1-\frac{1}{\psi}} (P^{m}_s)^{\frac{1}{\beta}-1}(P^{m}_s\vee\frac{1}{m})^{1-\frac{1}{\beta}} \, ds \right.\\
&\quad + \left. \frac{1}{2}(1-\frac{1}{\beta})(P_s^{m})^{\frac{1}{\beta}-2}|Q_s^{m}|^2 \, ds \,\bigg| \mathcal{F}_t \right].
\end{align*}
Since \(1 - \frac{1}{\beta} < 0\) and \(y^{\frac{1}{\beta}-1}(y\vee\frac{1}{n})^{1-\frac{1}{\beta}} \leq 1\), we further obtain
\[
\beta (P_t^{m})^{\frac{1}{\beta}} \leq E\left[\beta \eta^{\frac{1}{\beta}} + \int_t^T \delta\beta e^{-\delta s} C_s^{1-\frac{1}{\psi}} \, ds \bigg| \mathcal{F}_t\right].
\]
Since $\eta^{\frac{1}{\beta}}=e^{\delta T}X_T^{\frac{1-\gamma}{\beta}}=e^{\delta T}X_T^{1-\frac{1}{\psi}}$, it follows that
\begin{equation*}
P_t^{m} \leq \frac{1}{\beta^{\beta}} E\left[\beta e^{\delta T}X_T^{1-\frac{1}{\psi}} + \int_t^T \delta\beta e^{-\delta s} C_s^{1-\frac{1}{\psi}} \, ds \bigg| \mathcal{F}_t\right]^{\beta}:= \bar{P}_t<\infty
\end{equation*}

We can apply the localization argument analogous to those in
Theorem \ref{theorem:BSDE1} by defining the stopping time sequence \((\tau_j)_{j\geq 1}\):
\[
\tau_j = T \wedge \inf\{t \in [0,T]: \bar{P}_t > j,\  or\ E[\eta\wedge m|\mathcal{F}_t] < 1/j \}.
\]
Note that $\tau_j\rightarrow T$, since 
$$P_t^{m}\geq E[\eta\wedge m|\mathcal{F}_t] \geq 1/j>0,\quad t\leq \tau_{j},$$ which ensures (\ref{EZ_BSDE12}) is always well defined. This allows us to obtain a limiting pair \((P,Q)\) satisfying (\ref{EZ_BSDE_1}) with \(0 \leq P_t \leq \bar{P}_t\) and \(Q \in L^2(0,T)\). The remainder of the proof is similar to the case \(\beta > 1\) and is therefore omitted. 

Finally, since the generator $F^{EZ}(t,C,V)$ is decreasing in $V$, the uniqueness of the solution can be proved in a similar way as in Proposition 2.2 in Xing \cite{Xing2017}. See also the last section in \cite{bdh} for more general discussions. 

\end{proof}

\begin{remark} 
When \(\frac{1}{\psi} < \gamma < 1\), since the generator \(F^{EZ}(t,C,V)\) is not monotone decreasing in \(V\), it is unclear whether the solution is unique in the usual sense. Hence, we define the unique solution as the minimal solution. That is, \((V^{u,c},Z^{u,c})\) is the minimal solution to  (\ref{EZ_BSDE}), if for any other solution \((\bar{V}^{u,c},\bar{Z}^{u,c})\) satisfying  (\ref{EZ_BSDE}), it holds that \(V^{u,c} \leq \bar{V}^{u,c}\). When considering the Epstein-Zin utility maximization, this implies that we consider the maximin problem by maximizing \(V^{u,c}\) over \((u,v)\) among the minimal solutions \((V^{u,c}, Z^{u,c})\), resembling the robust utility maximization by considering the worst-case scenario.
\end{remark}

In this chapter, \textit{we aim to select an optimal $(u^{\star},c^{\star})\in\mathcal{A}_D^{EZ}$ such that}
\begin{equation}\label{EZ_formulation}
V(0,x)=\sup_{(u,c)\in\mathcal{A}_D^{EZ}}V_0^{x;u,c}=V_0^{x;u^{\star},c^{\star}}.
\end{equation}
The admissible set $\mathcal{A}_D^{EZ}$ is defined as follows:
\begin{definition}\label{admiss7}
	[{\bf Admissible strategies with constraints $\mathcal{A}_D^{EZ}$}]

Let $\mathcal{C}$ be a closed and convex set in $\mathbb R^{d}$ and $0\in\mathcal{C}$. 	
The admissible set $\mathcal{A}_D^{EZ}$ consists of all trading strategies $u$ and consumption processes $c\geq 0$ such that $u\in L^2(0,T)$ is predictable, satisfying $u_t\in\mathcal{C}$, $\mathbb{P}$-a.s., {for} $t \in
		[0,T]$; and $c\in L^1(0,T)$ is predictable. Moreover, they satisfy the integrability condition (\ref{integrability_condition1}) or (\ref{integrability_condition2}), and the following class (D) condition holds:
\begin{equation}\label{class_D_condition_EZ}
\{ V(\tau,X_{\tau}^{u,c}):\ \tau\ \mbox{is an $\mathbb{F}$-stopping time
		taking values in}\ [0,T]\}
\end{equation} is a uniformly integrable family, where the conditional value process $V(\cdot,\cdot)$ is defined as 
\begin{align}\label{expoopt_222_EZ}
	V\left(t,\xi\right):=& \esssup_{\substack{\{u_s,c_s\ s\in[t,
			T]\}\\\text{admissible}}}V_t^{\xi;u,c},
\end{align}
where $V^{\xi;u,c}$ is given as the first component of the unique solution of the Epstein-Zin BSDE (\ref{EZ_BSDE}). 
\end{definition}

\begin{remark}
In Theorem \ref{thm:optimal_EZ}, we will establish that the conditional value process has the structure of $V(t,\xi)=U(\xi)e^{Y_t}$, where $Y$ solves BSDE$(0,f+g)$ as specified in (\ref{BSDE_EZ}), with $f$ and $g$ given by (\ref{generator_EZ}) and (\ref{generator_EZ_2}), respectively.

Consequently, the condition of integrability is tantamount to affirming that $U(X_t^{u,c})e^{Y_t}$, $t\in[0,T]$, falls within Class (D). This condition remains unaffected by the optimization problem, which indicates that there is no circular dependency inherent in defining the admissible set.
\end{remark}

\section[Quadratic BSDE with unbounded coefficients]{Quadratic BSDE with unbounded coefficients and generator monotonic in $y$}

We consider the following quadratic BSDE with unbounded coefficients, which will be utilized to solve the optimization problem (\ref{EZ_formulation}):
\begin{equation}\label{BSDE_EZ}
Y_t=\int_t^T\left[f(s,Z_s)+g(Y_s)\right] ds-\int_t^TZ_s^{tr}dB_s,
\end{equation}
where the generator $f$ and $g$ are given as
\begin{equation}\label{generator_EZ}
f(t,z)=-\frac{1-\gamma}{2\gamma}\min_{u_t\in \mathcal{C}}\left|\gamma\sigma_t^{tr}u_t-(z+\theta_t)\right|^2+\frac{1-\gamma}{2\gamma}|z+\theta_t|^2+\frac{1}{2}|z|^2
\end{equation}
and
\begin{equation}\label{generator_EZ_2}
g(y)=\delta^{\psi}\frac{\beta}{\psi}\exp\left(-\frac{\psi}{\beta}y\right)-\delta\beta.
\end{equation}

There are  four scenarios in terms of the relative risk aversion $\gamma$ and EIS $\psi$:
(i) $\gamma>1, \psi>1$; (ii) $0<\gamma<1, \psi>1$; (iii) $\gamma>1, 0<\psi<1$; (iv) $0<\gamma<1, 0<\psi<1$.

In the literature, most works are primarily restricted within Case (i) alone (see, for instance, \cite{Xing2017} for the rationale). However, it is essential to note that this case excludes the traditional expected power utility where $\gamma=1/\psi$. We aim to tackle both Cases (i) and (ii), specifically with an EIS greater than $1$, to include the expected power utility case when $0<\gamma <1$. Additionally, we shall discuss the challenges hindering the resolution of Cases (iii) and (iv) for an EIS less than $1$.

Note that the ratio parameter introduced in (\ref{beta}) satisfies $\beta<0$ for Cases (i) and (iv), and $\beta>0$ for Cases (ii) and (iii).
It is clear that the generators $f$ and $g$ satisfy the following properties.

\begin{lemma}\label{Lemma_RDU_driver} Let $\mathcal{C}$ be a closed and convex set in $\mathbb{R}^d$ satisfying $0\in\mathcal{C}$.

(i) The generator $f$ has the following growth property:
When $\gamma>1$ (i.e. Cases (i) and (iii)),
\begin{equation}\label{EZ_generator_1}
-\frac{\gamma-1}{2}|\theta_t|^2\leq \frac{1}{2\gamma}|z|^2+\frac{1-\gamma}{\gamma}z^{tr}\theta_t+\frac{1-\gamma}{2\gamma}|\theta_t|^2\leq f(t,z)\leq \frac{1}{2}|z|^2;
\end{equation}
when $0<\gamma<1$ (i.e. Cases (ii) and (iv)),
\begin{equation}\label{EZ_generator_2}
\frac{1}{2}|z|^2\leq f(t,z)\leq \frac{1}{2\gamma}|z|^2+\frac{1-\gamma}{\gamma}z^{tr}\theta_t+\frac{1-\gamma}{2\gamma}|\theta_t|^2\leq \frac{1}{\gamma}|z|^2+\frac{1}{\gamma}|\theta_t|^2.
\end{equation}

(ii) The generator $f(t,z)$ is convex in $z$ and, moreover, 
\begin{equation}\label{derivative_generator_EZ}
(1-\gamma)\sigma_t^{tr}u_t^{\star}+z\in\partial_{z}f(t,z),
\end{equation}
where $u^{\star}$ is such that
\begin{equation}\label{proj_RDU}
\gamma\sigma_t^{tr}u_t^{\star}=\mbox{Proj}_{\gamma\sigma_t^{tr}\mathcal{C}} \left(z+\theta_t\right), \mbox{$\mathbb{P}$-a.s.}, \text{for}\ t\in [0,T].
\end{equation}

(iii) The generator $g$ is decreasing in $y$, and moreover, when $\gamma>1,\psi>1$(i.e. Case (i)) or $0<\gamma<1,0<\psi<1$ (i.e. Case (iv)), it holds that
$$g(y)\leq -\delta\beta;$$
when $0<\gamma <1, \psi>1$ (i.e. Case (ii)) or $\gamma>1,0<\psi<1$ (i.e. Case (iii)), it holds that
$$g(y)\geq -\delta\beta.$$
\end{lemma}

\begin{proof} (i) To prove the first two inequalities in (\ref{EZ_generator_1}), when $\gamma>1$, note that
\begin{align*}
f(t,z)&\geq \frac{1-\gamma}{2\gamma}|z+\theta_t|^2+\frac{1}{2}|z|^2\\
&=\frac{1}{2\gamma}|z|^2+\frac{1-\gamma}{\gamma}z^{tr}\theta_t+\frac{1-\gamma}{2\gamma}|\theta_t|^2\\
&\geq \frac{1}{2\gamma}|z|^2+\frac{1-\gamma}{\gamma}(\varepsilon|z|^2+\frac{1}{4\varepsilon}|\theta_t|^2)+\frac{1-\gamma}{2\gamma}|\theta_t|^2.
\end{align*}
The first inequality then follows by choosing $\varepsilon=\frac{1}{2(\gamma-1)}$. The rest inequalities for $f$ follow similarly.

(ii) When $\gamma>1$, the convexity of $f(t,z)$ in $z$ is evident. When $0< \gamma < 1$, note that the generator $f$ can be rewritten as
\begin{align*}
f(t,z)=\frac{1-\gamma}{2\gamma}\sup_{u_t\in\mathcal{C}}\left\{2\gamma u_t^{tr}\sigma_t(z+\theta_t)-\gamma^2|\sigma^{tr}_tu_t|^2\right\}+\frac{1}{2}|z|^2.
\end{align*}
The first term represents a supremum over linear functions of $z$ and the second term is a quadratic function of $z$, from which the convexity of $f(t, z)$ in $z$ follows. Additionally, the calculation of the subdifferential $\partial_zf(t, z)$ follows similar arguments as in (\ref{subdifferential}) of Lemma \ref{lemma_convex}.

Assertion (iii)  is immediate.
\end{proof}

\begin{remark}\label{remark_difficulty} The main challenge in solving the quadratic BSDE (\ref{BSDE_EZ}) and the optimization problem (\ref{EZ_formulation}) arises due to the unbounded nature of the market price of risk $\theta$, which leads to the quadratic BSDE (\ref{BSDE_EZ}) with \emph{unbounded coefficients}. Managing such unbounded coefficients relies on the specific parameter assumptions. On one hand, according to the growth property of $f$ in (i) of Lemma \ref{Lemma_RDU_driver}, it can be inferred that the solution component $Y$ is bounded from above, denoted as $Y_{\max}$, for Cases (i) and (iii), but unbounded from below. Consequently, if $g(y)\leq -\delta\beta$ (i.e., Cases (i) and (iv)), we can assert control over the generator $g$ using the inequality:
$$g(Y_{\max})\leq g(Y)\leq -\delta\beta.$$  This observation indicates the solvability of Case (i).

On the other hand, for Cases (ii) and (iv), the solution component $Y$ is bounded from below, denoted as $Y_{\min}$, but unbounded from above. Therefore, if $g(y)\geq -\delta\beta$ (i.e. Cases (ii) and (iii)), then we can gain control over the generator $g$ through the inequality:
$$-\delta\beta\leq g(Y)\leq g(Y_{\min}).$$
This observation suggests the solvability of Case (ii).

Finally, in Cases (iii) and (iv), we only observe a one-sided control of $g$. However, it becomes clear in the following sections that this one-sided control is inadequate for ensuring either the uniqueness of the solution to (\ref{BSDE_EZ}) or for meeting the Class (D) condition of the corresponding optimal density processes.
\end{remark}

\begin{theorem}\label{theorem:BSDE_RUD}
	Suppose that the market price of risk $\theta$ satisfies the following exponential integrability conditions:

(i) for $\gamma>1,\psi>1$ (i.e. Case (i)),
$$E\left[\exp\left(p\int_0^T\frac{\gamma-1}{2}|\theta_s|^2ds\right)\right]<+\infty$$
for some integer $p>1$;

(ii) for $\gamma<1, \psi>1$ (i.e. Case (ii)),
$$E\left[\exp\left(\frac{2p}{\gamma}\int_0^T\frac{1}{\gamma}|\theta_s|^2ds\right)\right]<+\infty$$
for some integer $p>2$.

Then BSDE$(0,f+g)$ with the generators $f$ and $g$  given respectively by (\ref{generator_EZ}) and (\ref{generator_EZ_2}) admits a unique solution $(Y,Z)$. Moreover, for Case (i), $Y$ is bounded from above, $e^{p
		Y^-}\in\mathcal{S}^{p}$, and $Z\in M^2$, i.e.
	\begin{equation*}
		E\left[e^{pY^{-}_{\star}}+\int_0^T|Z_s|^2ds\right]<+\infty.
	\end{equation*}
For Case (ii), $Y$ is bounded from below, $e^{\frac{2p}{\gamma}
		Y^+}\in\mathcal{S}^{p}$ and $Z\in M^2$, i.e.
	\begin{equation*}
		E\left[e^{\frac{2p}{\gamma} Y^{+}_{\star}}+\int_0^T|Z_s|^2ds\right]<+\infty;
	\end{equation*}
Herein, $Y_{\star}=\sup_{t\in[0,T]}Y_t$ is the running maximum of a
	stochastic process $Y$.
\end{theorem}

\begin{proof} We only prove the solution existence and leave its solution uniqueness in Chapter \ref{Section_convex_RDU}. Similar to the approach used in proving Theorem \ref{theorem:BSDE1}, where terminal data was truncated, we may employ an inf convolution procedure on the generator $f$ and a truncation procedure on $g$ in (\ref{BSDE_EZ}) by approximating them with Lipschitz continuous generators $f^{n,k}$. This ensures that the comparison theorem for Lipschitz BSDE applies to BSDE$(0, f^{n,k})$. 

For example, as in (\ref{inf_convolution}) in Remark \ref{remark_approx}, We define the inf convolution of \( f \) and the truncated function \( g \) as:
\[
f^{n,k}(t,y,z) = \inf_{q} \{ f(t,q) + n|q - z| \} + g(y \wedge k) \mathbf{1}_{\beta < 0} + g(y \vee (-k)) \mathbf{1}_{\beta > 0},
\]
which is Lipschitz continuous in \( (y, z) \), increasing in \( n \), and either decreasing or increasing in \( k \) depending on whether \( \beta < 0 \) or \( \beta > 0 \), respectively.
Moreover, we have
$$-\frac{\gamma-1}{2}|\theta_t|^2+g(y)\leq f^{n,k}(t,y,z)\leq f(t,z)+g(y\wedge k)\leq \frac{1}{2}|z|^2-\delta\beta$$
for Case (i) and
$$\frac{1}{2}|z|^2-\delta\beta\leq f^{n,k}(t,y,z)\leq f(t,z)+g(y)\leq \frac{1}{\gamma}|z|^2+\frac{1}{\gamma}|\theta_t|^2+g(y).$$
for Case (ii). They will provide an upper bound for $Y^{n,k}$ in Case (i) and a lower bound in Case (ii). The other side of bounds will then follow from the decreasing property of $g(y)$, which yields that $f^{n,k}(t,y,z)$ is also decreasing in $y$. Further details are provided below. 

We consider BSDE$(0,f^{n,k})$:
\begin{equation}\label{BSDE_EZ_A}
Y_t^{n,k}=\int_t^Tf^{n,k}(s,Y_s^{n,k},Z_s^{n,k})ds-\int_t^T(Z_s^{n,k})^{tr}dB_s.
\end{equation}
For \underline{Case (i)}, given that $\beta<0$ and $f^{n,k}$ has an upper bound $(\frac{1}{2}|z|^2-\delta\beta)$, we apply the comparison theorem for Lipschitz BSDE to
obtain an upper bound as
$Y_t^{n,k} \leq \overline{Y}_t$,
where $\overline{Y}$ solves
BSDE$(0,\frac{1}{2}|z|^2-\delta\beta)$, i.e.,
\begin{equation*}
	\overline{Y}_t=\int_t^T
	(\frac{1}{2}|\overline{Z}_s|^2-\delta\beta)ds-\int_t^T \overline{Z}_s^{tr}dB_s,\ t\in[0,T],
\end{equation*}
This equation obviously admits a unique solution $(\overline{Y}_t,\overline{Z}_t)=(-\delta\beta(T-t),0)$, $t\in[0,T]$.
We define $Y_{\max}:=-\delta\beta T\geq Y_t^{n,k}$.

Next, observe that $g$ is decreasing, we can therefor use $Y_{\max}$ to bound $g$ from below and consider the following BSDE to bound (\ref{BSDE_EZ_A}) from below:
\begin{align*}
	\underline{Y}_t&=\int_t^T
	(-\frac{\gamma-1}{2}|\theta_s|^2+g(Y_{\max}))ds-\int_t^T \underline{Z}_s^{tr}dB_s\\
&=g(Y_{\max})(T-t)-E\left[\int_t^T\frac{\gamma-1}{2}|\theta_s|^2ds|\mathcal{F}_t\right],\ t\in[0,T].
\end{align*}
The comparison theorem for Lipschitz BSDE then implies that $Y_t^{n,k}\geq \underline{Y}_t$. 

Using the localization argument analogous to those in the proof of Theorem \ref{theorem:BSDE1}, and defining the stopping time sequence $(\tau_j)_{j\geq 1}$:
$$\tau_j=T\wedge \{t\in[0,T]: -\underline{Y}_t >j\},$$
we deduce that there exists a solution \((Y, Z)\) to (\ref{BSDE_EZ}) that satisfies \(\underline{Y}_t \leq Y_t \leq {Y}_{\max}\) and $Z\in L^2(0,T)$.

To prove the exponential integrability of solutions, using $Y_t^{-}\leq -\underline{Y}_t$ and Jensen's inequality, we deduce that
\begin{align*}
		e^{\frac{p}{p^{\prime}}Y_t^{-}}&\leq
		C\exp\left(\frac{p}{p^{\prime}}E\left[\left.\int_t^T
		\frac{\gamma-1}{2}|\theta_s|^2ds\right|{\mathcal F}_t\right]\right)\\
		&\leq
		CE\left[\left.\exp\left(\frac{p}{p^{\prime}}(\int_t^T
		\frac{\gamma-1}{2}|\theta_s|^2ds)\right)\right|\mathcal{F}_t\right],
	\end{align*}
for $p^{\prime}>p>1$. Then, it follows from H\"older's inequality and Doob's inequality that
	\begin{align*}
		E[e^{p Y^{-}_{\star}}]&\leq CE\left[\sup_{t\in[0,T]}
		\left(E\left[\left.e^{\frac{p}{p^{\prime}}\int_t^{T}\frac{\gamma-1}{2}|\theta_s|^2dsds}\right|\mathcal{F}_t\right]\right)^{p^{\prime}}
		\right]\notag\\
		&\leq
		CE\left[\exp\left(p\int_0^T\frac{\gamma-1}{2}|\theta_s|^2ds\right)\right]<+\infty.
	\end{align*}
Moreover, the property that \( Z \in M^2 \) can be proved in a manner similar to that in Lemma \ref{lemma0}, and is thus omitted.

For \underline{Case (ii)}, recall that $\beta>0$ and $f^{n,k}$ admits a lower bound $(\frac{1}{2}|z|^2-\delta\beta)$. By applying the comparison theorem, we first obtain a lower bound as $Y^{n,k}_t\geq \underline{Y}_t$, where
$\underline{Y}$ solves BSDE$(0,\frac{1}{2}|z|^2-\delta\beta)$. It admits a unique solution $(\underline{Y}_t,\underline{Z}_t)=(-\delta\beta(T-t),0)$, $t\in[0,T]$. We define $Y_{\min}:=-\delta\beta T\leq Y_t^{n,k}$.

Since $g$ is decreasing, we can therefore use $Y_{\min}$ to bound $g$ from above and consider the following BSDE which dominates (\ref{BSDE_EZ_A}):
\begin{align*}
	\overline{Y}_t&=\int_t^T
	(\frac{1}{\gamma}|\overline{Z}_s|^2+\frac{1}{\gamma}|\theta_s|^2+g(Y_{\min}))ds-\int_t^T
	\overline{Z}_s^{tr}dB_s,\\
    &=\frac{\gamma}{2}\ln E\left[\exp\left(\int_t^T\frac{2}{\gamma}(\frac{1}{\gamma}|\theta_s|^2+g(Y_{\min}))ds\right)|\mathcal{F}_t\right],  \ t\in[0,T].
\end{align*}
Hence $Y_t^{n,k}\leq \overline{Y}_t$. 

Using the localization argument analogous to those in the proof of Theorem \ref{theorem:BSDE1}, 
and defining the stopping time sequence $(\tau_j)_{j\geq 1}$:
$$\tau_j=T\wedge \{t\in[0,T]: \overline{Y}_t >j\},$$
we deduce that there exists a solution \((Y, Z)\) to (\ref{BSDE_EZ}) that satisfies \({Y}_{\min} \leq Y_t \leq \overline{Y}_t\). Consequently, similar to Case (i), we deduce the exponential integrability of $Y^{+}$:
$$E[e^{\frac{2p}{\gamma}Y_{\star}^+}]\leq C E\left[\exp\left(\frac{2p}{\gamma}\int_0^T\frac{1}{\gamma}|\theta_s|^2ds\right)\right]<+\infty.$$
\end{proof}

\begin{remark}
{The exponential integrability condition on the market price of risk $\theta$ can be readily verified. We provide an example using the Heston stochastic volatility model. 
Consider a market with a single stock whose coefficients
depend on a single stochastic factor driven by a two-dimensional
Brownian motion, namely, $m=2$, $d=1$, and
\[
dS_{t} = \lambda V_t S_{t} \, dt + \sqrt{V_{t}} S_{t} \, dB_{1,t},
\]
\[
dV_{t} = \left(a - b V_{t}\right) \, dt + \sigma \sqrt{V_t} (\kappa_1 \, dB_{1,t} + \kappa_2 \, dB_{2,t}),
\]
where \( \lambda \), \( a \), \( b \), and \( \sigma \) are constant coefficients. We assume that \( \sigma > 0 \) and the constants \( \kappa_1 \) and \( \kappa_2 \) satisfy \( |\kappa_1|^2 + |\kappa_2|^2 = 1 \). The market price of risk is then given by \( \theta_t = \lambda \sqrt{V_t} \). Let \( \mathcal{C} = \mathbb{R} \). To verify the exponential integrability condition on \( \theta\), we apply Theorem 5.1 from Zeng and Taksar \cite{ZengTaksar2013}:}

\underline{Case (i):} The exponential integrability condition holds if and only if there is either inequality
\[
p \leq \frac{b^2}{(\gamma - 1) \lambda^2 \sigma^2}
\]
or both inequalities
\[
p > \frac{b^2}{(\gamma - 1) \lambda^2 \sigma^2}\quad  \text{ and } \quad T < \frac{1}{m_2} \text{arccot}\left(-\frac{m_1}{m_2}\right),
\]
where
$$ m_1 = \frac{1}{2}b, \quad  m_2 = \frac{1}{2} \sqrt{-b^2 + p(\gamma - 1) \lambda^2 \sigma^2}, $$
and \(\text{arccot}\) denotes the inverse cotangent function. Hence, if $b^2>(\gamma - 1) \lambda^2 \sigma^2$, there is $p>1$ such that the exponential integrability condition holds. 

\underline{Case (ii):} The exponential integrability condition holds if and only if there is either inequality
\[
p \leq \frac{\gamma^2 b^2}{4 \lambda^2 \sigma^2}
\]
or both inequalities 
\[
p > \frac{\gamma^2 b^2}{4 \lambda^2 \sigma^2} \quad \text{and} \quad T < \frac{1}{m_2} \text{arccot} \left(-\frac{m_1}{m_2}\right),
\]
where
$$ m_1 = \frac{1}{2}b \quad \text{and} \quad  m_2 = \frac{1}{2} \sqrt{-b^2 + \frac{4p}{\gamma^2} \lambda^2 \sigma^2}. $$
Hence, if $\gamma^2 b^2>8 \lambda^2 \sigma^2$, there is $p>2$ such that the exponential integrability condition holds. 
\end{remark}

\section{Epstein-Zin recursive utility maximization}

To solve the recursive utility maximization problem, we use the \emph{martingale optimality principle}, i.e. we aim to construct a conditional value process $V(t,\xi)$ such that
$$V(t,X_t^{u,c})+\int_0^{t}f^{EZ}(X_s^{u,c}c_s,V(s,X_s^{u,c}))ds,\ \ \ t\in[0,T],$$
is a supermartingale for any $(u,c)\in\mathcal{A}_D^{EZ}$ and a martingale for some $(u^{\star},c^{\star})$. Moreover, $$V(T,X_T^{u,c})=U(X_T^{u,c})=\frac{(X_T^{u,c})^{1-\gamma}}{1-\gamma}.$$ Then, for any $(u,c)\in\mathcal{A}_D^{EZ}$,
\begin{equation}\label{EZ_BSDE_comparison}
V(t,X_t^{u,c})=U(X_T^{u,c})+\int_t^TdK_s+\int_t^Tf^{EZ}(X_s^{u,c}c_s,V(s,X_s^{u,c}))ds-\int_t^T(Z_s^{u,c})^{tr}dB_s,
\end{equation}
for some increasing process $K$ with $K_0=0$ and $\mathbb{R}^m$-valued predictable density process $Z^{u,c}$. Moreover, with $(u^{\star},c^{\star})$,
$$V(t,X_t^{u^{\star},c^{\star}})=U(X_T^{u^{\star},c^{\star}})+\int_t^Tf^{EZ}(X_s^{u^{\star},c^{\star}}c_s^{\star},V(s,X_s^{u^{\star},c^{\star}}))ds-\int_t^T(Z_s^{u^{\star},c^{\star}})^{tr}dB_s,$$
for some $\mathbb{R}^m$-valued predictable density process $Z^{u^{\star},c^{\star}}$.

By comparing (\ref{EZ_BSDE}) with (\ref{EZ_BSDE_comparison}) and applying the comparison theorem in Lemma \ref{lemma_EZ}, we conclude that \( V(t, X_t^{u,c}) \geq V_t^{u,c} \) for any \((u,c) \in \mathcal{A}_D^{EZ}\). Furthermore, the uniqueness of the solution to BSDE (\ref{EZ_BSDE}) ensures that \( V(t, X_t^{u^{\star},c^{\star}}) = V_t^{u^{\star},c^{\star}} \).

We construct the conditional value process as
\begin{equation}
V(t,X_t^{u,c})=U(X_t^{u,c})e^{Y_t},
\end{equation}
with $Y$ solving BSDE (\ref{BSDE_EZ}). To verify the supermartingale and martingale properties, an application of It\^o's
formula to $U(X_t^{u,c})e^{Y_t}$, for $(u,c)\in\mathcal{A}_D^{EZ}$, gives
\begin{align}\label{value_process_EZ_2}	
&V(t,X_t^{\pi,C})+\int_0^tf^{EZ}(X_s^{u,c}c_s,V(s,X_s^{u,c}))ds\\
=&\ U(x)e^{Y_0}+\int_0^t (X_s^{u,c})^{1-\gamma}e^{Y_s}\left[A_s^{u,c}ds+(u_s^{tr}\sigma_s+\frac{Z_s^{tr}}{1-\gamma})dB_s\right],\notag
\end{align}
	where
\begin{align*}	
A_s^{u,c}=&\ -\frac12\gamma|\sigma_s^{tr}\pi_s|^2+\pi_s^{tr}\sigma_s(Z_s+\theta_s)
+\frac{|Z_s|^2}{2(1-\gamma)}\\
&+\frac{f^{EZ}(X_s^{u,c}c_s,U(X_s^{u,c})e^{Y_s})}{(X_s^{u,c})^{1-\gamma}e^{Y_s}}-c_s-
\frac{f(s,Z_s)+g(Y_s)}{1-\gamma}.
\end{align*}

We require that $A_s^{u,c}\leq 0$ for any admissible $(u,c)\in\mathcal{A}_D^{EZ}$ and $A_s^{u^{\star},c^{\star}}=0$ for some $(u^{\star},c^{\star})$. Hence, the generator must take the form
\begin{align*}	
f(s,Z_s)+g(Y_s)=& -\frac{1-\gamma}{2\gamma}\min_{u_s\in\mathcal{C}}
\left|\gamma\sigma_s^{tr}u_s-(Z_s+\theta_s)\right|^2 +\frac{1-\gamma}{2\gamma}|Z_s+\theta_s|^2+ \frac{1}{2}|Z_s|^2\\
&+(1-\gamma)\sup_{c_s\in\mathbb{R}_+}\left\{\frac{f^{EZ}(X_s^{u,c}c_s,U(X_s^{u,c})e^{Y_s})}{(X_s^{u,c})^{1-\gamma}e^{Y_s}}-c_s\right\}.
\end{align*}
Note that the second term can be rewritten as
\begin{align*}
&\delta e^{-\frac{Y_s}{\beta}}\max_{c_s\in\mathbb{R}_+}
\left\{\frac{\beta}{1-\gamma}c_s^{\frac{1-\gamma}{\beta}}-c_s\underbrace{\frac{1}{\delta}e^{\frac{Y_s}{\beta}}}_{P_s}\right\}-\frac{\delta\beta}{1-\gamma}\\
&\quad =\delta e^{-\frac{Y_s}{\beta}}\frac{\beta-1+\gamma}{1-\gamma}P_s^{\frac{\gamma-1}{\beta-1+\gamma}}-\frac{\delta\beta}{1-\gamma}\\
&\quad= \frac{1}{1-\gamma}\left(\delta^{\psi}\frac{\beta}{\psi}e^{-\frac{\psi}{\beta}Y_s}\right)-\frac{\delta\beta}{1-\gamma},
\end{align*}
where we utilized the convex duality between $\frac{1}{\alpha} x^{\alpha}$ and $\frac{1-\alpha}{\alpha}y^{\frac{\alpha}{\alpha-1}}$ with $\alpha=\frac{1-\gamma}{\beta}=1-1/\psi<1$ in the first equality. This lead to the generator $f$ and $g$ as defined in (\ref{generator_EZ}) and (\ref{generator_EZ_2}), respectively. Moreover, the optimal investment and consumption $(u^{\star},c^{\star})$ are given by
\begin{align}\label{pi_EZ}
\gamma\sigma_t^{tr}u_t^{\star}&=\mbox{Proj}_{\gamma\sigma_t^{tr}\mathcal{C}} \left(Z_t+\theta_t\right);\\
c_t^{\star}&=\delta^{\psi}e^{-\frac{\psi}{\beta}Y_t},\ \hbox{$\mathbb{P}$-a.s.,}\ \text{for}\ t \in	[0,T].\notag
\end{align}

With $(u^{\star},c^{\star})$ given in (\ref{pi_EZ}),
(\ref{value_process_EZ_2}) then yields that
\begin{align*}
&d(X_t^{u^{\star},c^{\star}})^{1-\gamma}e^{Y_t}\\
&=\
(X_t^{u^{\star},c^{\star}})^{1-\gamma}e^{Y_t}\left[((1-\gamma)\sigma_t^{tr}u_t^{\star}+Z_t)^{tr}dB_t
-\frac{f^{EZ}(X_t^{u^{\star},c^{\star}}c_t^{\star},U(X_t^{u^{\star},c^{\star}})e^{Y_t})}
{(X_t^{u{^{\star},c^{\star}}})^{1-\gamma}e^{Y_t}}dt\right],
\end{align*}	
from which we deduce that
\begin{equation}\label{optimal_state_RDU}
(X_t^{u^{\star},c^{\star}})^{1-\gamma}e^{Y_t}=x^{1-\gamma}e^{Y_0}\exp\left(\int_0^t\left(\delta\beta-\beta\delta^{\psi}e^{-\frac{\psi}{\beta}Y_s}\right)ds\right)L_t^{u{^{\star}}}
\end{equation}
where
\begin{equation}\label{optimal_density_RDU}
L_t^{u^{\star}}=\mathcal{E}_t(\int_0^{\cdot}((1-\gamma)\sigma_t^{tr}u_t^{\star}+Z_t)^{tr}dB_t).
\end{equation}

\begin{theorem} \label{thm:optimal_EZ}
Suppose that the parameters \( (\gamma, \psi) \) satisfy either Case (i) or Case (ii) and that the market price of risk \( \theta \) meets the Novikov's condition (\ref{Novikov}), the exponential integrability conditions in terms of $\gamma$ outlined in Theorem \ref{theorem:BSDE_RUD}. Additionally, \( \theta \) is required to satisfy the following exponential integrability condition in terms of \( \psi \):
\begin{equation}\label{exponential_psi}
E \left[\exp\left( \int_0^{T} \frac{1}{2} \frac{p}{p - 1} [p(\psi - 1)^2 + (\psi - 1)] |\theta_s|^2 ds \right)\right]<+\infty
\end{equation}
for some integer $p>1$.

Let $(Y,Z)$ be the unique solution
	to BSDE$(0,f+g)$ in (\ref{BSDE_EZ}) with the generators $f$ and $g$ given by (\ref{generator_EZ}) and (\ref{generator_EZ_2}), respectively. Then, the value function of
	the optimization problem (\ref{EZ_formulation}) with admissible set
	$\mathcal{A}_D^{EZ}$ is given by
	\begin{equation}\label{valuefunction_EZ}
		V(0, x)= U(x)e^{Y_0}=\frac{x^{1-\gamma}}{1-\gamma}e^{Y_0},
	\end{equation}
	and there exists an optimal trading strategy $u^{\star}$ and an optimal consumption process $c^{\star}$, both in the admissible set $\mathcal{A}_D^{EZ}$, given by (\ref{pi_EZ}).
\end{theorem}

\begin{proof} We have verified the supermartingale property and local martingale property of the conditional value process.  It remains to show that, with $(u^{\star}, c^{\star})$ as defined in (\ref{pi_EZ}), the process
$$V(\tau, X_{\tau}^{u^{\star},c^{\star}})=\frac{(X_{\tau}^{u^{\star}, c^{\star}})^{1-\gamma}}{1-\gamma}e^{Y_{\tau}}, \quad \tau\in [0,T]$$ is in Class (D), and the integrability condition (\ref{integrability_condition1}) or (\ref{integrability_condition2}) holds.  

Note that for Case (i), $\beta<0$ but $Y\leq Y_{\max}$, so we have 
\begin{equation}\label{bound_1}
\exp\left(\int_0^t\left(\delta\beta-\beta\delta^{\psi}e^{-\frac{\psi}{\beta}Y_s}\right)ds\right)\leq e^{(\delta\beta-\beta\delta^{\psi}e^{-\frac{\psi}{\beta}Y_{\max}})t}.
\end{equation}
For Case (ii), $\beta>0$, so we have
\begin{equation}\label{bound_2}
\exp\left(\int_0^t\left(\delta\beta-\beta\delta^{\psi}e^{-\frac{\psi}{\beta}Y_s}\right)ds\right)\leq e^{\delta\beta t}.
\end{equation}
Thus, it is sufficient to verify that the optimal density process $L^{u^{\star}}$ as defined in (\ref{optimal_density_RDU}) is in Class (D). According to Lemma \ref{Lemma_RDU_driver} (ii),
$$(1-\gamma)\sigma_t^{tr}u_t^{\star}+Z_t\in\partial_{z}f(t,Z_t).$$
Therefore, we aim to prove that $L_T^{q^{\star}}$ has a finite entropy for any $q^{\star}_t\in\partial_z f(t,Z_t)$, where
\begin{equation*}		L_{t}^{q^{\star}}=\mathcal{E}_t\left(\int_0^{\cdot}q_u^{\star}dB_u\right),\quad t\in[0,T].
	\end{equation*}
which will be established in the Lemma \ref{lemma_class_D_RDU} after the proof.

Next, we verify the integrability conditions (\ref{integrability_condition1}) or (\ref{integrability_condition2}). First, observe that by comparing the definition of \(g(y)\) in (\ref{generator_EZ_2}) with the expression of \(c^{\star}\) in (\ref{pi_EZ}), we have:
\[
c_t^{\star} = \frac{\psi}{\beta}(g(Y_t) + \delta \beta),
\]
which is always bounded in both Case (i) and Case (ii) (see Remark \ref{remark_difficulty}). Since \(C_t = X_t^{u,c} c_t\), it follows that conditions (\ref{integrability_condition1}) and (\ref{integrability_condition2}) are equivalent to the following:

\underline{Case (i)}: For \(\gamma > 1\) and \(\psi > 1\),
\begin{equation*}
E\left[(X_T^{u^{\star},c^{\star}})^{1-\gamma}\right] < \infty, \quad E\left[\int_0^T (X_t^{u^{\star},c^{\star}})^{1-\frac{1}{\psi}} dt \right] < \infty;
\end{equation*}

\underline{Case (ii)}: For \(\gamma < 1\) and \(\psi > 1\),
\begin{equation*}
E\left[(X_T^{u^{\star},c^{\star}})^{1-\gamma \wedge \frac{1}{\psi}}\right] < \infty, \quad E\left[\int_0^T (X_t^{u^{\star},c^{\star}})^{1-\gamma \wedge \frac{1}{\psi}} dt \right] < \infty,
\end{equation*}
respectively.

To verify the inequalities involving the exponent of \( (1-\gamma) \), by using (\ref{optimal_state_RDU}), we have 
\begin{align*}
(X_t^{u^{\star},c^{\star}})^{1-\gamma} &= x^{1-\gamma} e^{Y_0} e^{-Y_t} \exp\left( \int_0^t \left( \delta\beta - \beta \delta^{\psi} e^{-\frac{\psi}{\beta}Y_s} \right) ds \right) L_t^{u^{\star}}. 
\end{align*}
For Case (i), since \( Y_T = 0 \), we obtain \( E\left[ (X_T^{u^{\star},c^{\star}})^{1-\gamma} \right] < \infty \) by using (\ref{bound_1}) and the finite entropy condition of \( L^{u^{\star}} \). For Case (ii) with $\gamma < \frac{1}{\psi}$,
since \( Y_t \geq Y_{min} \), we use (\ref{bound_2}) to further obtain
\[
(X_t^{u^{\star},c^{\star}})^{1-\gamma} \leq x^{1-\gamma} e^{Y_0} e^{-Y_{min}} e^{\delta\beta t} L_t^{u^{\star}},
\]
so we verify the integrability condition for Case (ii) when \( \gamma < \frac{1}{\psi} \).

To verify the other inequalities involving the exponent of \( 1 - \frac{1}{\psi} \in (0,1) \), as in (\ref{EMM}), we define an equivalent minimal local
martingale measure (MLMM) \( \mathbb{Q}^{\theta} \) on \( \mathcal{F}_T \) by
\begin{equation*}
    \frac{d\mathbb{Q}^{\theta}}{d\mathbb{P}} := L_T^{\theta} = \mathcal{E}_T\left( -\int_0^{\cdot} \theta_t^{\text{tr}} dB_t \right),
\end{equation*}
where \( \mathcal{E}(\cdot) \) denotes the stochastic exponential. Note that this is permissible since \( \theta \) satisfies the Novikov's condition (\ref{Novikov}). Then, under \( \mathbb{Q}^{\theta} \), the wealth process \( X^{\pi^{\star},\mu^{\star}} \) becomes a supermartingale following
\begin{equation*}
    dX_t^{\mu^{\star},c^{\star}} = X_t^{\mu^{\star},c^{\star}} \left[ (u_t^{\star})^{\text{tr}} \sigma_t dB_t^{\theta} - c_t^{\star} \right],
\end{equation*}
where \( B_t^{\theta} := B_t + \int_0^{t} \theta_u du \), \( t \in [0,T] \), is an \( m \)-dimensional Brownian motion under the MLMM \( \mathbb{Q}^{\theta} \). Hence, by using H\"older's inequality and the supermartingale property of $X^{\mu^{\star},c^{\star}}$ under $\mathbb{Q}^{\theta}$, we obtain
\begin{align*}
    E\left[ (X_t^{u^{\star},c^{\star}})^{1 - \frac{1}{\psi}} \right]
    &= E^{\mathbb{Q}^{\theta}} \left[ \mathcal{E}_t\left( \int_0^{\cdot} \theta_t^{\text{tr}} dB_t^{\theta} \right) (X_t^{u^{\star},c^{\star}})^{1 - \frac{1}{\psi}} \right] \\
    &\leq E^{\mathbb{Q}^{\theta}} \left[ \mathcal{E}_t\left( \int_0^{\cdot} \theta_t^{\text{tr}} dB_t^{\theta} \right)^{\psi} \right]^{\frac{1}{\psi}} 
    E^{\mathbb{Q}^{\theta}} \left[  X_t^{u^{\star},c^{\star}} \right]^{1 - \frac{1}{\psi}} \\
    &\leq {E^{\mathbb{Q}^{\theta}} \left[ \mathcal{E}_t\left( \int_0^{\cdot} \theta_t^{\text{tr}} dB_t^{\theta} \right)^{\psi} \right]^{\frac{1}{\psi}}} x^{1-\frac{1}{\psi}},
\end{align*}
for any \( t \in [0,T] \). 
We conclude the proof by deriving an upper bound for 
$$
E^{\mathbb{Q}^{\theta}} \left[ \mathcal{E}_t\left( \int_0^{\cdot} \theta_t^{\text{tr}} dB_t^{\theta} \right)^{\psi} \right],
$$
which is independent of \(t\). Indeed, by applying H\"older's inequality, we obtain
\begin{align*}
E^{\mathbb{Q}^{\theta}} \left[ \mathcal{E}_t\left( \int_0^{\cdot} \theta_t^{\text{tr}} dB_t^{\theta} \right)^{\psi} \right] 
&= E \left[\mathcal{E}_t\left( -\int_0^{\cdot} \theta_t^{\text{tr}} dB_t \right) \mathcal{E}_t\left( \int_0^{\cdot} \theta_t^{\text{tr}} dB_t^{\theta} \right)^{\psi} \right] \\
&= E \left[\exp\left( \int_0^{t} (\psi - 1) \theta_s^{\text{tr}} dB_s \right) 
\exp\left( \int_0^{t} \frac{1}{2} (\psi - 1) |\theta_s|^2 ds \right)\right] \\
&= E \left[\mathcal{E}_t\left( \int_0^{\cdot} p(\psi - 1) \theta_t^{\text{tr}} dB_t^{\theta} \right)^{\frac{1}{p}}\right.\\
&\ \ \cdot\left.\exp\left( \int_0^{t} \frac{1}{2} [p(\psi - 1)^2 + (\psi - 1)] |\theta_s|^2 ds \right)\right] \\
&\leq E \left[\mathcal{E}_t\left( \int_0^{\cdot} p(\psi - 1) \theta_t^{\text{tr}} dB_t^{\theta} \right)\right]^{\frac{1}{p}}\\
&\ \  \cdot E \left[\exp\left( \int_0^{t} \frac{1}{2} \frac{p}{p - 1} [p(\psi - 1)^2 + (\psi - 1)] |\theta_s|^2 ds \right)\right]^{1 - \frac{1}{p}}.
\end{align*}
The first expectation is bounded by \(1\). By the exponential integrability condition \((\ref{exponential_psi})\), the second expectation is also bounded, from which we conclude.
\end{proof}

\begin{lemma}\label{lemma_class_D_RDU}
	The optimal density process $L_T^{q^{\star}}$ has a finite entropy.
 Hence, by De la Vall\'ee-Poussin theorem, $L^{q^{\star}}$ is in
	Class (D) and is therefore a uniformly integrable martingale.
 \end{lemma}

\begin{proof} The proof is similar to Lemmas~ \ref{lemma} and~\ref{lemma_consumption}, so we only highlight the main differences. With a localization sequence $(\sigma_j)_{j\geq 1}$, we rewrite BSDE$(0,f+g)$ under $\mathbb{Q}^{q^{\star}}$ as
$$dY_t=-\left[f(t,Z_t)-Z_t^{tr}q_t^{\star}+g(Y_t)\right]dt+Z_t^{tr}dB_t^{q^{\star}}.$$

For \underline{Case (i)}, we have
\begin{align*}
E\left[L_{\sigma_{j}}^{q^{\star}}
		Y_{\sigma_j}\right] &= E^{\mathbb{Q}^{q^{\star}}}[
		Y_{\sigma_{j}}]\\		
&=E^{\mathbb{Q}^{q^{\star}}}\left[
		Y_0+\int_0^{\sigma_{j}}
(-f(u,Z_u)+Z_u^{tr}q_u^{\star}-g(Y_u))du\right]\\
&=E^{\mathbb{Q}^{q^{\star}}}\left[
		Y_0+\int_0^{\sigma_{j}}
(f^{\star}(u,q_u^{\star})-g(Y_u))du\right],
	\end{align*}
where $f^{\star}$ is the convex dual of $f$:
\begin{equation}\label{dual_RDU}
	{f}^{\star}(t,q):=\sup_{z\in\mathbb{R}^m}\left({z}^{tr}q-f(t,z)\right), \quad (t,q)\in [0,T]\times \mathbb{R}^m, 
\end{equation}
with a lower bound of $f^{\star}(t,q)\geq \frac{1}{2}|q|^2$ due to the upper bound of $f$ derived in Lemma \ref{Lemma_RDU_driver}.
Then, utilizing $Y\leq Y_{\max}$ from Theorem \ref{theorem:BSDE_RUD} and $g(Y)\leq -\delta\beta$ from Lemma \ref{Lemma_RDU_driver}, we deduce that
$$Y_{\max}\geq Y_0+E^{\mathbb{Q}^{q^{\star}}}\left[\int_0^{\sigma^j}\frac{1}{2}|q_u^{\star}|^2du\right]+\delta\beta T,$$
and consequently
\begin{equation}\label{entropy_inequ_1}	
E\left[L_{\sigma_j}^{q^{\star}}\ln
L_{\sigma_j}^{q^{\star}}\right]=E^{\mathbb{Q}^{q^{\star}}}\left[\int_0^{\sigma_{j}}\frac12|q_u^{\star}|^2du\right]\leq Y_{\max}-Y_0-\delta\beta T.
\end{equation}

For \underline{Case (ii)}, we have
	\begin{align*}
E\left[L_{\sigma_{j}}^{q^{\star}}\frac{2}{\gamma}
		Y_{\sigma_j}\right] &= E^{\mathbb{Q}^{q^{\star}}}[\frac{2}{\gamma}
		Y_{\sigma_{j}}]\\		
&=E^{\mathbb{Q}^{q^{\star}}}\left[\frac{2}{\gamma}Y_0
		+\frac{2}{\gamma}\int_0^{\sigma_{j}}
(-f(u,Z_u)+Z_u^{tr}q_u^{\star}-g(Y_u))du\right]\notag\\
&=E^{\mathbb{Q}^{q^{\star}}}\left[\frac{2}{\gamma}Y_0
		+\frac{2}{\gamma}\int_0^{\sigma_{j}}
(f^{\star}(u,q_u^{\star})-g(Y_u))du\right].
	\end{align*}
Considering the upper bound of $f$ in Lemma \ref{lemma_EZ}, its convex dual $f^{\star}$ is bounded below by
$$f^{\star}(t,q)\geq \frac{\gamma}{4}|q|^2-\frac{1}{\gamma}|\theta_t|^2.$$
Also, from Theorem \ref{theorem:BSDE_RUD} and the decreasing property of $g$, we have $g(Y)\leq g(Y_{\min})$. Hence, we have
\begin{align}\label{inequality_RDU_0}
E\left[L_{\sigma_{j}}^{q^{\star}}\frac{2}{\gamma}
		Y_{\sigma_j}\right]\geq &\ \frac{2}{\gamma}Y_0+E^{\mathbb{Q}^{q^{\star}}}\left[\int_0^{\sigma_{j}}\frac12|q_u^{\star}|^2du\right]\notag \\[2mm]
&-E^{\mathbb{Q}^{q^{\star}}}\left[\int_0^{\sigma_{j}}\frac{2}{\gamma^2}|\theta_u|^2du\right]-\frac{2}{\gamma}g(Y_{\min})T.
\end{align}
Furthermore, using the Fenchel inequality, we have
	\begin{equation}\label{inequality_RDU_1}
		E\left[L_{\sigma_{j}}^{q^{\star}}\frac{2}{\gamma}Y_{\sigma_j}\right]\leq \ 
		\frac{1}{p}{E\left[L_{\sigma_j}^{q^{\star}}\ln L_{\sigma_j}^{q^{\star}}\right]}-
		\frac{1}{p} {E\left[L_{\sigma_j}^{q^{\star}}\right]\ln p}+E\left[e^{
			\frac{2p}{\gamma}Y^+_{\star}}\right]
	\end{equation}
and
\begin{align}\label{inequality_RDU_2}
	E^{\mathbb{Q}^{q^{\star}}}\left[\int_0^{\sigma_{j}}\frac{2}{\gamma^2}|\theta_u|^2du\right]\leq &\ 
		\frac{1}{p}{E\left [L_{\sigma_j}^{q^{\star}}\ln L_{\sigma_j}^{q^{\star}}\right]}-
		\frac{1}{p} {E\left[L_{\sigma_j}^{q^{\star}}\right]\ln p}\notag \\[2mm]
		&\displaystyle +E\left[e^{\frac{2p}{\gamma}\int_0^T\frac{1}{\gamma}|\theta_u|^2du}\right].
	\end{align}
Substituting (\ref{inequality_RDU_1}) and (\ref{inequality_RDU_2}) into (\ref{inequality_RDU_0}) and recalling $p>2$, we see
\begin{align}\label{entropy_inequ_2}
(1-\frac{2}{p})E\left[L_{\sigma_j}^{q^{\star}}\ln L_{\sigma_j}^{q^{\star}}\right]\ \leq \ &-\frac{2}{p}\ln p-\frac{2}{\gamma}Y_0+\frac{2}{\gamma}g(Y_{\min})T\notag\\[2mm]
&\displaystyle +E\left[e^{\frac{2p}{\gamma}Y^+_{\star}}\right]+E\left[e^{\frac{2p}{\gamma}\int_0^T\frac{1}{\gamma}|\theta_u|^2du}\right].
\end{align}

We conclude by sending $j\rightarrow\infty$ in the inequalities (\ref{entropy_inequ_1}) and (\ref{entropy_inequ_2}), and
	using Fatou's lemma.
\end{proof}

\section[A convex dual representation]{A convex dual representation of the Epstein-Zin utility maximization model}\label{Section_convex_RDU}

To conclude, we provide a convex dual representation of
the solution component $Y$ of BSDE$(0,f+g)$ in (\ref{BSDE_EZ}) with the generators $f$ and $g$ given by (\ref{generator_EZ}) and (\ref{generator_EZ_2}), respectively. This will complete the
proof of Theorem \ref{theorem:BSDE_RUD} for its solution uniqueness.

We first introduce the admissible set of the convex dual problem. For
an $\mathbb{R}^m$-valued $\mathbb{F}$-predictable process $q$,
we
define the stochastic exponential 
$$
L^{q}:=\mathcal{E}\left(\int_0^{\cdot}q_u^{tr}dB_u\right).
$$
If $L_T^{q}$ has a finite entropy,  i.e. $E[L_T^{q}\ln L_T^{q}]<+\infty$, then De la
Vall\'ee-Poussin theorem implies that $L^{q}$ is in Class (D) and
therefore a uniformly integrable martingale. We then define the
probability measure $\mathbb{Q}^{q}$ on $\mathcal{F}_T$ by
$d\mathbb{Q}^{q} :=L_T^{q} d\mathbb{P}$, and introduce the
admissible set
\begin{align*}
	\mathcal{A}_D^{EZ,\star}=\ \Biggl\{&q\in L^2(0,T):\
	L_T^{q}\ \mbox{has a finite
		entropy such} \\
	&\qquad \text{that}\ E^{\mathbb{Q}^q}\left[\int_0^T|f^{\star}(s,q_s)|ds\right]<+\infty\ \mbox{with}\ d\mathbb{Q}^{q} :=L_T^{q}d\mathbb{P}\ \Biggr\}
	\text{.}
\end{align*}%

\begin{theorem}\label{theorem_dual_RDU}
	Suppose that the parameters $(\gamma,\psi)$ satisfies either Case (i) or Case (ii), and that the market price of risk $\theta$ satisfies the exponential integrability conditions in Theorem \ref{theorem:BSDE_RUD}. Then, the solution component 
	$Y$ to BSDE$(0,f+g)$ admits the following
	convex dual representation
	\begin{equation}\label{dual_formula_RDU}
	Y_t=\esssup_{q\in\mathcal{A}_D^{EZ,\star}}Y_t^{q},
	\end{equation}
where $Y^q$ is the unique solution to BSDE$(0,-f^{\star}+g)$:
\begin{equation}\label{dual_BSDE_RDU}
Y_t^q=\int_t^T\left[-f^{\star}(s,q_s)+g(Y_s^q)\right]ds-\int_t^T(Z_s^{q})^{tr}dB_s^q.
\end{equation}
Moreover, there exists an
	optimal density process $q^{\star}\in\mathcal{A}_D^{EZ,\star}$
	such that $Y=Y^{q^{\star}}$.
\end{theorem}

\begin{proof}
The existence of BSDE (\ref{dual_BSDE_RDU}) can be proved in a manner similar to that for Lemma \ref{lemma_EZ} by applying the truncation and localization arguments, while the uniqueness follows from the decreasing property of the generator \( g(y) \) in \( y \).
On the other hand, BSDE$(0,f+g)$ reads
\begin{equation}\label{dual_BSDE_RDU_2}
Y_t=\int_t^T\left[f(s,Z_s)-Z_s^{tr}q_s+g(Y_s)\right]ds-\int_t^TZ_s^{tr}dB_s^{q},
\end{equation}
for any \( q \in \mathcal{A}_D^{EZ,\star} \), where \( B^q_t := B_t - \int_0^{t} q_u du \), \( t \in [0,T] \), is an \( m \)-dimensional Brownian motion under \( \mathbb{Q}^q \) defined at the beginning of this subsection.

With a localization sequence $(\sigma_{j})_{j\geq 1}$,
combining (\ref{dual_BSDE_RDU}) and (\ref{dual_BSDE_RDU_2}) yields
{
\begin{align*}
	&Y_{t}-Y_t^{q}\\
=&\ Y_{\sigma_{j}}-Y_{\sigma_{j}}^{q}+
\int_t^{\sigma_j}\left[f(s,Z_{s})-Z_s^{tr}q_s+{f}^{\star}(s,q_s)+g(Y_s)-g(Y_s^q)\right]ds\\
&\ -\int_t^{\sigma_j}(Z_s-Z_s^q)^{tr}dB_s^q\\
=&\ Y_{\sigma_{j}}-Y_{\sigma_{j}}^{q}+
\int_t^{\sigma_j}\left[f(s,Z_{s})-Z_s^{tr}q_s+{f}^{\star}(s,q_s)+h_s(Y_s-Y_s^q)\right]ds\\
&\ -\int_t^{\sigma_j}(Z_s-Z_s^q)^{tr}dB_s^q,
\end{align*}
}
where $$h_s=\frac{g(Y_s)-g(Y_s^q)}{Y_s-Y_s^{q}}1_{\{Y_s-Y_s^{q}\neq 0\}}\leq 0$$
by the decreasing property of $g$.
It follows that
\begin{align*}
	&Y_{t}-Y_t^{q}\\
=&\ e^{\int_t^{\sigma_j}h_udu}(Y_{\sigma_{j}}-Y_{\sigma_{j}}^{q})
+\int_t^{\sigma_j}e^{\int_t^{s}h_udu}\left[f(s,Z_{s})-Z_s^{tr}q_s+{f}^{\star}(s,q_s)\right]ds\\
&\ -\int_t^{\sigma_j}e^{\int_t^{s}h(u)du}(Z_s-Z_s^q)^{tr}dB_s^{q}\\[2mm]
=&\ E^{\mathbb{Q}^q}\biggl[e^{\int_t^{\sigma_j}h_udu}(Y_{\sigma_{j}}-Y_{\sigma_{j}}^{q})\\
&\ +\int_t^{\sigma_j}e^{\int_t^{s} h_udu}\left[f(s,Z_{s})-Z_s^{tr}q_s+{f}^{\star}(s,q_s)\right]ds\Bigm |\mathcal{F}_t\biggr].
\end{align*}

By the Fenchel-Moreau theorem, we then deduce that
\begin{align*}
	f(s,Z_{s})-Z_s^{tr}q_s+{f}^{\star}(s,q_s)&\geq 0,\quad \forall q\in\mathcal{A}_D^{EZ,\star}
\end{align*}
and the equality holds
for $q^{\star}_s\in\partial_z f(s,Z_s)$.
Hence,
$$Y_t-Y_t^{q}\geq E^{\mathbb{Q}^{q}}\left[e^{\int_t^{\sigma_j} h_udu}(Y_{\sigma_j}-Y_{\sigma_j}^{q})\Bigm |\mathcal{F}_t\right], \quad \forall q\in\mathcal{A}_D^{EZ,\star}
$$ 
 and $$Y_t-Y_t^{q^{\star}}= E^{\mathbb{Q}^{q^{\star}}}\left [e^{\int_t^{\sigma_j} h_udu}(Y_{\sigma_j}-Y_{\sigma_j}^{q^{\star}})\Bigm |\mathcal{F}_t\right].$$
For Case (i), $Y$ is bounded from above, and moreover,
$$E^{\mathbb{Q}^q}\left [Y^-_{\star}\right]\leq \frac{1}{p} E^{\mathbb{Q}^q}[\ln L^q_T]-\frac{1}{p}\ln p+E\left[e^{pY^-_{\star}}\right]<\infty$$
for $p>1$.
For Case (ii), $Y$ is bounded from below, and moreover,
$$E^{\mathbb{Q}^q}\left[Y^+_{\star}\right]\leq \frac{E^{\mathbb{Q}^q}[\ln L^q_T]}{2p/\gamma}-\frac{\ln (2p/\gamma)}{2p/\gamma}+E\left[e^{\frac{2p}{\gamma}Y^+_{\star}}\right]<\infty$$
for $p>2$. Furthermore,
note that $e^{\int_t^{\sigma_j}h_udu}\leq 1$, and $Y^{q}$ is uniformly integrable under $\mathbb{Q}^q$. Sending $j\rightarrow \infty$, we obtain that $Y_t\geq Y_t^{q}$ for any $q\in\mathcal{A}_D^{EZ,\star}$.

On the other hand, $Y_t=Y_t^{q^{\star}}$. To achieve this, it
remains to prove $q^{\star}\in\mathcal{A}_D^{EZ,\star}$.
Indeed,  the finite entropy of $L_{T}^{q^{\star}}$ has already been established in Lemma \ref{lemma_class_D_RDU}. Moreover, we can prove  $q^{\star}\in L^{2}[0,T]$ and
$$E^{\mathbb{Q}^{q^{\star}}}\left[\int_0^T|f^{\star}(s,q^{\star}_s)|ds\right]<+\infty$$
with a similar approach to the proof in Lemma \ref{lemma000}, and as such, the details are omitted.
\end{proof}

%
%
%


\chapter{Conclusions}

We have developed a systematic approach to solving utility maximization problems for an investor in constrained and unbounded financial markets. Our main methodology involves the theory of quadratic BSDEs with unbounded solutions and convex duality methods. A key step throughout the analysis is the verification of the finite entropy condition, which not only resolves the uniqueness issue of the unbounded solutions but also ensures the martingale property of the conditional value process and establishes its convex dual representation. Four applications are developed, namely, utility indifference valuation, regime switching, consumption, and Epstein-Zin utility maximization.

There are several possible directions for further research, which we list below and briefly discuss their difficulties and potential resolutions.

(i) When solving the exponential utility maximization problems throughout Chapters \ref{section:main} and \ref{Chapter7}, we imposed the boundedness assumption on the market price of risk $\theta$ while focusing on the unboundedness of the random endowment. Extending to both unbounded market price of risk and unbounded random endowments is actually straightforward. For example, one may impose the exponential integrability condition on the market price of risk, as considered in Theorem \ref{theorem:BSDE_RUD} in Chapter \ref{Chapter8}. We leave such a straightforward extension to the interested reader.   

(ii) Both expected exponential utility and Epstein-Zin utility (with expected power utility as a special case) exhibit the homothetic property for their conditional value processes. Consequently, the optimal trading strategy is independent of wealth. However, this property may not be satisfactory in certain scenarios, prompting the consideration of a general utility function $U$ that is solely increasing and concave, satisfying additional asymptotic conditions (refer to, for instance, Kramkov and Schachermayer \cite{KS1999, KS2003}). In cases where the random endowment $F$
is bounded and the portfolio constraint is a subspace of $\mathbb{R}^d$, the first definitive result was obtained by Horst et al.~\cite{HHIRZ}, who explored the martingale property of the marginal utility and proposed a coupled forward-backward stochastic differential equation (FBSDE) representation of $Y$ such that
$$U_x(X_t^{\pi^{\star}}-Y_t)=E\left[U_x(X_T^{\pi^{\star}}-F)|\mathcal{F}_t\right],$$
for any optimal trading strategy $\pi^{\star}$. See also Liang et al.~\cite{Liang2023} for an infinite horizon version related to forward utility maximization. Under further regularity conditions on $U$, the corresponding FBSDE has been solved by Fromm and Imkeller \cite{Fromm_Imkeller} employing the decoupling field technique. However, the general case where $F$ is unbounded and the portfolio constraint is convex and closed remains unknown to date. 

(iii)  Our asymptotic results on utility indifference valuation and its methodology appear to be new. The central idea is based on the observation that when the portfolio constraint takes the form of a convex and closed cone, a scaling property emerges for the characterizing BSDE governing the conditional value process (see Chapter \ref{section:asymptotic}). However, the more general case of convex and closed set constraints remains open, as the scaling property fails in such a general situation.

(iv) We have made progress on the consumption-investment problem with Epstein-Zin recursive utility by incorporating general convex and closed portfolio constraints and solving the cases for all possible parameter values of relative risk aversions as long as the elasticity of intertemporal substitution is greater than one, which includes the expected power utility case. However, the case where the elasticity of intertemporal substitution is less than one remains open, as the uniqueness of the corresponding characterizing BSDE solution is not clear, and there are challenges in verifying the finite entropy condition.

(v)  Finally, we would also like to point out a forward version of Epstein-Zin recursive utility maximization in the spirit of forward performance processes, a concept introduced by Musiela and Zariphopoulou \cite{Musiela2}. The aim is to construct a forward Epstein-Zin recursive utility $V(t,x)$ as a random field that is consistent over the entire time horizon: For any admissible $(u,c)$,
$$V(t,X_t^{u,c})+\int_0^{t}f^{EZ}(X_s^{u,c}c_s,V(s,X_s^{u,c}))ds,\ \ \ t\geq 0$$
is a supermartingale, and there exists an optimal $(u^{\star},c^{\star})$ such that
$$V(t,X_t^{u^{\star},c^{\star}})+\int_0^{t}f^{EZ}(X_s^{u^{\star},c^{\star}}c_s^{\star},V(s,X_s^{u^{\star},c^{\star}}))ds,\ \ \ t\geq 0$$
is a martingale. As shown in Liang and Zariphopoulou \cite{Liang2017} for the ergdic BSDE representation of the exponential and power utility cases, since the Epstein-Zin recursive utility is also homothetic, it is expected that there exists an ergodic BSDE representation for the  forward Epstein-Zin recursive utility $V(t,x)$.


\backmatter


\begin{thebibliography}{A}



\bibitem{Huang2023} Aurand, J. and Huang, Y. J. (2023). Epstein‐Zin utility maximization on a random horizon. {\it Math. Finance} \textbf{33(4)} 1370--1411.

\bibitem{Barrieu} Barrieu, P. and El Karoui, N. (2013).
Monotone stability of quadratic semimartingales with applications to
unbounded general quadratic BSDEs. {\it Ann. Probab.} \textbf{41} 1831-1863.

\bibitem{Bec0} Becherer, D. (2003). Rational hedging and valuation of
integrated risks under constant absolute risk aversion. {\it
	Insur. Math. Econ.} \textbf{33(1)} 1--28.

\bibitem{Bec} Becherer, D. (2006).
Bounded solutions to backward SDEs with jumps for utility
optimization and indifference hedging. {\it Ann. Appl. Probab.}
\textbf{16} 2027--2054.

\bibitem{Biagini} Biagini, S., Frittelli, M. and Grasselli, M. (2011).
Indifference price with general semimartingales. {\it Math. Finance}
{\bf 21(3)} 423-446.

\bibitem{Bremaud_0}  Br\'emaud, P. (1981). {\it Point processes and queues: Martingale dynamics}. 
{Springer-Verlag, Berlin and Heidelberg}.

\bibitem{Bremaud} Br\'emaud, P. (1999). {\it Markov chains: Gibbs fields, Monte Carlo
simulation, and queues.} {Springer-Verlag, New York}.

\bibitem{bdh} Briand, P., Delyon, B., Hu, Y., Pardoux, E. and
Stoica, L. (2003). $L^p$ solutions of backward stochastic
differential equations. {\it Stochastic Process. Appl.} {\bf  108}
109--129.

\bibitem{Briand} Briand, P. and Elie, R. (2013).
A simple constructive approach to quadratic BSDEs with or without
delay. {\it Stochastic Process. Appl.} \textbf{123} 2921--2939.

\bibitem{BH1} Briand, P. and Hu, Y. (2006).
BSDE with quadratic growth and unbounded terminal value. {\it
Probab. Theory Related Fields} {\bf 136} 604--618.

\bibitem{BH2} Briand, P. and Hu, Y. (2008).
Quadratic BSDEs with convex generators and unbounded terminal
conditions. {\it Probab. Theory Related Fields} {\bf 141} 543--567.

\bibitem{Carmona} Carmona, R. (ed.) (2009). {\it Indifference pricing: theory and applications}.
Princeton University Press, Princeton. 

\bibitem{CH} Cheridito, P. and Hu, Y. (2011). Optimal consumption and investment in incomplete markets with general constraints. {\it Stoch. Dyn.} {\bf 11} 283--299.

\bibitem{CKS} Cvitanic, J., Karatzas, I. and Soner, H. M. (1998).
Backward stochastic differential equations with constraints on the
gains-process. {\it Ann. Probab.} {\bf 26} 1522-1551.

\bibitem{Davis}
Davis, M. H. A. (1997). Option pricing in incomplete markets. {\it In:
Dempster, M.A.H. and Pliska, S.R. (eds), Mathematics of Derivative
Securities. Cambridge University Press, Cambridge}. 216--227

\bibitem{dgr} Delbaen, F., Grandits, P., Rheinlander, T., Samperi, D., Schweizer, M., and Stricker, C. (2002).
Exponential hedging and entropic penalties. {\it Math. Finance} {\bf
12} 99--123.

\bibitem{Delbaen0} Delbaen, F., Hu, Y. and Bao, X. (2011).
Backward SDEs with superquadratic growth. {\it Probab. Theory
Related Fields} \textbf{150} 145--192.

\bibitem{dhr} Delbaen, F., Hu, Y. and Richou, A. (2011).
On the uniqueness of solutions to quadratic BSDEs with convex
generators and unbounded terminal conditions. {\it Ann. inst. Henri Poincare (B) Probab. Stat.} {\bf 47(2)} 559--574.

\bibitem{dhr2} Delbaen, F., Hu, Y. and Richou, A. (2015).
On the uniqueness of solutions to quadratic BSDEs with convex
generators and unbounded terminal conditions: The critical case.
{\it Discrete Contin Dyn Syst Ser A} {\bf
35(11)} 5273--5283.

\bibitem{Meyer1978} Dellacherie, C. and Meyer, P. A. (1978). {\it Probabilities and potential}, North-Holland Mathematics Studies 29. North-Holland Publishing Co.

\bibitem{Matoussi2016} El Karoui, N., Matoussi, A. and Ngoupeyou, A. (2016). Quadratic exponential semimartingales and application to BSDEs with jumps. {\it arXiv preprint arXiv:1603.06191.}

\bibitem{El} El Karoui, N. and Rouge, R. (2000).
Pricing via utility maximization and entropy. {\it Math. Finance}
{\bf 10} 259--276.

\bibitem{EpsteinZin}
Epstein, L. G. and Zin, S. E. (1989). Substitution, risk aversion, and the temporal behavior of consumption and asset returns: A theoretical framework. {\it Econometrica} \textbf{57(4)} 937–-969.

\bibitem{Fan}
Fan, S., Hu, Y., and Tang, S. (2023). A user's guide to 1D nonlinear backward stochastic differential equations with applications and open problems. {\it arXiv preprint arXiv:2309.06233}.

\bibitem{Feng2023}
Feng, Z. and Tian, D. (2023). Optimal consumption and portfolio selection with Epstein–Zin utility under general constraints. {\it Probab. Uncertain. Quant. Risk} \textbf{8(2)} 281--308.

\bibitem{FMS} Frei, C., Malamud, S. and Schweizer, M. (2011).
Convexity bounds for BSDE solutions, with applications to
indifference valuation. {\it Probab. Theory and Related Fields}
\textbf{150(1-2)} 219--255.

\bibitem{FS}
Frei, C. and Schweizer, M. (2009). Exponential utility indifference
valuation in a general semimartingale model. \emph{In: Delbaen, F.,
Rasonyi, M. and Stricker, C (eds), Optimality and Risk--Modern
Trends in Mathematical Finance. The Kabanov Festschrift, Springer}.
49--86.

\bibitem{Fromm_Imkeller}
Fromm, A. and Imkeller, P. (2020). Utility maximization via decoupling fields. {\it Ann. Appl. Probab.} {\bf 30(6)} 2665--2694.

\bibitem{Lin2023} Gu, Z., Lin, Y. and Xu, K. (2023). Quadratic exponential BSDEs driven by a marked point process. {\it arXiv preprint arXiv:2310.14728.}

\bibitem{Henderson1} Henderson, V. (2000). Valuation of claims on nontraded assets using
utility maximization. {\it Math. Finance} {\textbf{12}} 351--373.

\bibitem{Henderson3} Henderson, V. and Liang, G. (2014).
Pseudo linear pricing rule for utility indifference valuation.
{\it  Finance Stoch.} {\bf 18(3)} 593--615.

\bibitem{Herdegen2023(1)} 
Herdegen, M., Hobson, D. and Jerome, J. (2023). The infinite-horizon investment–consumption problem for Epstein–Zin stochastic differential utility. I: Foundations. {\it Finance Stoch.} \textbf{27(1)}, 127--158.

\bibitem{Herdegen2023(2)}
Herdegen, M., Hobson, D. and Jerome, J. (2023). The infinite-horizon investment–consumption problem for Epstein–Zin stochastic differential utility. II: Existence, uniqueness and verification for $\nu\in(0, 1)$. {\it Finance Stoch.} \textbf{27(1)}, 159--188.

\bibitem{Heunis}
Heunis, A. J. (2015). Utility maximization in a regime switching model with convex portfolio constraints and margin requirements: optimality relations and explicit solutions. {\it SIAM J Control Optim.} {\bf 53(4)}, 2608--2656.

\bibitem{HodgesNeuberger} Hodges, S. and Neuberger, A. (1989). Optimal replication of contingent claims under transactions
costs. {\it Rev. Futures Mark.} \textbf{8} 222--239.

\bibitem {HHIRZ} Horst, U., Hu, Y., Imkeller, P., Reveillac, A. and Zhang, J. (2014).
Forward-backward systems for expected utility maximization, {\it Stochastic
Process. Appl.} \textbf{124(5)} 1813--1848.

\bibitem{him} Hu, Y., Imkeller P. and Muller, M. (2005).
Utility maximization in incomplete markets. {\it Ann. Appl. Probab.}
{\bf 15} 1691--1712.

\bibitem{HuPeng}
Hu, Y. and Peng, S. (2006). On the comparison theorem for
multidimensional BSDEs. {\it Comptes Rendus Mathematique} {\bf
343(2)} 135--140.

\bibitem{HLT}
Hu, Y., Liang, G. and Tang, S. (2020). Systems of ergodic BSDEs arising in regime switching forward performance processes. {\it SIAM J Control Optim.} \textbf{58(4)} 2503--2534.

\bibitem{HSX2022}
Hu, Y., Shi, X. and Xu, Z. Q. (2022). Optimal consumption-investment with coupled constraints on consumption and investment strategies in a regime switching market with random coefficients. {\it arXiv preprint arXiv:2211.05291.}

\bibitem{KLS1987}
Karatzas, I., Lehoczky, J. P. and Shreve, S. E. (1987). Optimal portfolio and consumption decisions for a small investor on a finite horizon. {\it SIAM J. Control Optim.} {\bf 25(6)} 1557--1586.

\bibitem{KLSX1991}
Karatzas, I., Lehoczky, J. P., Shreve, S. E. and Xu, G. L. (1991). Martingale and duality methods for utility maximization in an incomplete market. {\it SIAM J. Control Optim.} {\bf 29(3)} 702--730.

\bibitem{KS} Karatzas, I. and Shreve, S. (1998). {\it Brownian motion and stochastic
calculus}. Springer-Verlag, New York. 

\bibitem{Kaz} Kazamaki, N. (1994). {\it Continuous exponential martingales and BMO.} Lecture Notes
in Math. 1579. Springer-Verlag, Berlin Heidelberg. 

\bibitem{Kobylanski} Kobylanski, M. (2000).
{Backward stochastic differential equations and partial differential
equations with quadratic growth}. {\it Ann. Probab.} {\bf 28}
558--602.

\bibitem{Kraft2017} Kraft, H., Seiferling, T. and Seifried, F. T. (2017). Optimal consumption and investment with Epstein–Zin recursive utility. {\it Finance Stoch.} \textbf{21} 187--226.

\bibitem{KS1999}
Kramkov, D. and Schachermayer, W. (1999). The asymptotic elasticity of utility functions and optimal investment in incomplete markets. {\it Ann. Appl. Probab.}  {\bf 9}
904--950.

\bibitem{KS2003}
Kramkov, D. and Schachermayer, W. (2003). Necessary and sufficient conditions in the problem of optimal investment in incomplete markets. {\it Ann. Appl. Probab.} {\bf 13} 1504--1516.


\bibitem{Liang2023}
Liang, G., Sun, Y. and Zariphopoulou, T. (2023). Representation of forward performance criteria with random endowment via FBSDE and application to forward optimized certainty equivalent. {\it arXiv preprint arXiv:2401.00103.}

\bibitem{Liang2017}
Liang, G. and Zariphopoulou, T. (2017). Representation of homothetic forward performance processes in stochastic factor models via ergodic and infinite horizon BSDE. {\it SIAM J. Financ. Math.}, {\bf 8(1)} 344--372.

\bibitem{Lepeltier1997}
Lepeltier, J. P. and San Martin, J. (1997). Backward stochastic differential equations with continuous coefficient. {\it Statistics \& Probability Letters}, {\bf 32(4)}, 425--430.

\bibitem{MS} Mania M. and Schweizer, M. (2005).
{Dynamic exponential utility indifference valuation}. {\it Ann.
Appl. Probab.} {\bf 15} 2113--2143.


\bibitem{Matoussi2019} Matoussi, A. and Salhi, R. (2019). Exponential quadratic bsdes with infinite activity jumps. {\it arXiv preprint, arXiv:1904.08666.}

\bibitem{Matoussi2018}
Matoussi, A. and Xing, H. (2018). Convex duality for Epstein–Zin stochastic differential utility. {\it  Math. Finance} \textbf{28(4)} 991--1019.


\bibitem{Mocha}
Mocha, M. and Westray, N. (2012). Quadratic semimartingale BSDEs
under an exponential moments condition. {\it S\'eminaire de
{P}robabilit\'es {XLIV}, Lecture Notes in Math., Springer}. 105--139.

\bibitem{Morlais} Morlais, M. A. (2009).
Quadratic BSDEs driven by a continuous martingale and applications
to the utility maximization problem. {\it Finance Stoch.}
\textbf{13(1)} 121--150.

\bibitem{Musiela} Musiela, M. and Zariphopoulou, T. (2004). An example of
indifference prices under exponential preferences. {\it Finance Stoch.} {\textbf{8}} 229--239.

\bibitem{Musiela2}
Musiela, M. and Zariphopoulou, T. (2008) Optimal asset allocation under forward exponential performance criteria, {\it Markov Processes and Related Topics: A Festschrift for T. G. Kurtz, IMS Lecture Notes-Monograph Series 4}. 285--300.

\bibitem{Gordan} Owen, M. P. and Zitkovic, G. (2009).
Optimal investment with an unbounded random endowment and
utility-based pricing. {\it Math. Finance} {\bf 19(1)}
129-159.


\bibitem{Peng} Peng, S. (1999).
Monotonic limit theorem of BSDE and nonlinear decomposition theorem
of Doob-Meyers type. {\it Probab. Theory and Related Fields} {\bf
113(4)} 473-499.

\bibitem{Sethi} Sethi, S. P. (2012). {\it Optimal consumption and investment with bankruptcy}. Springer-Verlag, New York.


\bibitem{Tang2015} Tang, S. (2015), 
Dynamic programming for general linear quadratic optimal stochastic control with random coefficients.{\it 
SIAM J. Control Optim.} {\bf53(2)} 1082-1106.

\bibitem{Tevzadze}
Tevzadze, R. (2008). {Solvability of backward stochastic
differential equations with quadratic growth}. {\it Stochastic
Process. Appl.} \textbf{118} 503--515.

\bibitem{Xing2017}
Xing, H. (2017). {Consumption–investment optimization with Epstein–Zin utility in incomplete markets}. {\it Finance Stoch.} \textbf{21} 227--262.

\bibitem{Yin}
Yin, G. G. and Zhu, C. (2009). {\it Hybrid switching diffusions: properties and applications}. Springer-Verlag, New York.

\bibitem{ZengTaksar2013}
Zeng, X. and Taksar, M., (2013). {\it A stochastic volatility model and optimal portfolio selection}. {\it Quant. Finance}, 13(10), pp.1547-1558.


\bibitem{ZhangJF}
Zhang, J. (2017). {\it Backward stochastic differential equations: From linear to fully nonlinear theory}. Springer-Verlag, New York.

\end{thebibliography}


\bibliographystyle{amsalpha}


%
%
%

\begin{theindex}

\item Admissible strategies with constraints 
\subitem ${\mathcal A}_D$, 14
\subitem ${\mathcal A}_D^{conv}$, 16
\subitem $\mathcal{A}_D^{M}$, 31
\subitem $\mathcal{A}_D^{cone}$, 39
\subitem $\mathcal{A}_D^{\mathbb{G}}$, 53
\subitem $\mathcal{A}_D^C$, 65
\subitem $\mathcal{A}_D^{EZ}$, 80
\item Asymptotics for the risk aversion parameter, 39

\indexspace

\item B-D-G inequality, 25, 33

\indexspace

\item Class (D) condition, 14, 16, 28, 53, 65, 81
\item Cone condition, 39
\item Consumption, 65
\item Convex dual representation, 36, 61, 71, 92

\indexspace

\item De la Vall\'ee-Poussin theorem, 28, 36, 60, 61, 70, 71, 90, 92   
\item Doob's inequality 23, 85

\indexspace

\item Fenchel inequality, 29, 60, 70, 91, 
\item Fenchel-Moreau theorem, 18, 37, 63, 72, 93 
\item Finite entropy condition, 28, 36, 60, 61, 70, 71, 90, 92

\indexspace

\item Inf convolution, 20, 56, 77, 84
\item It\^o-Tanaka  formula, 23, 24

\indexspace

\item Lipschitz BSDE
\subitem comparison theorem, 78, 79, 83, 84

\indexspace

\item Localization, 22, 43, 58, 67, 78, 80, 84, 85

\indexspace

\item Market price of risk, 13, 53, 76
\item Martingale optimality principle, 16, 27, 58, 68, 86
\item Martingale measure
\subitem minimal entropy martingale measure, 51
\subitem minimal local martingale measure, 31, 76
\subitem minimal martingale measure, 31
\item Minimal entropy representation, 49
\item Monotonic limit theorem, 44 
\item Multidimensional comparison theorem, 56
\indexspace

\item Quadratic BSDE 
\subitem comparison theorem, 21, 42, 46, 67, 73
\subitem stability property, 20
\subitem exponential quadratic BSDE, 55
\subitem with bounded terminal data, 20
\subitem with unbounded coefficients and generator monotonic in $y$, 81
\subitem with unbounded terminal data, 21
\subitem with unbounded terminal data and generator linear in $y$, 66
\subitem with unbounded terminal data in multi dimension, 54
\indexspace

\item Regime-switching market, 53

\indexspace

\item Scaling property, 39
\item Stochastic exponential, 19, 31, 36, 61, 71, 76, 89, 92
\item Subdifferential, 18
\item Super-replication price, 51

\indexspace

\item Truncation, 20, 57, 84 

\indexspace

\item Utility
\subitem Epstein-Zin recursive utility, 75
\subitem Exponential utility, 14
\subitem Power utility, 75
\subitem Utility indifference valuation, 35

\indexspace

\item Value function, 14
\item Conditional value process, 16

\indexspace

\end{theindex}



\end{document}